\documentclass[a4paper,10pt]{article} 
\usepackage{amsmath,amssymb,amsthm,amsfonts,epsfig,color,url}
\usepackage{latexsym,mathrsfs,cite} 
\usepackage{hyperref} 
\usepackage{graphicx}
\usepackage{float}
\usepackage[font=small]{caption,subfig}
\usepackage[margin=3cm]{geometry}

\graphicspath{{figures/}}

\newtheorem{theorem}{Theorem}[section]
\newtheorem{lemma}[theorem]{Lemma}
\newtheorem{corollary}[theorem]{Corollary}

\newtheorem{definition}[theorem]{Definition}

\newtheorem{hypothesis}[theorem]{Hypothesis}

\newtheorem{remark}[theorem]{Remark}

\numberwithin{equation}{section}
\numberwithin{figure}{section}

\def\discr{{\delta_{\calM}}}
\def\discrh{{\delta_{\calH}}}

\def\R{\mathbb{R}}
\def\C{\mathbb{C}}

\def\Cspace{{\sf C}}
\def\Lspace{{\sf L}}

\def\crit{\mathrm{crit}}

\def\E{E}
\def\rms{{\rm s}}

\def\kap{{\tilde\kappa}}

\def\tel{{\tilde\ell}}

\def\ep{{\rm ep}}
\def\tp{{\rm tp}}
\def\mp{{\rm mp}}
\def\ex{{\rm ex}}
\def\nl{{\rm nl}}
\def\crit{{\rm c}}

\def\scal{\varepsilon}

\def\hK{k_0}
\def\hq{q_2}
\def\tq{q_0}
\def\tA{{\tilde A}}

\def\rmi{\mathrm{i}}
\def\rmd{\mathrm{d}}
\def\rme{\mathrm{e}}

\def\calO{\mathcal{O}}
\def\calL{\mathcal{L}}

\def\calE{\mathcal{E}}
\def\calB{\mathcal{B}}
\def\calH{\mathcal{H}}
\def\calM{\mathcal{M}}
\def\calQ{\mathcal{Q}}

\def\hex{\mathrm{hex}}

\def\Id{{\rm Id}}
\def\A{L}
\def\Lc{\widetilde{\mathcal{L}}}
\def\M{M}
\def\B{B}
\def\Q{Q}
\def\K{K}

\def\x{{\bf x}}
\def\kc{{{\bf k}_{\rm c}}}
\def\kcsq{{{\bf k}_{\rm c}^2}}
\def\ta{{b}}

\def\tQ{\Q_0}
\def\hQ{\Q_2}
\def\k{{\bf k}}
\def\krec{{\bf k}^{\rm sq}}
\def\bbeta{{\beta\beta}}

\def \talpha {{\tilde\alpha}}
\def \calpha {{\check\alpha}}

\def\Re{\mathrm{Re}}

\def \sgn {\mathrm{sgn}}

\def \diag {\mathrm{diag}}

\def\rhomb{{\rm qh}}
\def\Rhomb{{\rm qH}}
\def\rhb{{\rm rh}}
\def\rect{{\rm rect}}
\def\sq{{\rm sq}}
\def\zz{{\rm zz}}
\def\eh{{\rm eh}}

\begin{document}
\title{The impact of advection on the stability of stripes on lattices near planar Turing instabilities}
\author{Jichen Yang\footnote{University of Bremen, Faculty 3 -- Mathematics, Bibliothekstrasse 5, 28359 Bremen, Germany  ({\tt jyang@uni-bremen.de})}
\and Jens D. M. Rademacher\footnote{University of Bremen, Faculty 3 -- Mathematics, Bibliothekstrasse 5, 28359 Bremen, Germany  ({\tt jdmr@uni-bremen.de})}
\and Eric Siero\footnote{Carl von Ossietzky University of Oldenburg, Institute for Mathematics, Carl von Ossietzky Str.9-11, 26111 Oldenburg, Germany  ({\tt eric.siero@uni-oldenburg.de})}}
\date{February 27, 2020}

\maketitle

\begin{abstract}
Striped patterns are known to bifurcate in reaction-diffusion systems with differential isotropic diffusions at a supercritical Turing instability. In this paper we study the impact of weak anisotropy by directional advection on the stability of stripes with respect to various lattice modes, and the role of quadratic terms therein. We focus on the generic form of planar reaction-diffusion systems with two components near such a bifurcation. Using centre manifold reduction we derive a rigorous parameter expansion for the critical eigenvalues for lattice mode perturbations, specifically nearly square and nearly hexagonal ones. This provides detailed formulae for the loci of stability boundaries under the influence of the advection and quadratic terms. In particular, the well known destabilising effect of quadratic terms can be counterbalanced by advection, which leads to intriguing arrangements of stability boundaries. We illustrate these results numerically by a specific example. Finally, we show numerical computations of these stability boundaries in the extended Klausmeier model for vegetation patterns and show stripes bifurcate stably in the presence of sufficiently strong advection.
\end{abstract}

\section{Introduction}\label{s:intro} 
 
The ubiquitous isotropic pattern forming Turing instabilities are known to generate various solutions, that are dominated in one-dimensional space by spatially periodic solutions. These trivially extend to striped solutions in two-dimensional space, where they are in competition with hexagonal and square shaped states, e.g., \cite{Hoyle2006}. This naturally leads to the question which pattern is selected at onset of the Turing instability -- we consider supercritical Turing instabilities only. In the isotropic situation and for generic quadratic terms, it is well known that stripes are unstable with respect to modes on the hexagonal lattice, as discussed in \cite{Gowda2014} in the context of vegetation patterns. However, in \cite{Siero2015} it was found that striped vegetation patterns in a sloped terrain are stable at onset and are connected to large amplitude stripes within the `Busse balloon' of stable stripes in wavenumber-parameter space. Here the slope is modelled by an advective term in the water component, which breaks the spatial isotropy. From a symmetry perspective for weakly anisotropic perturbations this has been predicted already in \cite{Callahan2000} and the destabilising effect of advection terms on homogeneous steady states have been broadly studied in the context of differential flows, e.g., \cite{Rovinsky1992,Merkin2000,Carballido-Landeira2012} and also appear in ecology, e.g., \cite{Wang2009,Cangelosi2015,Bennett2019}.

\medskip
In this paper we show that advection can have a stabilising effect on relevant lattice modes for stripes, counteracting in particular the destabilising effect of quadratic terms, and explaining the observations in \cite{Siero2015}. This study extends our analysis in \cite{Yang2019a} by considering lattice modes or equivalently stability on certain rectangular domains with periodic boundary conditions. Specifically, in \cite{Yang2019a} a detailed expansion for the bifurcation of stripes and the stability against large-wavelength modes, also called sideband modes, was studied. The consideration of periodic domains adapted to suitable wavevectors is a classical theme in amplitude equations, and is a standard tool in the context of Turing instabilities, see \cite{Hoyle2006,Siero2015,Pena2003} and the references therein. However, the analysis of weak anisotropy seems scarce.

\medskip
In our analysis we employ centre manifold reductions on domains that are nearly square and nearly `hexagonal', i.e., with the hexagonal lattice for wavevectors. We expand the critical eigenvalues of the stripes in perturbed spatial scalings as well as the system parameters. Herein we can conveniently use the existence of stripes from \cite{Yang2019a}. The advantages of this approach are that it is fully rigorous and that we gain direct access to all relevant characteristic quantities in terms of the advection, the quadratic terms, stretching and compressing.  A particular motivation is to bridge the discussion of stripe stability in \cite{Siero2015} for a variant of the Klausmeier model with rather large advection to the results from \cite{Gowda2014} for zero advection.

The approach applies to arbitrary number of components, but the parameter spaces and determination of signs of relevant characteristics become analytically less accessible for more than two components. Hence we restrict our attention to this case.

\medskip
Upon changing coordinates, the generic form of such a system up to cubic nonlinearity reads
\begin{align}\label{e:RDS}
u_t = D\Delta u + \A u + \calpha \M u + \beta \B u_x + \Q[u,u] + \K[u,u,u], \; \x\in\R^2
\end{align}
with multilinear functions $\Q, \K$ and diagonal diffusion matrix $D>0$; higher order nonlinear terms can be added without change to our results near bifurcation. We assume that for $\calpha=\beta=0$ the zero steady state is at a Turing instability with wavenumber $\kc$, cf.\ Definition~\ref{def:Turing} below, and that $\calpha$ moves the spectrum through the origin. The isotropy is broken for $\beta\neq 0$ and we assume, without loss of generality, differential advection 
\[
\B=\B(c)=\begin{pmatrix}1+c & 0\\ 0 & c\end{pmatrix}, \; c\in\R.
\]
Note that $\beta c\partial_x$ appears in both equations as a comoving frame in the $x$-direction, and positive (negative) $\beta$ implies the advection of the first component in negative (positive) $x$-direction.

\medskip
Our main results may be summarised as follows. Here the parameters are $\mu=(\alpha,\beta,\kap,\tel)$, where $\alpha=\lambda_M\calpha$, for certain $\lambda_M\neq0$ determined in \S\ref{s:Turing}, and $\kap=\kappa-\kc$ is the deviation of the stripe's nonlinear wavenumber from $\kc$, i.e., the stripe's spatial period is $2\pi/\kappa$. Lastly, $\tel$ is the deviation of the domain's spatial extent from a square or `hexagonal' domain along the stripe, i.e., in $y$-direction. The velocity parameter $c$ is determined from $\mu$. Throughout we consider $|\mu|\ll 1$, and consider stripes  $U_\rms(x;\mu)$ that are constant in $y$ with amplitude parameter $A=\|\widehat{U_\rms}(\kc;\mu)\|$ the norm of the first Fourier mode. In order to simplify the discussion, we assume the scaling relation 
\begin{align}\label{e:scalingintro}
(A,\alpha, \beta, \kap, \tel) = (\scal A', \scal^2 \alpha', \scal \beta', \scal \kap' ,\scal \tel'), 
\end{align}
with a scaling parameter $\scal> 0$. This scaling is homogeneous for $\mu$ with respect to the relevant terms in the expansion of stripes, cf.\ Theorem~\ref{t:bif} from \cite{Yang2019a} below.

As is well known from the isotropic case, the quadratic terms enter at lower order into stability on the hexagonal lattice and thus should be small in order to discuss changes of stability. A convenient, though not necessary, implementation of this is the following uniform smallness hypothesis that we shall adopt.
 
\begin{hypothesis}\label{h:Qscale}
$\Q[\cdot,\cdot]=\scal\Q'[\cdot,\cdot]$.
\end{hypothesis}
In the standard amplitude/modulation equation approach this assumption is required a priori, while in our approach it enters only a posteriori in order to obtain non-trivial stability boundaries.

\medskip
As mentioned, in this paper we are concerned with finite wavenumber stability and instability; large-wavelength in/stability, i.e., $\kappa,\ell\approx 0$, was analysed in \cite{Yang2019a}. It was shown in~\cite[Corollary 2.6]{Yang2019a} that in the anisotropic case $\beta\neq 0$ stripes are spectrally stable near bifurcation if they are stable against large-wavelength modes, which is reflected in the results of this paper as well, cf. Remark~\ref{r:1Dstab}. It is natural to consider domains whose Fourier wavevectors form periodic lattices and the symmetric ones are square (rotation by $\pi/2$) and hexagonal (rotation by $\pi/3$). We refer to the lattice modes considered on the (nearly) square and (nearly) `hexagonal' domains as the {\em (quasi-)square} and {\em (quasi-)hexagonal modes}, respectively. It turns out that certain quasi-hexagonal modes are more unstable than others, and therefore the resulting stability boundaries are briefly illustrated next. In order to build the foundation for this case, in the body of the paper we begin by discussing quasi-square modes in \S\ref{s:square} and exact hexagonal modes in \S\ref{s:hex}.

\begin{figure}[!t]
\centering
\subfloat[$0\leq|\beta|<\beta_\ep$]{\includegraphics[width=0.35\linewidth]{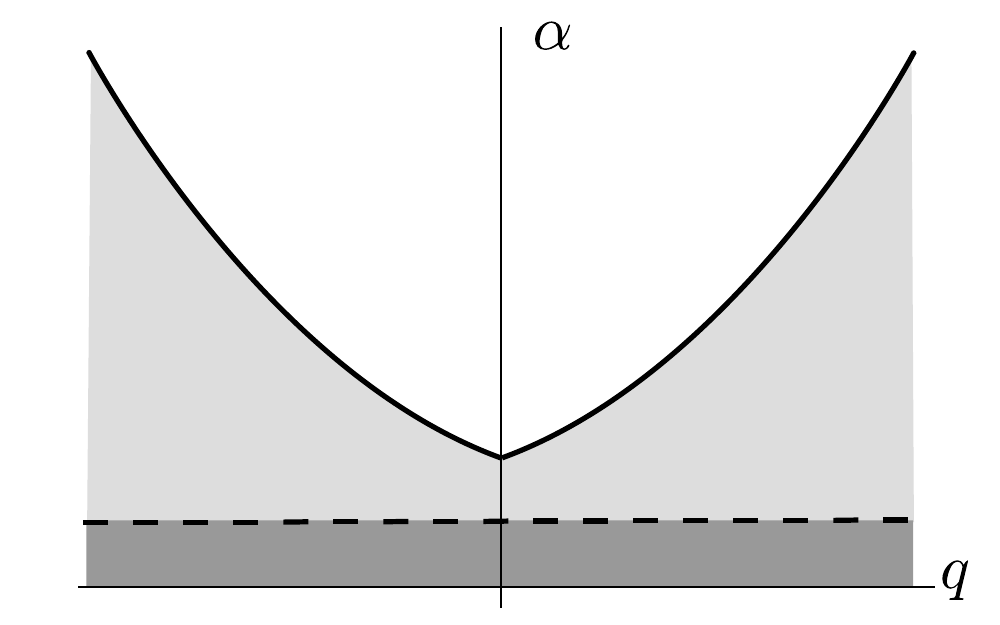}\label{f:rhombetn0kapn0le}}
\hfil
\subfloat[$|\beta|=\beta_\ep$]{\includegraphics[width=0.35\linewidth]{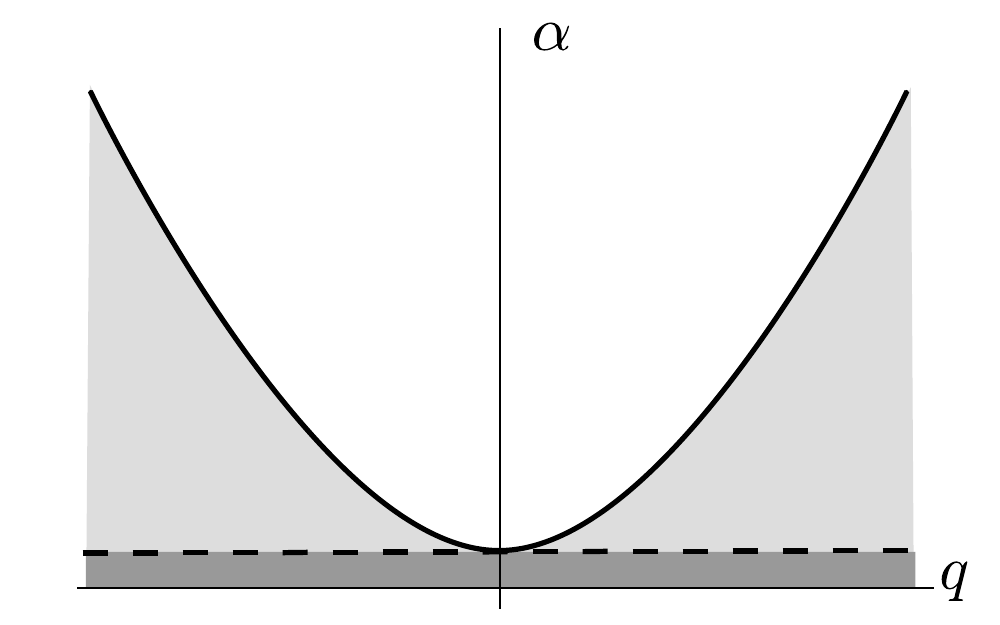}\label{f:rhombetn0kapn0eq}}
\hfil
\subfloat[$|\beta|>\beta_\ep$]{\includegraphics[width=0.35\linewidth]{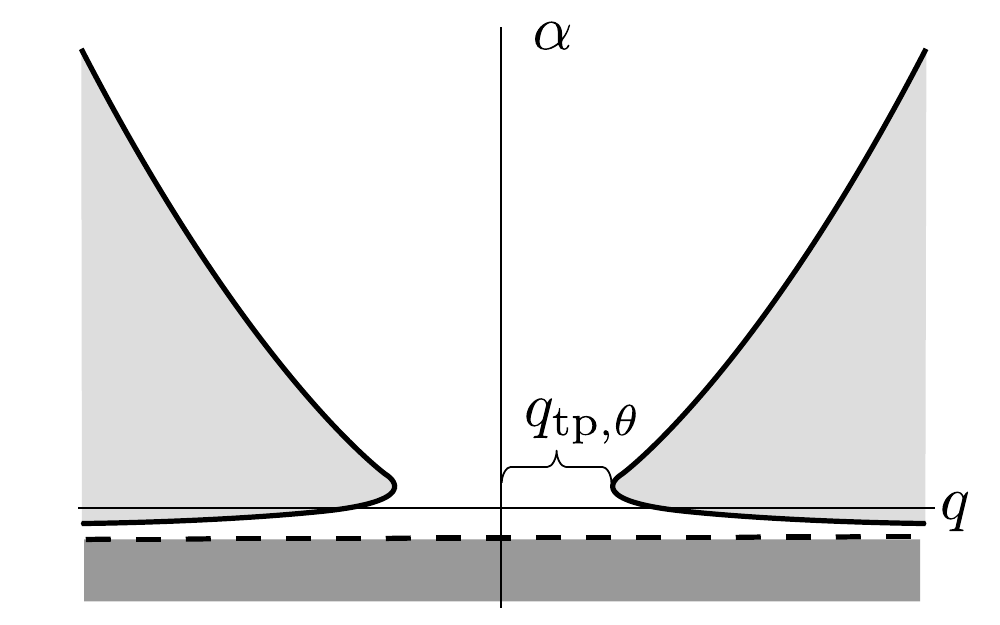}\label{f:rhombetn0kapn0gr}}
\hfil
\subfloat[]{\includegraphics[width=0.35\linewidth]{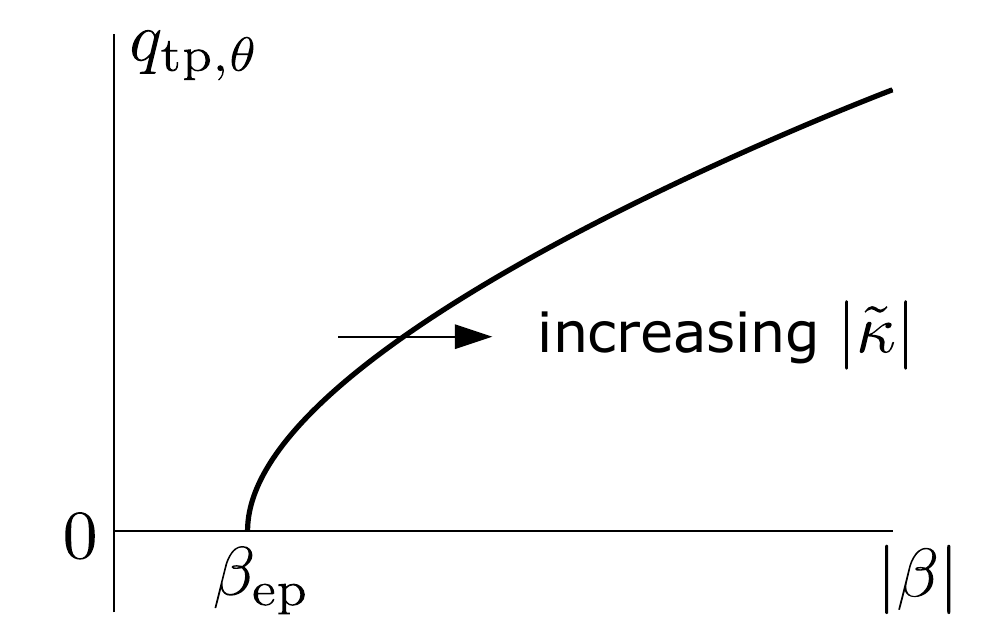}\label{f:betaq}}
	\caption{In (a)--(c) we plot sketches of the quasi-hexagonal in/stability regions in the $(q,\alpha)$-plane for fixed $\kap\neq0$ and $\theta\in(0,1]$, the quadratic term $q=q(\Q)=\calO(\scal)$ measures the effect of the quadratic nonlinearity. Stripes exist in the complement of the dark grey regions; light grey: quasi-hex-unstable; white: quasi-hex-stable. Stripe bifurcation~\eqref{e:bifbnd} (dashed), quasi-hexagonal boundaries~\eqref{e:rhombbnd1} (solid). In (d) we sketch the half width $q_{\tp,\theta}$ of the stable connection marked in (c), cf.~\eqref{e:qtpthe}. The curve intersects the $\beta$-axis at $\beta_\ep$ which is linearly increasing with $|\kap|$. The stripes are quasi-hex-stable below the curve for any $\alpha$. This highlights that for larger $|\kap|$ the connection of stable regions requires larger $|\beta|$, and that the connection is wider for larger $|\beta|$. Note that $\beta_\ep=0$ at $\kap=0$.}
\label{f:rhombetn0kapn0}
\end{figure}

\begin{figure}[!t]
\centering
\subfloat[$\beta=0$]{\includegraphics[width=0.33\linewidth]{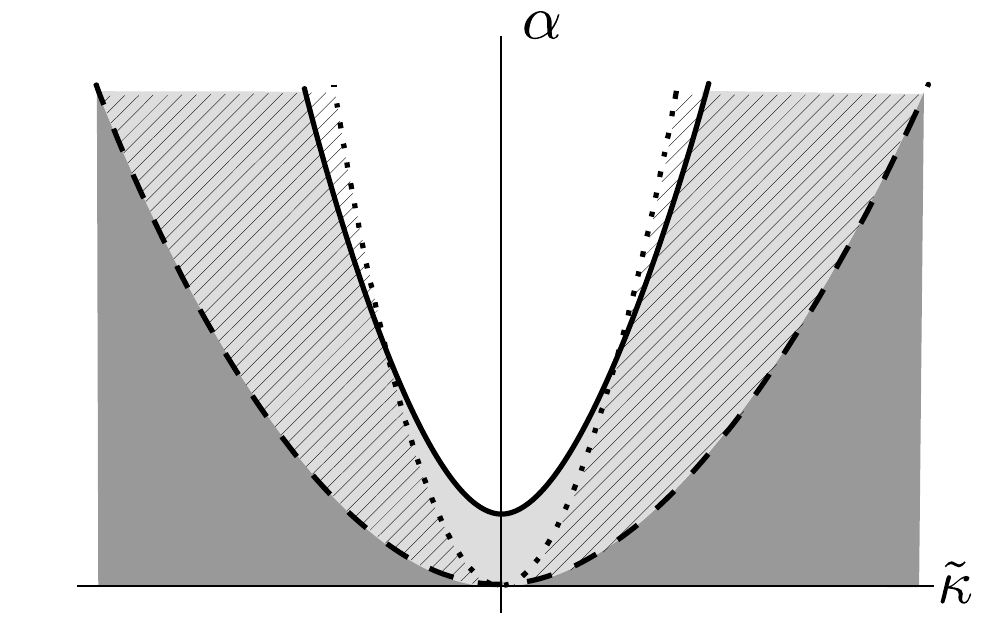}\label{f:rhombkapalp1}}
\hfil
\subfloat[$|\beta|<\beta_\tp$]{\includegraphics[width=0.33\linewidth]{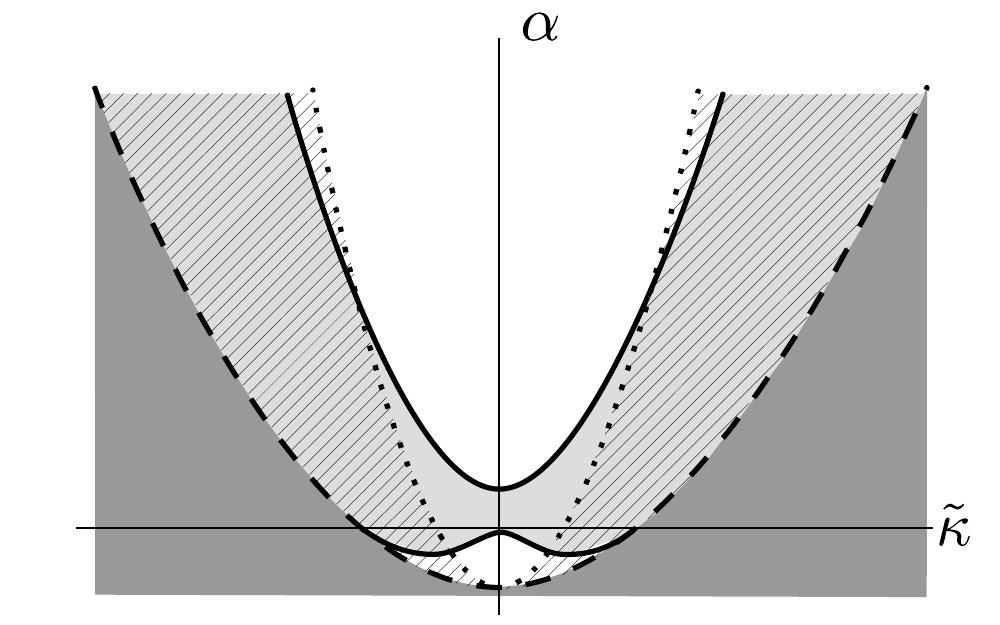}\label{f:rhombkapalp2}}
\hfil
\subfloat[$|\beta|=\beta_\tp$]{\includegraphics[width=0.33\linewidth]{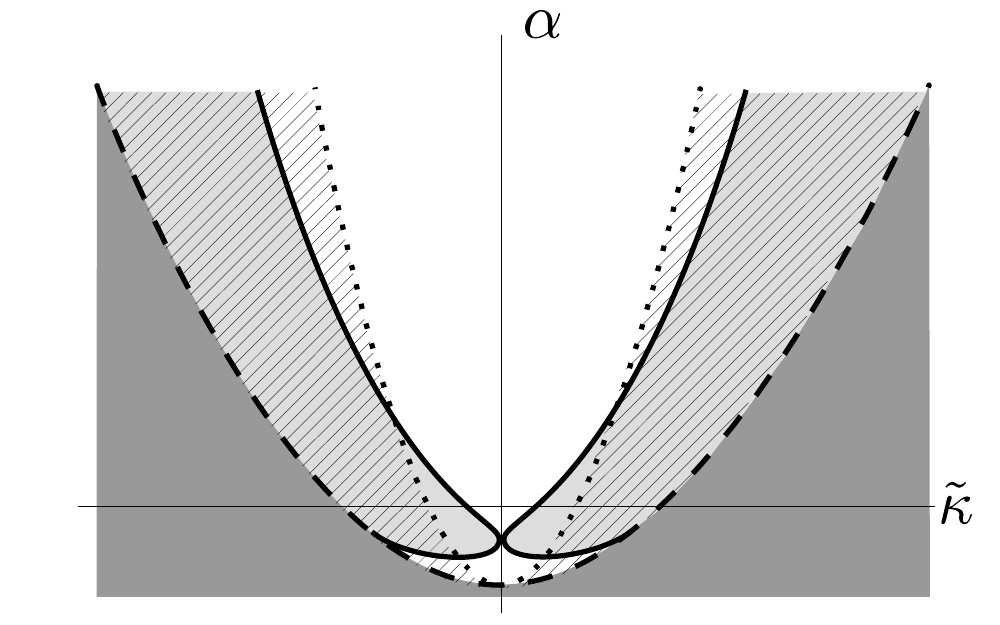}\label{f:rhombkapalp3}}
\hfil
\subfloat[$\beta_\tp<|\beta|<\beta_\ex$]{\includegraphics[width=0.33\linewidth]{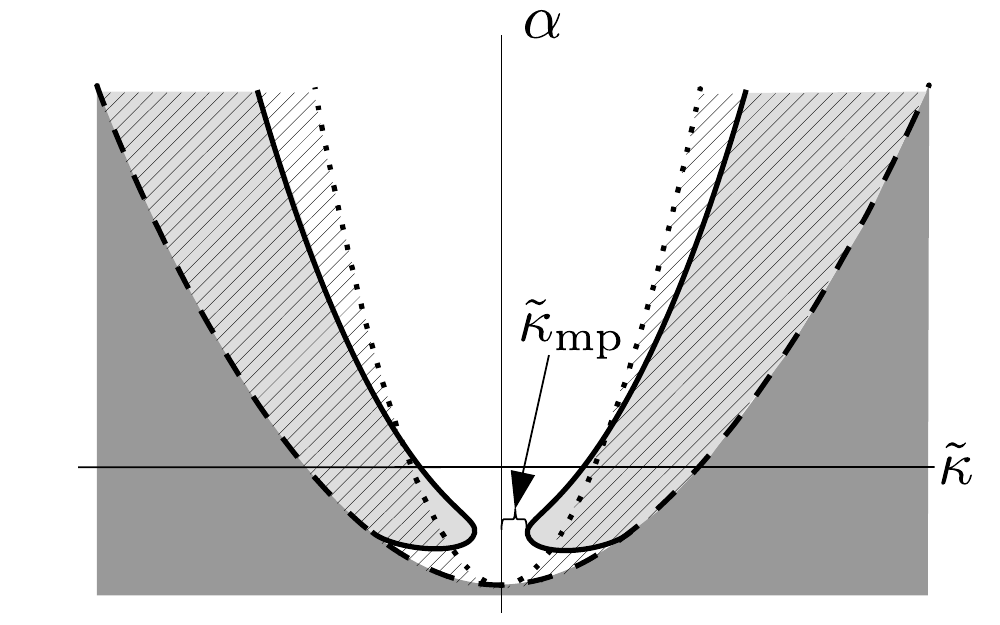}\label{f:rhombkapalp4}}
\hfil
\subfloat[$|\beta|\geq\beta_\ex$]{\includegraphics[width=0.33\linewidth]{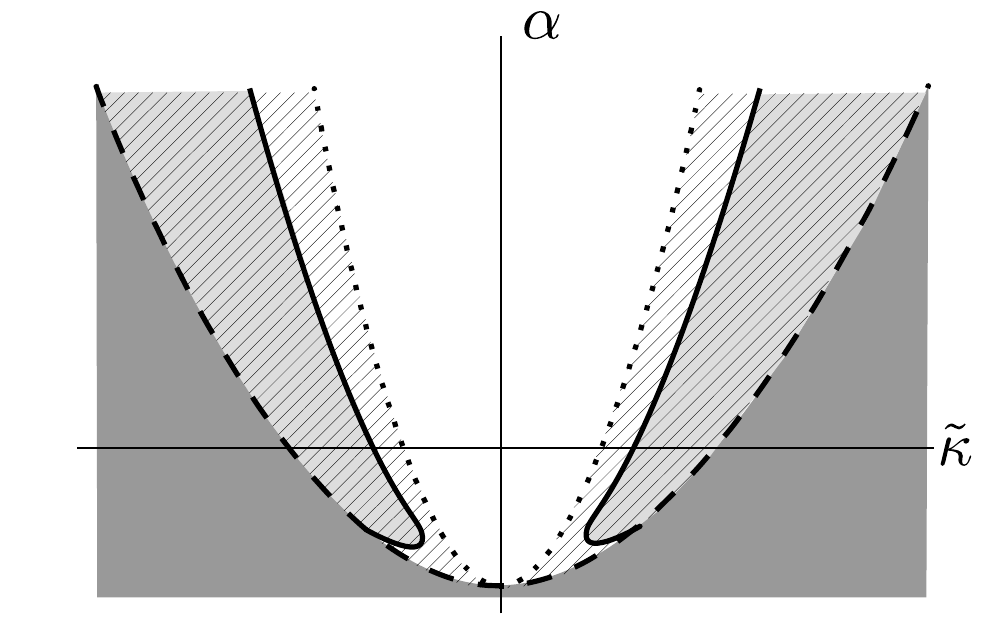}\label{f:rhombkapalp5}}
\hfil
\subfloat[]{\includegraphics[width=0.33\linewidth]{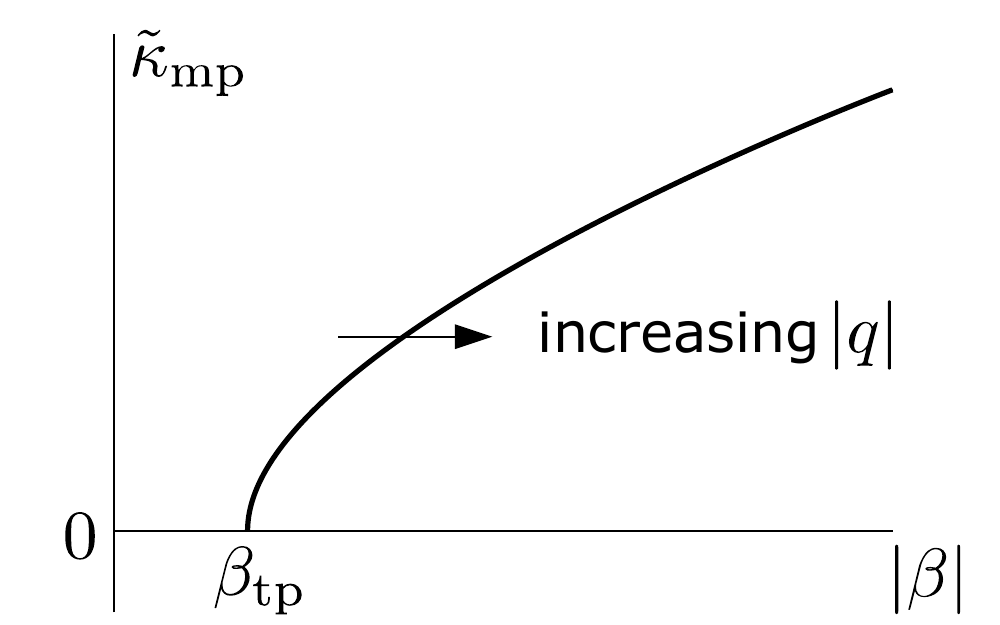}\label{f:betakap}}
	\caption{In (a)--(e) we plot sketches of the quasi-hexagonal stability regions in the $(\kap,\alpha)$-plane for fixed quadratic coefficient $q\neq0$, $q=q(\Q)=\calO(\scal)$ and fixed $\theta\in(0,1]$. Stripes exist in the complement of the dark grey regions. Light grey: quasi-hex-unstable; hatched regions: Eckhaus-unstable; stripes are zigzag-unstable for $\kap<0$. We plot stripe bifurcation~\eqref{e:bifbnd} (dashed), Eckhaus boundaries~\eqref{e:eckhaus} (dotted), quasi-hexagonal boundaries~\eqref{e:rhombbnd1} (solid). In (f) we sketch the half width $\kap_\mathrm{mp}$ \eqref{e:kapmp} of the stable connection marked in (d). The curve intersects the $\beta$-axis at $\beta_\tp$ which is linearly increasing with $|q|$. The stripes are quasi-hex-stable below the curve for any $\alpha$. This illustrates that the stable regions connect later for larger $|q|$ and the connection is wider for larger $|\beta|$.}
\label{f:rhombkapalp}
\end{figure}

\paragraph{Quasi-hexagonal stability boundaries.} 
We consider periodic boundary conditions on the rectangular domains $\x\in\Omega_\rhomb:=[0,4\pi/\kappa]\times[0,4\pi/(\sqrt 3 \ell)]$, $\kappa:=\kc+\kap$, $\ell:=\kc+\tel$, $\tel\neq\kap$. We prove that the ratio $\kap/\tel=-3$ yields the most unstable modes near onset -- it is the scale ratio on which the hexagonal modes of the homogeneous steady state are critical. For generic quadratic term, in the isotropic case $\beta=0$ the stripes are unstable near the onset of Turing instability (Fig.~\ref{f:rhombetn0kapn0le}~\&~\ref{f:rhombkapalp1}). In the anisotropic case, $\beta\neq0$, any advection strength stabilises the stripes with wavenumbers close to the Turing critical wavenumber (Fig.~\ref{f:rhombkapalp2}), but this `small' stability region is not connected to the stable region of larger amplitude stripes. However, the size of the small stability region increases with advection strength and eventually connects to that of larger amplitude stripes (Fig.~\ref{f:rhombetn0kapn0gr}~\&~\ref{f:rhombkapalp4}). Notably, the thresholds are of the form $\beta_\ep= c_\ep |\kap|$, $\beta_\tp =c_\tp |q|$ with explicit constants $c_\ep, c_\tp>0$, cf.\ Fig.~\ref{f:betaq} and \ref{f:betakap}, respectively.  This transitioning from two disconnected stability regions of small and larger amplitude stripes to a connected `Busse balloon' explains the consistency with the results from \cite{Siero2015}: there it was observed that, disregarding zigzag-instability, under increasing advection strength the 2D stability region transitions to an effective 1D stability region determined by the Eckhaus boundary. Recall that the Eckhaus instability is a large-wavelength instability orthogonal to the stripe and the dominant instability mechanism for wavetrains in 1D; we do not further discuss the leading order zigzag instability region $\{\kap<0\}$.


More specifically, the quasi-hexagonal instability compares with the Eckhaus instability as follows. In the presence of a generic quadratic term, the quasi-hexagonal instability is dominant near the onset in the isotropic case (Fig.~\ref{f:rhombkapalp1}) while the Eckhaus instability is dominant near the onset in the anisotropic case (Fig.~\ref{f:rhombkapalp2} to \ref{f:rhombkapalp5}). In particular, the Eckhaus instability is completely dominant for relatively strong advection $\beta>\beta_\ex = c_\ex |q|$ for an explicit constant $c_\ex>0$ (Fig.~\ref{f:rhombkapalp5}).

\medskip
In our analysis, we consider the leading order bifurcation and stability boundaries of stripes with the scaling relation \eqref{e:scalingintro}. This leads to the reflection symmetric bifurcation and stability boundaries, cf. Fig.~\ref{f:rhombkapalp}. Relaxing the scaling relation and including the higher order terms will generically break such symmetry as shown for the zigzag stability boundary in~\cite{Yang2019a}. Indeed, we observe asymmetry of the stability diagrams in the numerical computations of Klausmeier model in \S\ref{s:Klausmeier}. These latter results refine and complete the study of \cite{Siero2015} for small advection and explain how, in a concrete model, the Busse balloon for stripes can connect to a Turing-Hopf point.

\medskip
This paper is organised as follows: In \S\ref{s:pre} we recall the results from \cite{Yang2019a} on linear stability of the homogeneous steady state near the Turing instability, the existence and large wavelength in/stabilities of stripes. The stabilities against lattice modes are discussed in \S\ref{s:stability}. In \S\ref{s:example}, we illustrate these results by a concrete example of the form \eqref{e:RDS} and in \S\ref{s:Klausmeier} we study the lattice stability numerically for the extended Klausmeier model that was used in~\cite{Siero2015}.

\subsection*{Acknowledgements}
J.Y was funded by the China Scholarship Council,  and Degree completion stipend from University of Bremen. J.Y is grateful for the hospitality and support from Faculty 3 -- Mathematics, University of Bremen as well as travel support through an Impulse Grants for Research Projects by University of Bremen. J.R. acknowledges this paper is a contribution to the project M2 (Systematic multi-scale modelling and analysis for geophysical flow) of the Collaborative Research Centre TRR 181 ``Energy Transfer in Atmosphere and Ocean" funded by the Deutsche Forschungsgemeinschaft (DFG, German Research Foundation) under project number 274762653. E.S. was supported by a postdoctoral fellowship from the Alexander von Humboldt Foundation.

\section{Preliminaries}\label{s:pre} We place into context some fundamental results for the striped solutions to \eqref{e:RDS}, including the Turing instability, the bifurcation of stripes, Eckhaus and zigzag instabilities. More details of this section can be found in~\cite{Yang2019a}.

\subsection{Turing instability}\label{s:Turing}

The linearisation of \eqref{e:RDS} in $u_{\rm hom}=0$ is 
\[
\calL:=D\Delta+\A + \calpha \M + \beta \B\partial_x,
\]
whose spectrum is most easily studied via the Fourier transform
\[
\hat{\calL}(k,\ell) = -(k^2+\ell^2)D+\A + \calpha \M + \rmi k\beta \B,
\]
with Fourier-wavenumbers $k$ in $x$-direction and $\ell$ in $y$-direction. It is well known, e.g., \cite{Sandstede2002},  that in the common function spaces such as $\Lspace^2(\R^2)$ the spectrum $\Sigma(\calL)$ of $\calL$ equals that of $\hat\calL$ and is the set of roots of the (linear) dispersion relation
\begin{align}\label{e:disp}
d(\lambda,k,\ell) = \det(\hat\calL(k,\ell) -\lambda\Id).
\end{align}
Let $S_{\kc}\subset\R^2$ be the circle of radius $\kc$.

\begin{definition}\label{def:Turing}
We say that $\calpha=\beta=0$ is a (non-degenerate) Turing instability point for $u_{\rm hom}$ in \eqref{e:RDS} with wavelength $\kc$ if 
\begin{itemize}
\item[(1)]  $\A$ has strictly stable spectrum $\Sigma(\A)\subset\{\lambda\in\C : \Re(\lambda)<0\}$,
\item[(2)] The spectrum of $\calL$ is critical for wavevectors $(k,\ell)$ of length $\kc>0$: 
\[
d(\lambda,k,\ell) =0 \ \&\ \Re(\lambda)\geq 0\quad \Leftrightarrow\quad \lambda=0,\ (k,\ell)\in S_{\kc}
\]
which in particular means $\Sigma(\calL)\cap \{z\in \C: \Re(z)\geq 0\} = \{0\}$,
\item[(3)] $\partial_\lambda d\neq 0$ at $\lambda=0$ and $(k_\crit,\ell_\crit)\in S_{\kc}$. We denote the unique continuation of these solutions to \eqref{e:disp} by $\lambda_\crit(k,\ell;\calpha,\beta)$, i.e., $(k,\ell)$ in a neighboorhood of $S_{\kc}$.
\end{itemize}
\end{definition}

Writing $\A=\begin{pmatrix}a_1 & a_2\\ a_3 & a_4\end{pmatrix}$, condition (1) implies negative trace of $\A$, $a_1+a_4<0$, and positive determinant $a_1a_4>a_2a_3$, and (3) implies the well known condition $d_1a_4+d_2a_1>0$, which together imply $a_2a_3<a_1a_4<0$, e.g., \cite{Murray2003}.

As a first step to understand the impact of advection, we next quote basic results from~\cite{Yang2019a}. In particular, for this two-component case the unfolding by $\beta$ is only to quadratic order.

\begin{lemma}[\cite{Yang2019a}]\label{l:Turbeta} 
For the critical eigenvalues near a Turing instability of \eqref{e:RDS} as in Definition~\ref{def:Turing} it holds for any  $(k_\crit,\ell_\crit)\in S_{\kc}$ that
\[
\lambda_\crit(k_\crit,\ell_\crit;\beta)=\rmi k_\crit (\lambda_{\beta}+c) \beta + k_\crit^2\lambda_{\bbeta} \beta^2 + \calO(|k_\crit\beta|^3),
\]
where $\lambda_{\beta}=\frac{a_4-\kcsq d_2}{a_1+a_4-\kcsq(d_1+d_2)}$, $\lambda_{\bbeta}= \frac{(a_1-\kcsq d_1)(a_4-\kcsq d_2)}{(a_1+a_4-\kcsq(d_1+d_2))^3}>0$. In particular, the real part grows fastest for 1D-modes with $\ell_\crit=0$ and remains zero for transverse modes with $k_\crit=0$.
\end{lemma}

\begin{remark}\label{r:Turing}
It is well known that for a two-component system $\kcsq = \frac{d_1 a_4+d_2 a_1}{2d_1d_2}$ and $a_2a_3= (a_1-\kcsq d_1)(a_4-\kcsq d_2)$, cf.\ \cite{Murray2003}.
\end{remark}

\begin{remark}[\cite{Yang2019a}]\label{r:ev}
For the non-trivial matrix $\begin{pmatrix}\ta_1 & \ta_2\\ \ta_3 & \ta_4\end{pmatrix}$ with $b_1b_4-b_2b_3=0$, $b_1\neq0$ and $b_1+b_4\neq0$, we can choose the kernel eigenvector $\E_0$ and the adjoint kernel eigenvector $\E_0^*$ with $\langle \E_0,\E_0\rangle=1$ and $\langle \E_0,\E_0^* \rangle=1$ as
\begin{align*}
E_0 = (b_2,-b_1)^T/c_0, \quad E_0^* = (b_3,-b_1)^T/c_0^*,
\end{align*}
with $c_0 := \sqrt{b_2^2 + b_1^2}$, $c_0^* := (b_1b_4 + b_1^2)/c_0$ and $c_0^*\neq0$.
\end{remark}

\medskip
As expected and observed in \cite{Yang2019a}, in contrast to $\beta$, the change of real parts of the critical eigenvalue through $\calpha$, with matrix $M=(m_{ij})_{1\leq i,j\leq 2}$, is linear with coefficient 
\begin{equation}
\begin{aligned}
\lambda_\M&:=-\left.\frac{\partial_\calpha d}{\partial_\lambda d}\right|_{\calpha=0,\lambda=0}
=\frac{m_{11}(a_4-\kcsq d_2)-m_{12}a_3-m_{21}a_2+m_{22}(a_1-\kcsq d_1)}{a_1+a_4-\kcsq (d_1+d_2)}\neq 0,
\end{aligned}
\end{equation}
where we \emph{assume} $\lambda_M\neq 0$ throughout this paper. Notably, $\lambda_\M=1$ if $\M=\Id$ in which case $\calpha$ just rigidly moves the real part of the spectrum.

In the following we therefore change parameters and use the effective impact on the real part given by 
\[
\alpha:=\lambda_\M\calpha
\]
as the new parameter so that
\begin{equation}\label{e:lamcrit1}
\begin{aligned}
\lambda_\crit(k_\crit,\ell_\crit;\alpha,\beta) = \ &\alpha + \rmi (k_\crit (\lambda_{\beta}+c) + a_\M\lambda_{\M\beta}\alpha)\beta + k_\crit^2\lambda_{\bbeta} \beta^2\\
& + \calO(a_\M\alpha^2 + |k_\crit\beta|^3),
\end{aligned}\end{equation}
with $\lambda_{\M\beta} :=  k_\crit\frac{m_{22}-\lambda_\M-(2\lambda_\M-m_{11}-m_{22})\lambda_{\beta}}{\lambda_\M(a_1+a_4-\kcsq(d_1+d_2))}$, and we emphasise the special case $M=\Id$ through the factor $a_M$, where $a_M=0$ if $M=\rm{Id}$ and $a_M=1$ otherwise. 

\medskip
We illustrate the region where $\Re(\lambda)\geq0$ in the wavevector space $(k,\ell)$ in Fig.~\ref{f:homspec} based on the example \eqref{e:example} below. For relatively small $\beta=0.2$, the Fourier modes with (the leading order) wavevectors $(\pm\kc,0)$ become unstable for $\alpha>-k_\crit^2\lambda_\bbeta\beta^2\approx-0.00448$, where $\alpha=12.24\calpha$, thus the stripes with the wavelength $2\pi/\kc$ may bifurcate from the homogeneous steady state. Then the modes with (the leading order) wavevectors $(\pm\kc/2,\pm\sqrt{3}\kc/2)$ become unstable at $\alpha=-k_\crit^2\lambda_\bbeta\beta^2/4\approx-0.00112$ so that the hexagons may bifurcate. Lastly, the modes with wavevectors $(0,\pm\kc)$ become unstable for $\alpha=0$ and thus squares may bifurcate. 
\begin{figure}[!t]
\centering
\includegraphics[width=0.4\linewidth]{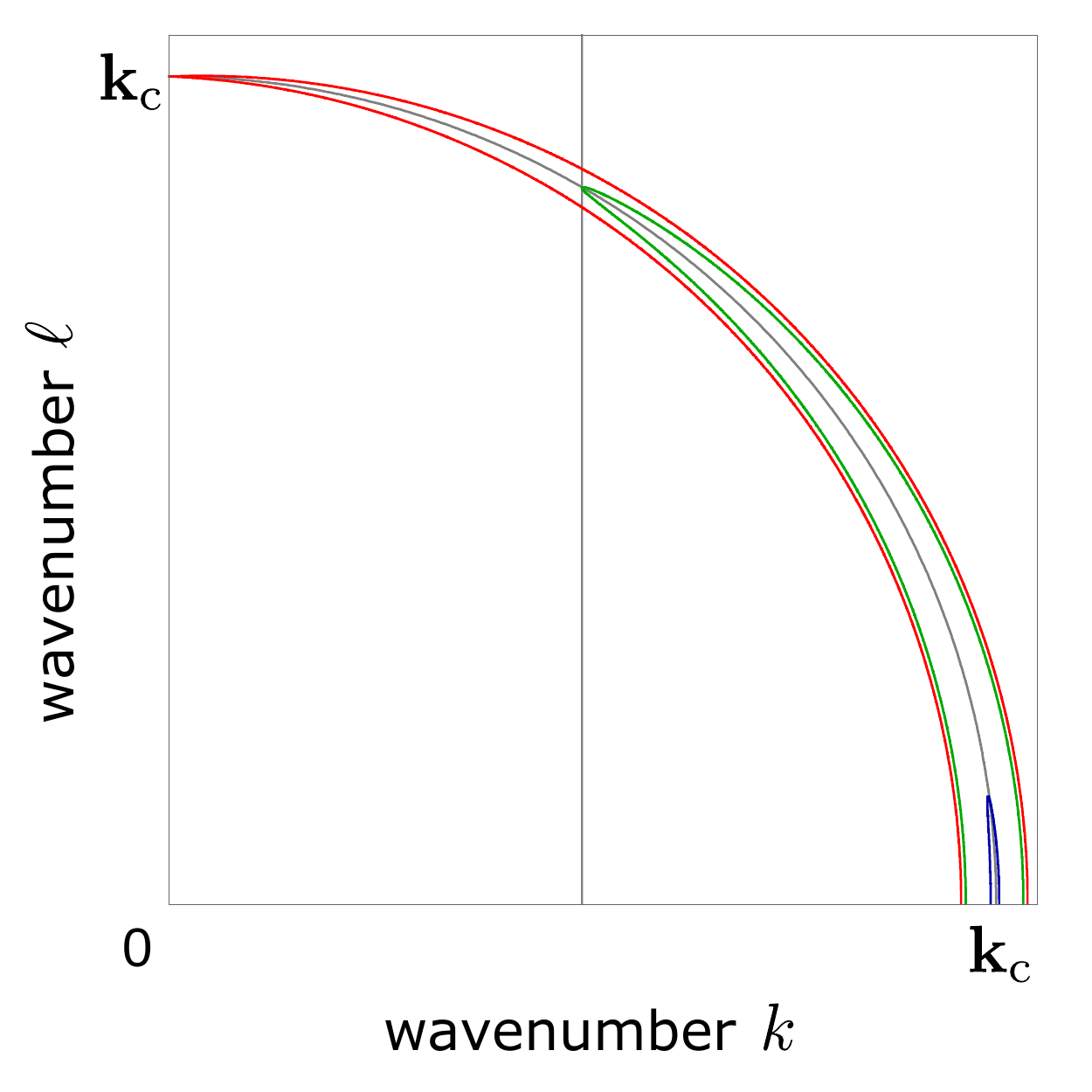}
\caption{Based on the example \eqref{e:example} below, we illustrate the locations of the critical spectrum of homogeneous state, i.e., $\Re(\lambda(k,\ell;\alpha,\beta))=0$, on the $(k,\ell)$-plane for relatively small $\beta=0.2$, and $\Re(\lambda)>0$ inside each horn-shaped region. The unfolding parameter $\alpha=12.24\calpha$. Grey curve: wavevectors $(k_\crit,\ell_\crit)\in S_\kc$ with radius $\kc=1$; grey vertical line: $k=\kc/2=1/2$; blue: $\calpha = -3.6\times10^{-4}$ ($\alpha\approx-0.00441$); green: $\calpha = -9.15\times10^{-5}$ ($\alpha\approx-0.00112$); red: $\calpha=\alpha=0$. These contours are reflection symmetric with respect to both axes.}
\label{f:homspec}
\end{figure}

\subsection{Bifurcation of stripes}\label{s:bif}

With the directional advection, the system \eqref{e:RDS} possesses striped solutions perpendicular to the direction of the advection as proven in~\cite{Yang2019a}. We next briefly recall this existence result and introduce notations that will be used in the stability analysis. We later perform a centre manifold reduction just for the stability analysis. 

\medskip
Stripes are travelling wave solutions of \eqref{e:RDS} that are constant in $y$ and periodic in $x$ for any $t$. It is therefore sufficient for establishing existence to consider the 1D case $\x=x\in[0,2\pi/\kappa]$ with periodic boundary conditions and wavenumber $\kappa$. A Turing instability point as defined above implies that the 1D realisation of $\calL$ possesses a kernel at $\alpha=\beta=0$ on spaces of $2\pi/\kc$-periodic functions, which motivates rescaling space to $[0,2\pi]$ so the linear part \eqref{e:RDS} becomes 
\[
\calL_\mu:= \kappa^2 D\partial_x^2  + \A + \calpha \M  + \beta \kappa\B \partial_x,
\]
with the parameter $\kap$ in $\kappa= \kc+\kap$ that allows to detects stripes with nearby wavenumber. Hence, the main parameters are $\mu=(\alpha,\beta,\kap)$ in vector form. The continuation of the zero eigenvalue of $\calL_\mu$ with $\kap=0$ was given in \eqref{e:lamcrit1}. In the presence of non-trivial $\kap$, such continuation can be expressed in the following form with the notations used throughout this paper
\begin{equation}\label{e:evlinearop}
\begin{aligned}
	\lambda_\mu =&\; \alpha + \rho_\beta\beta^2+\rho_\kap\kap^2 + \rmi(\gamma_\beta+\gamma_{\kap\beta}\kap + a_\M\lambda_{\M\beta}\alpha)\beta \\ 
	& +a_\M \lambda_{\M\kap}\alpha\kap+ \calO(a_\M\alpha^2 + |\kap|^3+ |\beta|^3),
\end{aligned}
\end{equation}
with $a_M=0$ if $M=\rm{Id}$ and $a_M=1$ otherwise. The coefficient $\rho_\beta$ corresponds to that in \eqref{e:lamcrit1} with $\rho_\beta>0$, and $\rho_\kap$ is given by $\rho_\kap = -\partial_k^2d/(2\partial_\lambda d)|_{k=\kc,\lambda=0}<0$; the other coefficients will not be relevant for the leading order stability and  explicit formulae can be found in \cite[Theorem 3.1]{Yang2019a}.

Let us denote by $A\geq 0$ the amplitude of the bifurcating striped solutions. Throughout this paper, we use the scaling
\begin{align}\label{e:scaling}
(A,\alpha, \beta, \kap) = (\scal A', \scal^2 \alpha', \scal \beta', \scal \kap'), 
\end{align}
with a scaling parameter $\scal> 0$ and consider primed quantities $A'$ and $\mu'=(\alpha', \beta', \kap')$ bounded with respect to $\scal$. This scaling is homogeneous for $\mu$ with respect to the relevant first three terms in \eqref{e:evlinearop} and the scaling $A=\scal A'$ is natural due to the relation between the parameters and the amplitude of the striped solutions, cf.~\cite{Yang2019a}. Notably, the impact of $M\neq \Id$ is now at higher order and highlights that leading order results in the following will have additional symmetry.

Using the scaling \eqref{e:scaling}, vanishing the real part $\Re(\lambda_\mu) = 0$ thus occurs to leading order on a hyperbolic paraboloid 
\begin{align}\label{e:bifbnd}
\alpha =\calB(\kap,\beta) := -(\rho_\kap\kap^2 +\rho_\beta\beta^2)
\end{align}
in $\mu$-space. Since the eigenvalue is stable (unstable) for $\alpha<\calB(\kap,\beta)$ ($\alpha>\calB(\kap,\beta)$), this constitutes the bifurcation surface at leading order. Notably, it includes $\mu'=0$ since the signs of $\rho_\beta$ and $\rho_\kap$ are opposite.

In preparation of formulating the existence theorem for this scaling, and for completeness we define a number of quantities that appear in the expansions near $\mu=0$ and $0\leq \scal\ll1$ evaluated at $\mu=0$; most frequently used are the first four.
\begin{equation}\label{e:defs} 
\begin{aligned}
\tq&:= \langle \Q[\E_0, \tQ],\E_0^* \rangle,\qquad\quad
\hq:= \langle \Q[\overline{\E_0}, \hQ],\E_0^* \rangle,\\
\hK &:= \langle \K[\E_0,\E_0,\overline{\E_0}], \E_0^*\rangle, \quad 
\rho_\nl:= 3\hK + 2\tq + \hq,\\
\tQ &:= -2\A^{-1} \Q[\overline{\E_0},\E_0],\quad
\hQ := -2(-4 \kcsq D + \A)^{-1} \Q[\E_0,\E_0],\\
w_{A\calpha} &:= (-\kcsq D+\A)^{-1} (\langle \M \E_0, \E_0^*\rangle- \M) \E_0,\\
w_{A\beta} &:= \kc(-\kcsq D+\A)^{-1}(\langle \B\E_0,\E_0^*\rangle -\B)\E_0,\\
w_{A\kap} &:= 2\kc(-\kcsq D+\A)^{-1} D \E_0,\\
w_{A\beta\beta}&:= 2\kc(-\kcsq D+\A)^{-1} (\B w_{A\beta} - \langle \B w_{A\beta},E_0^*\rangle\E_0),\\
e_\mu(x)&:=(\E_0+\calpha  w_{A\calpha} + \rmi\beta w_{A\beta} + \kap w_{A\kap} + \beta^2 w_{A\beta\beta}) \rme^{\rmi x}
\end{aligned}
\end{equation}
Here $-\kcsq D+\A$ has a one-dimensional generalised kernel spanned by $\E_0$, and thus has an inverse from its range to the kernel of the projection $\langle \cdot, \E_0^*\rangle \E_0$, cf. \cite{Yang2019a}.

Concerning the velocity parameter $c$ as given in the following theorem, the evaluation at $\mu=0$ gives $c=-\lambda_\beta$ where $\langle \B\E_0,\E_0^*\rangle =0$.

\begin{theorem}[Stripe existence~\cite{Yang2019a}]\label{t:bif}
Up to spatial translation, non-trivial stripe solutions to \eqref{e:RDS} with parameters $\mu=\scal \mu'$, and sufficiently small $|\mu|, A = \scal A'$ with $\|U_\rms(\cdot;\mu)\|_{\Lspace^2}=\calO(\scal)$ on $[0,2\pi/\kappa]$, are in 1-to-1 correspondence with solutions $A> 0$ to
\begin{equation}\label{e:stripeeqn}
\scal^2\left(\alpha' + \rho_\beta \beta'^2 + \rho_\kap \kap'^2 + \rho_\nl A'^2 + \calO(\scal)\right)= 0.
\end{equation}
Stripes have velocity $\beta c$ with 
\begin{align}\label{e:stripev}
c = -\lambda_{\beta} + \calO(\scal)
\end{align}
and, in this comoving frame, are up to translation of the form 
\begin{align}\label{e:Stripes}
U_\rms(x;\mu) = \scal A'(e_\mu(x) + \overline{e_\mu(x)}) + \scal^2A'^2\left( \frac{1}{2} \hQ\left( \rme^{2 \rmi x} +  \rme^{-2 \rmi x}\right)+ \tQ\right) + \calO(\scal^3).
\end{align} 
Moreover, the coefficients in \eqref{e:stripeeqn} satisfy 
\begin{equation}\label{e:stripecoeffalt}
\begin{aligned}
\rho_\beta =  -\kc\langle \B w_{A\beta}, \E_0^*\rangle, \;
\rho_\kap = -2\kc\langle D w_{A\kap}, \E_0^* \rangle.
\end{aligned}
\end{equation}
\end{theorem}
\begin{proof}
Using \eqref{e:evlinearop} the bifurcation equation in \cite[Theorem 3.1]{Yang2019a} expands as 
\begin{equation}\label{e:stripeeqn2}
\alpha + \rho_\beta \beta^2 + \rho_\kap \kap^2 + \rho_\nl A^2 + \calO\left(A^3+|\mu|(|\alpha| + \beta^2+\kap^2)\right)= 0.
\end{equation}
Then substituting the homogeneous scalings \eqref{e:scaling} into the above equation yields the claimed result.
\end{proof}

Throughout this paper, we assume $\rho_\nl<0$ so that the bifurcation is a generic supercritical pitchfork.

\subsection{Large-wavelength stability}\label{s:sideband}
In order to include all available stability information, in this section we recall the results on large wavelength stability of stripes from \cite{Yang2019a} with the scalings \eqref{e:scaling}. In the present paper, we will compare this with the stability on the lattices under consideration.

\paragraph{Zigzag instability} The stability of stripes against large-wavelength perturbation parallel to the stripes is referred to as the zigzag stability. It is determined by the following curve of spectrum attached to the origin, cf.~\cite[Corollary 4.2]{Yang2019a}
\[
\lambda_\zz(\ell) = \scal\left(\kc\rho_\kap\kap' + \calO(\scal)\right)\ell^2,
\]
where $\ell$, $|\ell|\ll1$, is the wavenumber of the perturbation in $y$-direction. Hence, the zigzag stability boundary with the scalings \eqref{e:scaling} is to leading order given by $\kap'=0$ and $\sgn(\Re(\lambda_\zz)) = -\sgn(\kap')$. This coincides with the well-known result for the Swift-Hohenberg equation, that stretched stripes are zigzag unstable, while stripes are not as sensitive to compression.

\paragraph{Eckhaus instability} The Eckhaus instability (also called sideband instability) arises from large-wavelength perturbations orthogonal to the stripes. It is well known that a supercritical Turing bifurcation for fixed $\beta=\kap=0$ leads to  Eckhaus-stable stripes, and unstable ones for $\kap\neq0$. As shown in \cite[Theorem 4.4]{Yang2019a}, in competition with $\beta$, the Eckhaus in/stability with the scalings \eqref{e:scaling} is determined by the following real parts of a curve in the spectrum attached to the origin
\[
\Re(\lambda_\eh) = -\kcsq\frac{\rho_\kap}{\rho_\nl A'^2}\left(\alpha'+\rho_\beta\beta'^2+3\rho_\kap\kap'^2 +\calO(\scal))\right)\gamma^2,
\]
with  wavenumber $\gamma, |\gamma|\ll1$ of the perturbation in $x$-direction. Vanishing real part gives the Eckhaus stability boundary in terms of unscaled parameters to leading order 
\begin{align}\label{e:eckhaus}
\alpha = \calE(\kap,\beta) = -(3\rho_\kap\kap^2 + \rho_\beta\beta^2).
\end{align}
This boundary is attached to the bifurcation surface at $\kap=0$ since $\calE(0,\beta) = -\rho_\beta\beta^2 = \calB(0,\beta)$ and lies within the existence region since $\calE(\kap,\beta)\geq\calB(\kap,\beta)$.

\section{Stability of stripes on lattices}\label{s:stability}

In this section we analyse the stability of stripes on rectangular domains with periodic boundary conditions that are nearly square or nearly `hexagonal' in the sense that Fourier modes with wavevectors on a nearly hexagonal lattice are permitted. Indeed, in the Fourier picture these domains have wavevectors on a lattice, and the stability can be studied by centre manifold reduction. While this reduction also allows to study other solutions and nonlinear interactions, here we consider the stability of stripes only. Recall that the lattice modes considered on the (nearly) square and (nearly) `hexagonal' domains are referred to as the {\em (quasi-)square} and {\em (quasi-)hexagonal modes}, respectively.

\begin{remark}\label{r:1Dstab}
Corollary 2.6 from \cite{Yang2019a} has the following a priori consequences for the upcoming detailed stability analysis on lattices: on the one hand, for $\beta\neq0$ stripes near bifurcation, i.e., $0<\alpha-\calB(\kap,\beta)\ll 1$, are stable against quasi-square and quasi-hexagon modes. On the other hand, under the scalings \eqref{e:scaling}, the stripes are zigzag-stable for $\kap>0$ and Eckhaus-stable for $\alpha>\calE(\kap,\beta)$ to leading order, cf.~\S\ref{s:sideband}. In the $(\kap,\alpha)$-plane this connects to the point $\kap=0$, $\alpha=\calB(0,\beta)$ so that stripes are $\Lspace^2$-spectrally stable for the parameters in $\{(\kap,\alpha):0<\kap\ll1,\,0<\alpha-\calE(\kap,\beta)\ll1\}$. Indeed, this is reflected in the results plotted in Figures~\ref{f:rectbetn0tel0},~\ref{f:hexEckbetn0qn0gr},~\ref{f:Eckrhombetn0qn0gr}. 
\end{remark}

\subsection{Centre manifold reduction}\label{s:cmf}

In preparation of the concrete cases, we first consider somewhat abstractly centre manifold reductions for \eqref{e:RDS}. Let us denote
\[
\Lc(\mu):= \calL_\mu - \calL_0 = \calpha\M + \beta(\kc+\kap) \B\partial_x + (2 \kc \kap + \kap^2)D\Delta.
\]

\begin{theorem}[Centre manifold reduction]\label{t:cmf}
Consider \eqref{e:RDS} posed on the interval $\Omega_1=[0,2\pi]$ or on a square $\Omega_2=[0,2\pi]^2$ or a rectangle $\Omega_3= [0,4\pi]\times[0, 4\pi/\sqrt{3}]$ on the space $X=(\Lspace^2(\Omega_j))^2$ with periodic boundary conditions and assume a Turing instability occurs at $\mu=0$. The generalised kernel $N$ of the associated realisation of $\calL_0$  and its co-kernel $Y$ have dimension $2j$ on $\Omega_j$, $j=1,2,3$. In all cases, a $2j$-dimensional centre manifold exists for $|\mu|\ll 1$, which is the graph of $\Psi\in \Cspace^2(N\times \Lambda, Y)$ with $\Psi(0,0)=0$, $\partial_u\Psi(0,0)=0$, and the reduced ODE for $u_c(t)\in N$ is of the form 
\[
\dot u_c = f(u_c;\mu):= P\Lc(\mu)(u_c+ \Psi(u_c,\mu)) + P F(u_c+ \Psi(u_c,\mu)),
\]
where $P:X\to N$ is the projection with kernel $Y$. In particular,
\[
\partial_u f(u_c;\mu) = P\big(\Lc(\mu) + \partial_u F(u_c+\Psi(u_c;\mu))\big)(\Id + \partial_u\Psi(u_c;\mu)) + \calO(|u_c|^3).
\]
\end{theorem}
\begin{proof}
It suffices to show the claimed dimension of the kernel depending on $j$; the result then follows from standard centre manifold theory, e.g., \cite{Haragus2010}, by the definition of Turing instability. For $\Omega_1$ this was already discussed in the previous section. From Lemma~\ref{l:Turbeta} the critical eigenmodes of $\calL_0$ are explicitly known, in particular their wavevectors satisfy $(k_j,\ell_j)\in S^1$. Hence, on $\Omega_2$ these are the four choices $\krec_1:=(1,0)$, $\krec_2:=(0,1)$ and their negatives, and on $\Omega_3$ the six choices $\k_1:=(1,0)$, $\k_2:=(-1/2,\sqrt{3}/2)$, $\k_3:=-(1/2,\sqrt{3}/2)$ and their negatives.
\end{proof}

\begin{remark}\label{r:resonances}
As to nonlinear terms of $f$ we note that $Pv=0$ if $v$ consists of Fourier modes whose wavevectors are not in $S^1$, which leads to the following resonance condition. Since  wavevectors are added in products, any nonlinear term must stem from products of terms for which the sum of wavevectors from $S^1$ lies again in $S^1$. Such resonant interactions require at least three terms, and on $\Omega_2$ are possible only among wavevectors in the same spatial direction. In contrast, $\Omega_3$ allows for the so-called resonance triads (or three-wave interactions) $\k_1+\k_2+\k_3=0$.
\end{remark}

Next, we expand the linearisation on the centre manifold somewhat abstractly in order to be conveniently used for different settings later.

Let us denote $\Psi_{j \ell}:=\partial_u^j \partial_\mu^\ell \Psi(0;0)/(j!\ell!)$ so that $\Psi_{00}=\Psi_{10}=0$ in general and due to the zero equilibrium for all parameters also $\Psi_{0 j}=0$ for all $j\geq 0$.

\begin{corollary}\label{c:cmf}
Assume the conditions and notations of Theorem~\ref{t:cmf} and the scaling \eqref{e:scaling} so that $u_c=\scal A' u_1\in N$, $\mu=\scal\mu_1+ \scal^2\mu_2$ and $u_1, \mu_1,\mu_2=\calO(1)$ with respect to $\scal$. 
We have $\Psi(u_c;\mu)=\scal^2u_2 +\calO(\scal^3)$ with $u_2:=A'^2\Psi_{20}[u_1,u_1]+A'\Psi_{11}[\mu_1,u_1]$, and it holds that
\begin{equation}\label{e:genlin}\begin{aligned}
\partial_u f(u_c;\mu) =\ & 2\scal A'PQ[u_1,\cdot] \\
& + \scal^2 P\Big(\Lc(\mu_2)+(\Lc(\mu_1)+2A'Q[u_1,\cdot])(2A'\Psi_{20}[u_1,\cdot] + \Psi_{11}[\mu_1,\cdot])\\
& + 2Q[u_2,\cdot]  + 3A'^2 K[u_1,u_1,\cdot] )\Big) + \calO(\scal^3).
\end{aligned}\end{equation}
\end{corollary}

\begin{proof}
Substituting $u_c$, $\mu$ as assumed gives $\Lc(\mu)= \scal\Lc(\mu_1) + \scal^2\Lc(\mu_2)$ and Taylor expanding  $\Psi(u_c;\mu)=\scal^2u_2 +\calO(\scal^3)$ as well as 
\begin{align*}
\partial_u\Psi(u_c,\mu) &= \scal(2A'\Psi_{20}[u_1,\cdot] + \Psi_{11}[\mu_1,\cdot]) + \calO(\scal^2),\\
\partial_u F(u_c+\Psi(u_c,\mu)) &= 2\scal A'Q[u_1,\cdot] + 2\scal^2Q[u_2,\cdot] + 3\scal^2A'^2K[u_1,u_1,\cdot]+\calO(\scal^3).
\end{align*}
Combining these, $\partial_u f$ from Theorem~\ref{t:cmf} and using that $\langle B\E_0 ,\E_0^*\rangle|_{\mu=0}=0$ and $\langle D\E_0,\E_0^* \rangle=0$, which removes $P\Lc(\mu_1)$, we obtain the claimed form.
\end{proof}

\subsection{Stability in one space-dimension}

We first note that due to lack of triads, cf. Remark~\ref{r:resonances} a number of terms in  \eqref{e:genlin} vanish:
$U=U_0\rme^{\rmi x}$ with any $U_0\in\C^{2}$ gives $PQ[U,\cdot]=0$ on $N$. Analogously, $Q[u_1,\cdot]$, $Q[u_1,\Psi_{11}[\mu_1,\cdot]]$, $\Lc(\mu_1)\Psi_{20}[u_1,\cdot]$, $\Q[\Psi_{11}[\mu_1,u_1],\cdot]$ vanish so that  \eqref{e:genlin} simplifies to  
\begin{equation}\label{e:genlin1D}
\begin{aligned}
\partial_u f(u_c;\mu) =\ &\scal^2 P\Big(\Lc(\mu_2)+\Lc(\mu_1)\Psi_{11}[\mu_1,\cdot] + 2A'^2 \Q[\Psi_{20}[u_1,u_1],\cdot]\\ 
& + 4A'^2 \Q[u_1,\Psi_{20}[u_1,\cdot]] + 3A'^2 \K[u_1,u_1,\cdot] \Big) + \calO(\scal^3).
\end{aligned}
\end{equation}

Next we infer the matrix form of the linearisation from the existence result. It is convenient to also span the centre eigenspace by $\sin$ and $\cos$, i.e., $u_c=u_0 \cos + u_1 \sin$ for $u_0,u_1\in\R$; the projection in these coordinates is given by $P:=\Id-P_h = \langle \cdot, E_0^*\cos\rangle\cos +  \langle \cdot, E_0^*\sin\rangle\sin$; and, up to translation in $x$, stripes are given by
\begin{align*}
	U_\rms(x;\mu) =\ &2\scal A'\E_0\cos(x)\\
	& + 2\scal^2A'\Big(\kap'w_{A\kap}\cos(x) - \beta'w_{A\beta}\sin(x) + A'(\Q_2\cos(2x)+\Q_0)\Big)+\calO(\scal^3).
\end{align*}

\begin{theorem}\label{t:cmfstripe}
Assume the conditions and notations of Theorem~\ref{t:cmf} for the domain $\Omega=[0,2\pi]$ with periodic boundary conditions and velocity parameter $c=c(\mu)$ as in \eqref{e:stripev}. Stripes $U_\rms$ are in 1-to-1 correspondence with equilibria $u_c\in N$, $f(u_c,\mu)=0$, $\mu$ solving \eqref{e:stripeeqn} and, up to translation in $x$, $U_\rms=u_c + \Psi(u_c;\mu)$ for $u_c = 2 AE_0\cos(x)$. 
The linearisation in stripes satisfies $\partial_u f(u_c;\mu)E_0\sin = 0$ as well as $\partial_u f(u_c;\mu)E_0\cos = 2 A^2 \rho_\nl + \calO(\scal^3)$ with \eqref{e:scaling}, and, up to this order, has the matrix forms 
\[
A^2\begin{pmatrix}2\rho_\nl & 0\\ 0 & 0\end{pmatrix}, \quad A^2 \begin{pmatrix}\rho_\nl & \rho_\nl\\ \rho_\nl & \rho_\nl\end{pmatrix}, 
\]
in the coordinates $\cos, \sin$ and $e_0, \overline{e_0}$, respectively.
\end{theorem}

\begin{proof}
Centre manifold equilibria $f(u_c;\mu)=0$ correspond to equilibria near bifurcation and, due to Theorem~\ref{t:bif}, these are stripes so that $u_c=2A\cos$. By translation symmetry $P\partial_x U_\rms$ lies in the kernel of $\partial_u f(2A\cos;\mu)$ and in particular each expansion order with respect to $A$ of the linearisation has the corresponding order of $P\partial_x U_\rms$ as its kernel. In fact, due to the translation symmetry of \eqref{e:RDS}, the ODE in Theorem~\ref{t:cmf} is independent of the translation direction, cf.\ \cite{Haragus2010}. Hence, the matrix is diagonal in $(\cos,\sin)$-coordinates and it remains to determine the second eigenvalue. In this reduced equation, the bifurcation of stripes is a generic pitchfork with $\lambda_\mu$ the normal form unfolding parameter, and it is well known that the eigenvalue of the bifurcating branch is to leading order $-2\lambda_\mu=2\rho_\nl A^2$ \cite{Haragus2010}.
\end{proof}

\begin{remark}
The proof for the matrix form in Theorem~\ref{t:cmfstripe} does not rely on the detailed expansion of the linearisation  \eqref{e:genlin1D}, but can of course be derived from it. This is somewhat tedious since $\Psi_{20}, \Psi_{11}$ enter in general, and we do this for the hexagonal lattice in Appendix~\ref{s:rectmat}. 
\end{remark}

\subsection{Stability against (quasi-)square perturbations}\label{s:square} 
We start with the simplest case, the stability against (quasi-)square perturbations. Although it turns out that these are not the dominant instability mechanisms among planar modes, it is instructive and adds to completeness of the analyses of lattice modes.

\medskip
We consider the problem \eqref{e:RDS} with periodic boundary conditions on the (quasi-)square domain 
\[
\Omega_\sq:=[0,2\pi/\kappa]\times[0,2\pi/\ell],\; \kappa:=\kc+\kap,\, \ell:=\kc+\tel,
\]
with the scaling $\tel = \scal\tel'$ in accordance with \eqref{e:scaling}, so that $\tel = \calO(\scal)$. In particular, the quasi-square domain reduces to the square domain when $\kappa=\ell$. Rescaling the spatial variables with $\tilde x=x/\kappa$ and $\tilde y=y/\ell$, so that the scaled domain is given by $\Omega_2=[0,2\pi]^2$ with dual lattice wavevectors $\krec_j = (k_j, \ell_j)\in \R^2$, where
\[
\krec_1=(1,0), \; \krec_2=(0,1),
\]
and for convenience $\krec_{-j} = -\krec_j$, $j=1,2$. As noted in Theorem~\ref{t:cmf} this leads to a four-dimensional centre manifold for
\[
u_c(x) = U_\sq(x) = \sum_{j=-2, j\neq 0}^2 u_j e_j,
\]
where $u_j=\overline{u_{-j}}\in \C$ and $e_j:=\rme^{\rmi \krec_j\cdot \x} \E_0$ are the four linearly independent kernel eigenvectors that appear for $\Omega_2$; we also denote $e_j^*:=\rme^{\rmi \krec_j\cdot \x} \E_0^*$. 

\begin{theorem}\label{t:cmfrect}
Assume the conditions and notations of Theorem~\ref{t:cmf} for the domain $\Omega_2$ with periodic boundary conditions and the scaling \eqref{e:scaling}. Let the velocity parameter $c=c(\mu)$ be as in \eqref{e:stripev}. The subspace $\{ u_j=0, j=\pm 2\}$ is invariant for the reduced ODE and contains the stripes as equilibria. The linearisation in stripes in the index ordering $(1, -1, 2, -2)$ has a block diagonal matrix of the form $L_\sq = \diag(L_1, L_2^\sq)+\calO(\scal^3)$ with $2\times 2$-subblocks
\begin{align*}
L_1 =
A^2\begin{pmatrix}\rho_\nl & \rho_\nl\\ \rho_\nl & \rho_\nl\end{pmatrix}, \quad
L_2^\sq=\scal^2
\begin{pmatrix}
\lambda'_\tel + A'^2\xi & 0\\
0 & \lambda'_\tel + A'^2\xi
\end{pmatrix},
\end{align*}
where $\lambda'_\tel:= \alpha' + \rho_\kap\tel'^2 + \calO(\scal)$ is real and 
\begin{align*}
\xi := 6 \hK + 2\tq + 8 q_{11},\; q_{11} :=  \langle \Q[E_0,\Q_{11}],E_0^* \rangle, \;\Q_{11}&:=-(-2\kcsq D+L)^{-1} Q[\E_0,\E_0].
\end{align*}
\end{theorem}

See Appendix~\ref{s:rectmat} for the proof.

The eigenvalues of the matrix $L_1$ are $0$ and $2\rho_\nl A^2<0$, as in Theorem~\ref{t:cmfstripe}, which reflects that the stripes are stable against perturbations in the $x$-direction on this domain.

Concerning the subblock $L_2^\sq$, we first note the general form of eigenvalues.
\begin{lemma}\label{l:rectevals}
Under the assumptions of Theorem~\ref{t:cmfrect}, the double eigenvalue of the matrix $L_2^\sq$ is real and given by
\[
\lambda = \scal^2\left(A'^2(3\hK - \hq + 8 q_{11}) - \rho_\beta \beta'^2 + \rho_\kap(\tel'^2-\kap'^2)\right) + \calO(\scal^3),
\]
where $A'=\sqrt{-(\alpha'+\rho_\beta \beta'^2 + \rho_\kap \kap'^2)/\rho_\nl} + \calO(\scal)$. 
\end{lemma}
\begin{proof}
The two eigenvalues are the same diagonal term which, by \eqref{e:stripeeqn}, read
\begin{align*}
\scal^2(\lambda'_\tel + A'^2\xi) = \scal^2 (-A'^2\rho_\nl -\rho_\beta \beta'^2 - \rho_\kap\kap^2+\rho_\kap \tel^2 + A'^2\xi )
+ \calO(\scal^3),
\end{align*}
and \eqref{e:stripeeqn} gives $A'$ as claimed; that $\lambda_\tel'$ is real was already stated in Theorem~\ref{t:cmfrect}.
\end{proof}

We note that the signs of $\hq,q_{11}$ depend on $\Q$. For the sake of simplicity and comparison with (quasi-)hexagonal stabilities discussed below, we consider Hypothesis~\ref{h:Qscale}. This immediately gives the following stability result.

\begin{corollary}[(Quasi-)square lattice stabilities]\label{c:rectev}
Under the assumptions of Theorem~\ref{t:cmfrect} and Hypothesis~\ref{h:Qscale} the double eigenvalue of matrix $L_2^\sq$ is given by
\begin{align}\label{e:evrect}
\lambda_\sq = \scal^2\left(-\alpha' - 2\rho_\beta\beta'^2 + \rho_\kap(\tel'^2-2\kap'^2) \right)+o(\scal^2).
\end{align}
\end{corollary}

\begin{figure}[!t]
\centering
\subfloat[$\beta=0$]{
\includegraphics[width=0.35\linewidth]{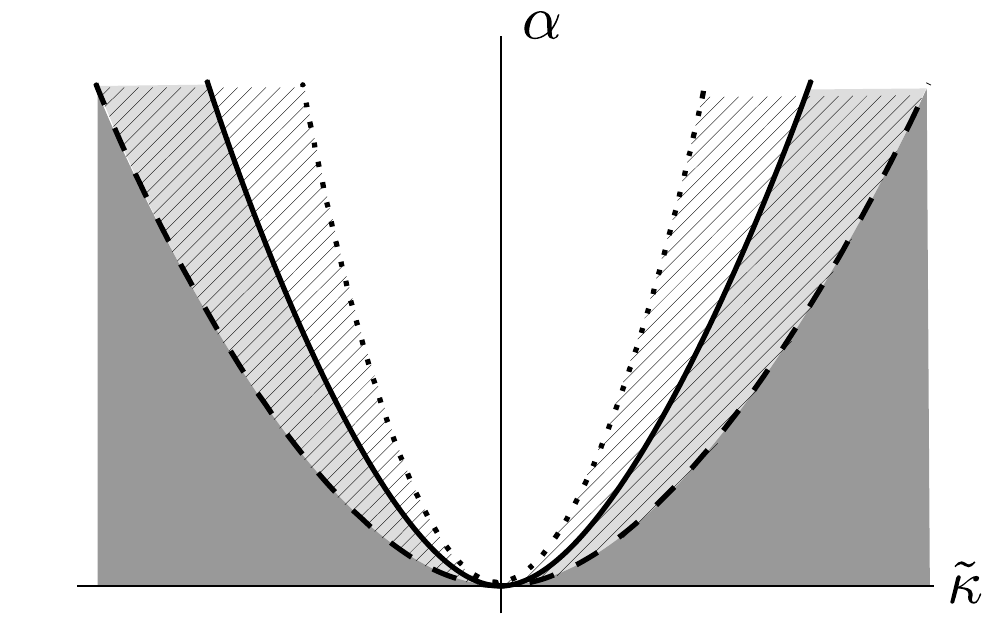}\label{f:rectbet0tel0}}
\hfil
\subfloat[$\beta\neq0$]{\includegraphics[width=0.35\linewidth]{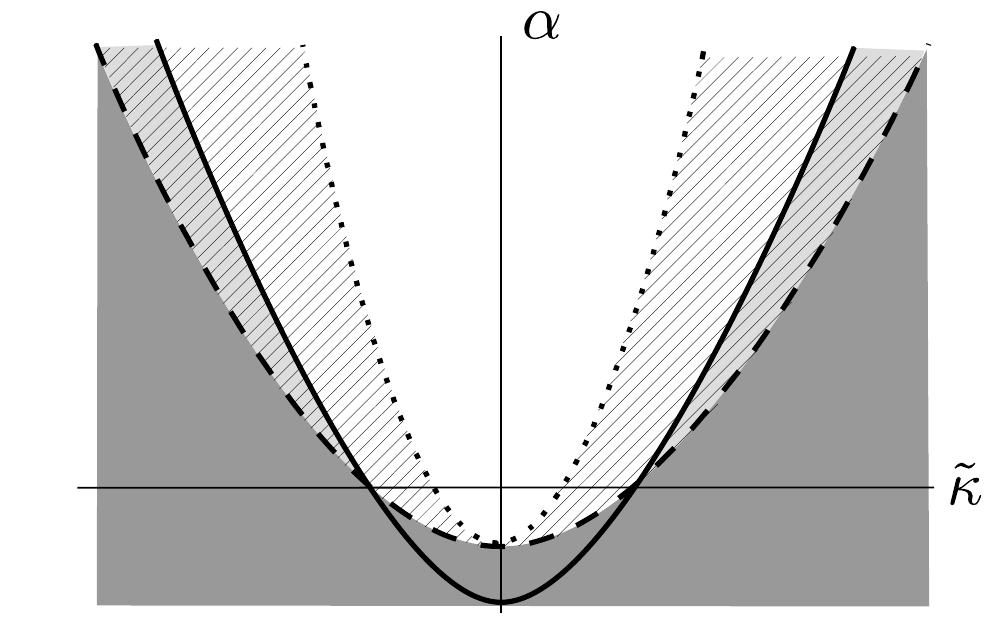}\label{f:rectbetn0tel0}}
\caption{Sketches of the stability regions in the $(\kap,\alpha)$-plane for $\tel=0$. Stripes exist in the complements of the dark grey regions. Light grey: quasi-square-unstable; hatched region: Eckhaus-unstable; white: stable; dashed curve: bifurcation curve $\alpha=\calB(\kap,\beta)$; dotted curve: Eckhaus boundary $\alpha = \calE(\kap,\beta)$; black solid curves: quasi-square stability boundary $\alpha = \calQ(\kap,\beta;0)$.}
\label{f:rectbnd}
\end{figure}

The stability boundary $\lambda_\sq=0$ in terms of unscaled parameters is given by
\begin{align}\label{e:rectbnd}
\alpha = \calQ(\kap,\beta;\tel) := -2\rho_\beta \beta^2 + \rho_\kap (\tel^2-2\kap^2).
\end{align}
For any $\kap$ and fixed $\beta$, if $|\tel|\geq|\kap|$, then $\calQ(\kap,\beta;\tel)\leq\calB(\kap,\beta)$ (see \eqref{e:bifbnd}) and thus the stripes are stable against the square perturbation and the quasi-square perturbations with $|\tel|>|\kap|$. 
The curvature of $\calQ$ with respect to $\kap$ is larger than that of $\calB$ for $|\tel|<|\kap|$, which causes unstable stripes against quasi-square perturbations. The most unstable quasi-square perturbation occurs at $\tel=0$, cf.~Fig.~\ref{f:rectbnd}. However, note that the Eckhaus boundary is dominant since $\alpha = \calE(\kap,\beta)\geq\calQ(\kap,\beta;0)$.

\subsection{Stability against hexagonal perturbations}\label{s:hex}

Concerning the six-dimensional lattice modes, we first study the exact hexagonal perturbation as a basis for the more unstable  quasi-hexagonal perturbations. On the one hand, it is natural and relatively easy to consider the hexagonal case as, e.g., in the amplitude equations approach. On the other hand, the stability proof is neat and can be extended to the quasi-hexagonal case.

\medskip
We consider the system \eqref{e:RDS} with periodic boundary condition on the rectangular domains 
\[
\Omega_\hex = [0,4\pi/\kappa]\times[0,4\pi/(\sqrt 3 \kappa)],\;\kappa=\kc+\kap,
\]
and isotropically rescale to   
$\Omega_3=[0,4\pi]\times [0,4\pi/\sqrt{3}]$ with dual lattice wavevectors $\k_j = (k_j, \ell_j)\in \R^2$, where 
 $\k_{-j} = -\k_j$, $j=1,2,3$, cf.\ Remark~\ref{r:resonances},
\[
\k_1=(1,0), \; \k_2=(-1/2,\sqrt{3}/2), \; \k_3=-(1/2,\sqrt{3}/2).
\]
As noted in Theorem~\ref{t:cmf} this leads to a six dimensional centre manifold for
\[
u_c(x) = U_\hex(x) = \sum_{j=-3, j\neq 0}^3 u_j e_j,
\]
where $u_j=\overline{u_{-j}}\in \C$ and $e_j:=\rme^{\rmi \k_j\cdot \x} \E_0$ are the six linearly independent kernel eigenvectors that appear for $\Omega_3$; we also denote $e_j^*:=\rme^{\rmi \k_j\cdot \x} \E_0^*$. For convenience, here we use the same notation for the wavevectors and (adjoint) eigenvectors as in \S\ref{s:square}.

\begin{theorem}\label{t:cmfhex}
Assume the conditions and notations of Theorem~\ref{t:cmf} for the domain $\Omega_3$ with periodic boundary conditions, and the parameter scaling \eqref{e:scaling} for $\mu$. Let the velocity parameter $c=c(\mu)$ be as in \eqref{e:stripev}. The subspace $\{ u_j=0, j=\pm 2,\pm 3\}$ is invariant for the reduced ODE and contains the stripes as equilibria. The linearisation in stripes in the index ordering $(1, -1, 2, -3, 3, -2)$ has a block diagonal matrix of the form $L_\hex = \diag(L_1, L_2^\hex, L_2^\hex)+\calO(\scal^3)$ with $2\times 2$-subblocks
\begin{align*}
L_1 =
A^2\begin{pmatrix}\rho_\nl & \rho_\nl\\ \rho_\nl & \rho_\nl\end{pmatrix}, \quad
L_2^\hex =\scal^2
\begin{pmatrix}
\lambda'_{\mu,2} + A'^2\eta & 2A'\frac{q}{\scal} + A' p(\mu_1)\\
2A'\frac{q}{\scal} + A' \overline{p(\mu_1)} & \lambda'_{\mu,2} + A'^2\eta
\end{pmatrix},
\end{align*}
where $\lambda'_{\mu,2}:= \alpha' + \frac{1}{4}\rho_\beta\beta'^2 + \rho_\kap\kap'^2 + \calO(\scal)$ and 
\begin{align*}
q&:=\langle Q[E_0,E_0],E_0^*\rangle,\quad \eta := 6 \hK + 2\tq + 8 q_1,\quad q_1 :=  \langle \Q[E_0,\Q_1],E_0^* \rangle, \\
\Q_1&:=(-\kcsq D+L)^{-1}(\langle Q[\E_0,\E_0],E_0^*\rangle E_0 - Q[\E_0,\E_0]),\\
p(\mu_1) &:= \langle \Q[\rmi\beta'w_{A\beta} + 4\kap'w_{A\kap},\E_0], \E_0^* \rangle + \langle(-4\kap' \kc D - \rmi \beta'\kc B )\Q_1,E_0^*\rangle.
\end{align*}
\end{theorem}
See Appendix~\ref{s:hexmat} for the proof.

\medskip
Since $L_1$ concerns perturbations in the $x$-direction, i.e., orthogonal to stripe, from Theorem~\ref{t:cmfstripe} we know that $L_1$ has the eigenvalues $0$ and $2\rho_\nl A^2<0$.

Concerning the subblock $L_2^\hex$, we first note the general form of eigenvalues.

\begin{lemma}\label{l:hexevals}
Under the assumptions of Theorem~\ref{t:cmfhex}, the eigenvalues of the matrix $L_2^\hex$ are
\[
\lambda_\pm = \scal^2\left(A'^2(3\hK - q_2 + 8 q_1) - \frac 3 4 \rho_\beta \beta'^2 \pm A'\left|\frac{2q}{\scal}+A'p(\mu_1)\right|\right) + \calO(\scal^3),
\]
where $A'=\sqrt{-(\alpha'+\rho_\beta \beta'^2 + \rho_\kap \kap'^2)/\rho_\nl} + \calO(\scal)$.
\end{lemma}
\begin{proof}
The matrix $L_2^\hex$ is of the form $\begin{pmatrix} a & b\\ \bar b & a \end{pmatrix}$ with $a\in\R,\,b\in\C$, and such a matrix possesses the two real eigenvalues $\lambda_\pm=a\pm|b|$. For $\lambda_\pm$, $b$ is as in the matrix unchanged, and for $a$ we have
\begin{align*}
a &= \scal^2(\lambda'_{\mu,2} + A'^2\eta) = \scal^2 \left(-A'^2\rho_\nl - \frac 3 4 \rho_\beta \beta'^2 + A'^2\eta \right)
+ \calO(\scal^3)
\end{align*}
and using \eqref{e:stripeeqn} gives the claimed form. 
\end{proof}

The lemma shows that for small $\scal$ and $q=\calO(1)$ with respect to $\scal$ we have $\lambda_+>0$, and the stripe thus unstable. In order to study destabilisation of stripes near onset, and thus the competition of the quadratic term and advection, we therefore assume $q=\scal q'$ with $q'=\calO(1)$. This is most easily realised by Hypothesis~\ref{h:Qscale}, which assumes the entire quadratic term has a prefactor $\scal$, though we note that $q=\scal q'$ can be realised by a scaling assumption on certain parts of $Q$ only.

\medskip
In this case we can rewrite the entries in $L_2^\hex$ related to $\Q$ as follows
\begin{align*}
q& =\scal q',&
q_1 & = \scal^2 q_1', &
\eta & = 6\hK + \scal^2(2\tq' + 8q_1'),&
\Q_1 & = \scal\Q_1',&
p(\mu_1) & = \scal p'(\mu_1),
\end{align*}
with bounded primed quantities. Moreover, we recall 
\[
\rho_\nl = 3\hK + \scal^2(2\tq' + \hq') < 0
\]
with sign due to the assumed supercriticality of the stripe bifurcating so that also $\hK<0$. This gives the following hexagonal in/stability result.
\begin{theorem}[Hexagonal lattice stability]\label{t:evhex}
Under the assumptions of Theorem~\ref{t:cmfhex} and Hypothesis~\ref{h:Qscale} the eigenvalues of the matrix $L_2^\hex$ are given by
\begin{align}
\lambda_\hex^\pm = \scal^2\left(3\hK \tA'^2 -\frac{3}{4}\rho_\beta\beta'^2\pm 2\tA'|q'| + \calO(\scal)\right),
\end{align}
where $\tA':=\sqrt{-(\alpha'+\rho_\beta\beta'^2+\rho_\kap\kap'^2)/(3\hK)}$. In particular, $\lambda_\hex^\pm\in\R$.
\end{theorem}
\begin{proof}
Using Lemma~\ref{l:hexevals}  the claim directly follows from Hypothesis~\ref{h:Qscale} and the resulting factors of $\scal$ as noted above. The term $\tA'$ stems from the leading order of $A'$, i.e., $A' = \tA'+\calO(\scal)$.
\end{proof}

In particular, under these assumptions, $q'$ is the only relevant quantity that relates to $\Q$. In case $\Q=o(\scal)$ we have $q'=0$ so that $\lambda_\hex^\pm<0$, i.e., stripes are always stable on the hexagonal lattice, since $\hK<0$ and $\rho_\beta>0$.

\medskip
In Theorem~\ref{t:evhex} the eigenvalue $\lambda_\hex^-$ is stable for all $\mu$ and $q$ such that the striped solution \eqref{e:Stripes} exists. The sign of $\lambda_\hex^+$, however, depends on both $\mu$ and $q$. A critical eigenvalue $\lambda_\hex^+=0$ to leading order requires $3\hK \tA'^2-\frac{3}{4}\rho_\beta\beta'^2<0$ or equivalently
\[
\alpha > -\frac{7}{4}\rho_\beta\beta^2 - \rho_\kap\kap^2
\]
in terms of unscaled parameters. Since $-\frac{7}{4}\rho_\beta\beta^2 - \rho_\kap\kap^2<\calB(\kap,\beta)$ the above condition is automatically fulfilled for $\mu$ such that the stripes exist. 

Solving $\lambda_\hex^+=0$ yields the hex-stability boundaries to leading order. In terms of the unscaled parameters this reads
\begin{align} \label{e:hexstabbnd}
\{\alpha=\calH_\pm(\kap,\beta,q): \alpha\geq\calB(\kap,\beta),\;\discrh\geq0\},
\end{align}
where
\begin{align}
\calH_\pm(\kap,\beta,q) &:= -\frac{7}{4}\rho_\beta\beta^2 - \rho_\kap\kap^2 -\frac{1}{3\hK}\left(2q^2 \pm \sqrt{\discrh}\right),\label{e:hexbnd}\\
\discrh&:=4q^4 + 9\hK q^2\rho_\beta\beta^2.\label{e:discrh}
\end{align}
We remark that since $\alpha\in\R$, the condition $\discrh\geq0$ appears. The stripes are hex-unstable for $\discrh>0$ and $\alpha\in(\calH_-,\calH_+)$, and hex-stable otherwise.

In order to simplify notations, we formulate the hex-stability boundaries in terms of the unscaled parameters in \S\ref{s:hexcrit} and \S\ref{s:hexsidemodes}.

\subsubsection{Stripes with critical wavenumber}\label{s:hexcrit}
We first consider the stripes with the Turing critical wavenumber, i.e., $\kap=0$.

\begin{figure}[!t]
\centering
\subfloat[$\beta=0$, $\kap=0$]{\includegraphics[width=0.35\linewidth]{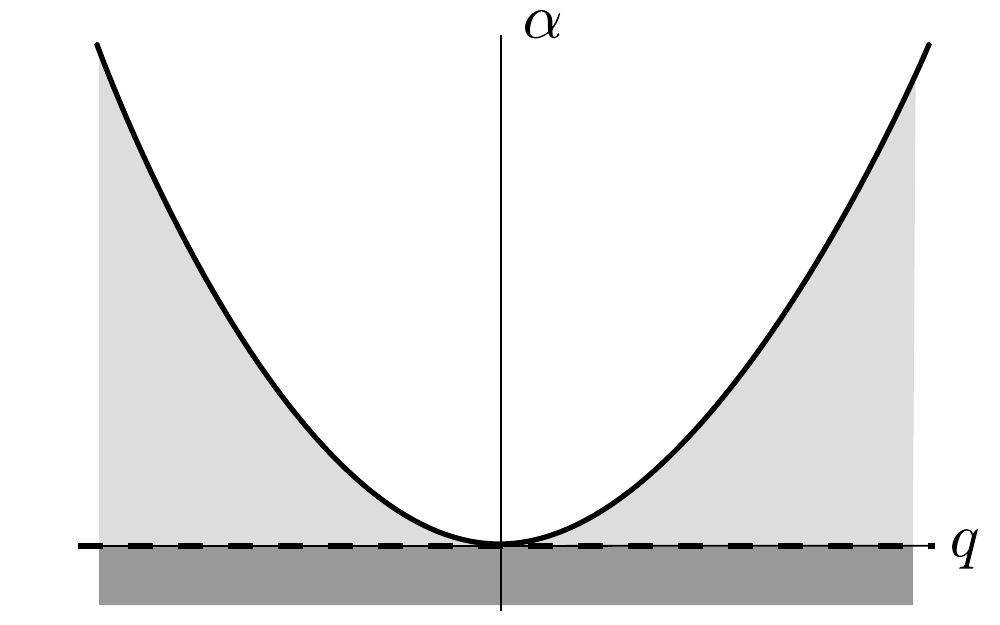}\label{f:hexkap0bet0}}
\hfil
\subfloat[$\beta\neq0$, $\kap=0$]{\includegraphics[width=0.35\linewidth]{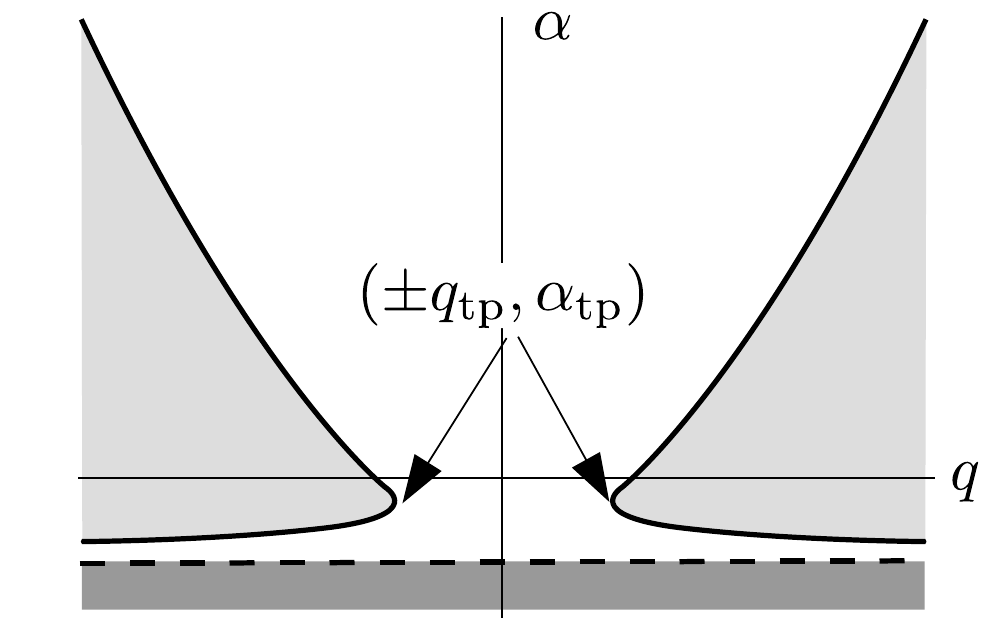}\label{f:hexkap0betn0}}
\hfil
\subfloat[$\beta=0$, $\kap\neq0$]{\includegraphics[width=0.35\linewidth]{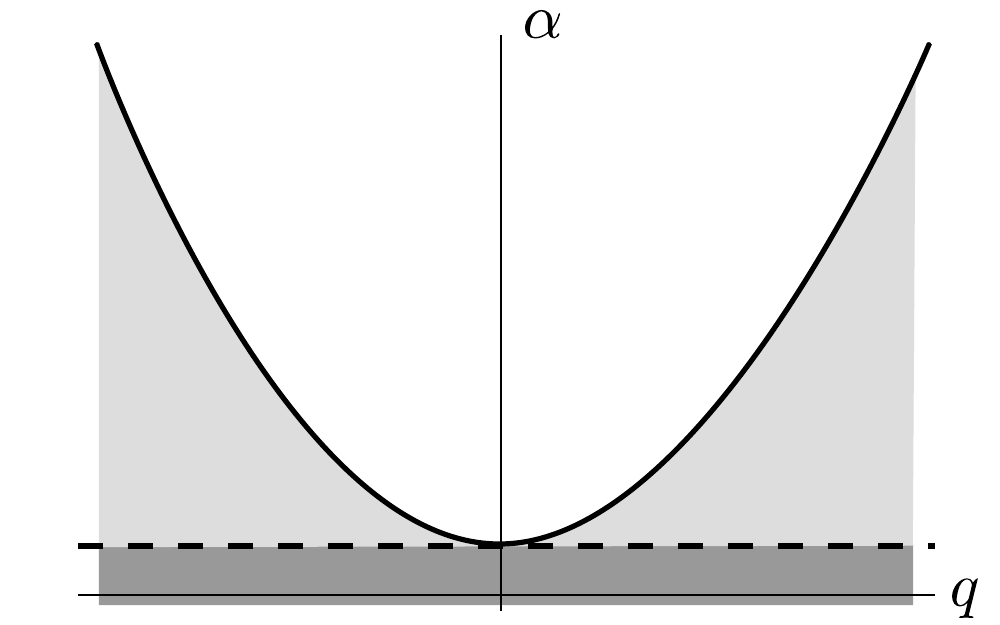}\label{f:hexkapn0bet0}}
\hfil
\subfloat[$\beta\neq0$, $\kap\neq0$]{\includegraphics[width=0.35\linewidth]{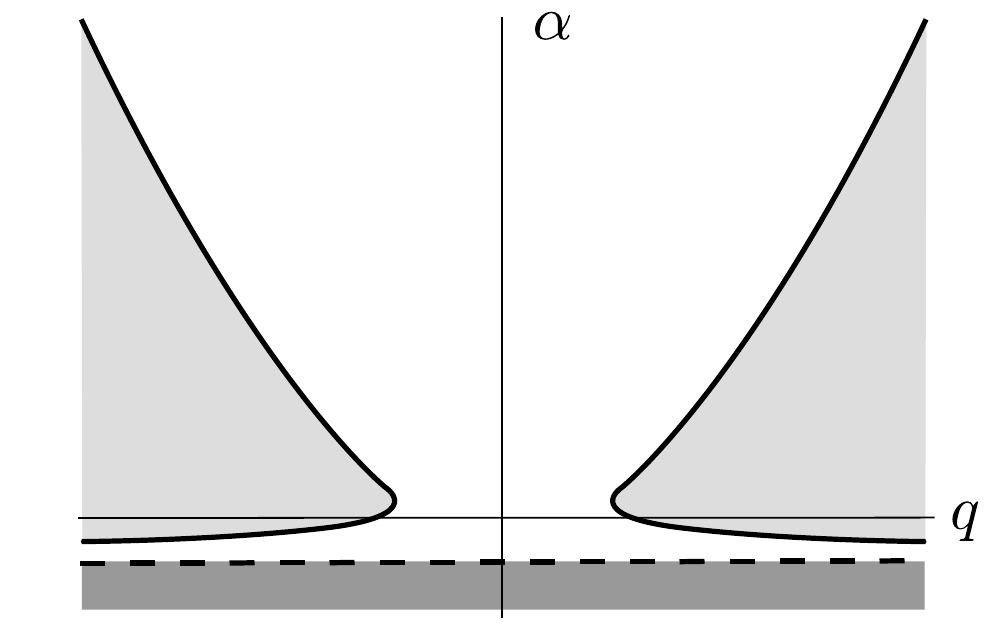}\label{f:hexkapn0betn0}}
\caption{Sketches of the hexagonal stability regions of stripes in the $(q,\alpha)$-plane. Stripes exist in the complement of the dark grey regions. White: hex-stable; grey: hex-unstable. Dashed line: bifurcation line $\alpha=\calB(0,\beta)$; solid curves: hex-stability boundaries (a) $\alpha = \calH_+(0,0,q)$, cf. \eqref{e:stabbndiso} (b) $\alpha = \calH_\pm(0,\beta,q)$, cf. \eqref{e:stabbndaniso} (c) $\alpha = \calH_+(\kap,0,q)$, cf. \eqref{e:hexbndbet0qn0large} (d) $\alpha = \calH_\pm(\kap,\beta,q)$, cf. \eqref{e:hexbnd}.}
\label{f:hexalpq}
\end{figure}

\paragraph{Case $\beta=0$, $\kap=0$ (Fig.~\ref{f:hexkap0bet0})} The hex-stability boundary reduces to a parabola
\begin{align}\label{e:stabbndiso}
\alpha = \calH_+(0,0,q)= -\frac{4}{3\hK}q^2.
\end{align}
This coincides with the well-known result that the stripes are hex-unstable near the onset of Turing bifurcation except for $q=0$~\cite{Gowda2014}. The other curve $\alpha=\calH_-(0,0,q)=0$ overlaps the bifurcation curve $\alpha=\calB(0,0)=0$.

\paragraph{Case $\beta\neq0$, $\kap=0$ (Fig.~\ref{f:hexkap0betn0})} The hex-stability boundaries are given by
\begin{align}\label{e:stabbndaniso}
\alpha = \calH_\pm(0,\beta,q) := -\frac{7}{4}\rho_\beta\beta^2 -\frac{1}{3\hK}\left(2q^2 \pm \sqrt{4q^4 + 9\hK q^2\rho_\beta\beta^2}\right).
\end{align}
There exist two turning points $(\pm q_\tp,\alpha_\tp)$ given by
\begin{align}\label{e:qtp}
	q_\tp = \frac{3}{2}|\beta|\sqrt{-\hK\rho_\beta},\quad \alpha_\tp = -\frac{1}{4}\rho_\beta\beta^2.
\end{align}
The boundaries below the turning points are given by  $\calH_-$ which decreases to zero for increasing $|q|$. Hence there exists a hex-stable region near the bifurcation. In particular, for $|q|<q_\tp$, the stripes are hex-stable for all $\alpha$. These indicate that the advection stabilises the stripes: for $\beta\neq 0$ stripes bifurcate stably in accordance with Remark~\ref{r:1Dstab}, and advection connects the regions of stable small and larger amplitude stripes for small quadratic effects. Nevertheless, the hex-unstable region becomes larger for larger $|q|$, which highlights the destabilising effect of the quadratic term.

\begin{remark}
	For fixed $|q|>q_\tp$, the stable stripes lose the stability when $\alpha$ increases to $\alpha_*$ where $\alpha_*<\alpha_\tp$. In fact, at $\alpha=\alpha_\tp$ the homogeneous steady state becomes unstable against hexagonal modes, cf.~\eqref{e:evlinearop}, also see Fig.~\ref{f:homspec} (green curve). Hence, the stable stripes lose stability `before' the bifurcation of hexagons.
\end{remark}

\subsubsection{Stripes with off-critical wavenumber}\label{s:hexsidemodes}
Now we turn to hex-stability of stripes with off-critical wavenumber $\kappa=\kc+\kap$, $\kap\neq0$. We also compare the hex-instability with Eckhaus instability in $(\kap,\alpha)$-plane. Recall that the stripes are zigzag unstable (stable) for $\kap<0$ ($\kap>0$).

\medskip
In the $(q,\alpha)$-plane, for any fixed $\beta$, the stability boundaries are shifted upwards compared with $\kap=0$, cf.\ Fig.~\ref{f:hexkapn0bet0} \& \ref{f:hexkapn0betn0}.

\begin{figure}[!t]
\centering
\subfloat[$\beta=0$, $q=0$]{\includegraphics[width=0.33\linewidth]{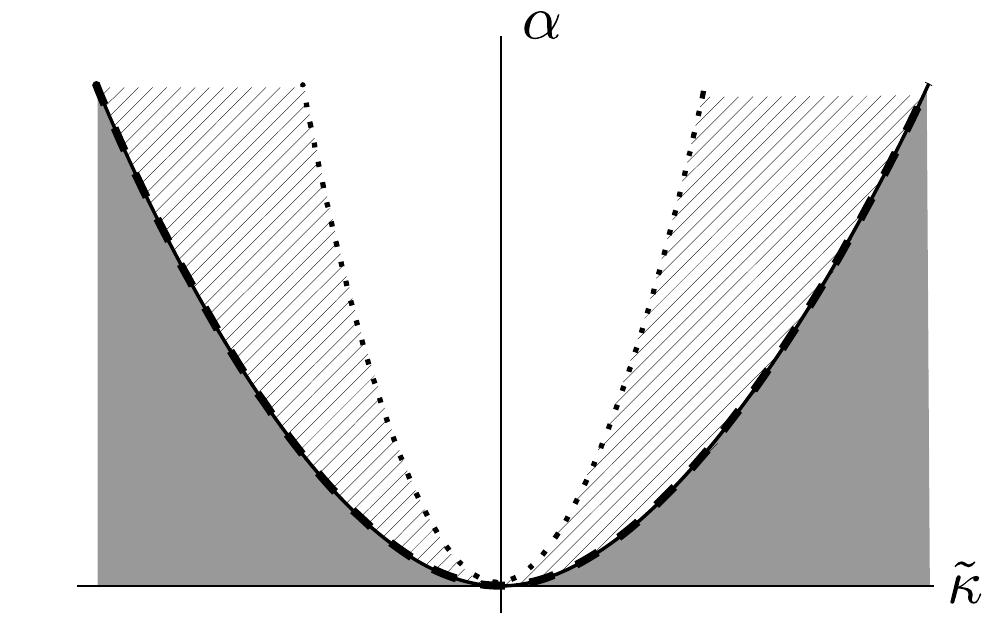}\label{f:hexEckbet0q0}}
\hfil
\subfloat[$\beta=0$, $q\neq0$]{\includegraphics[width=0.33\linewidth]{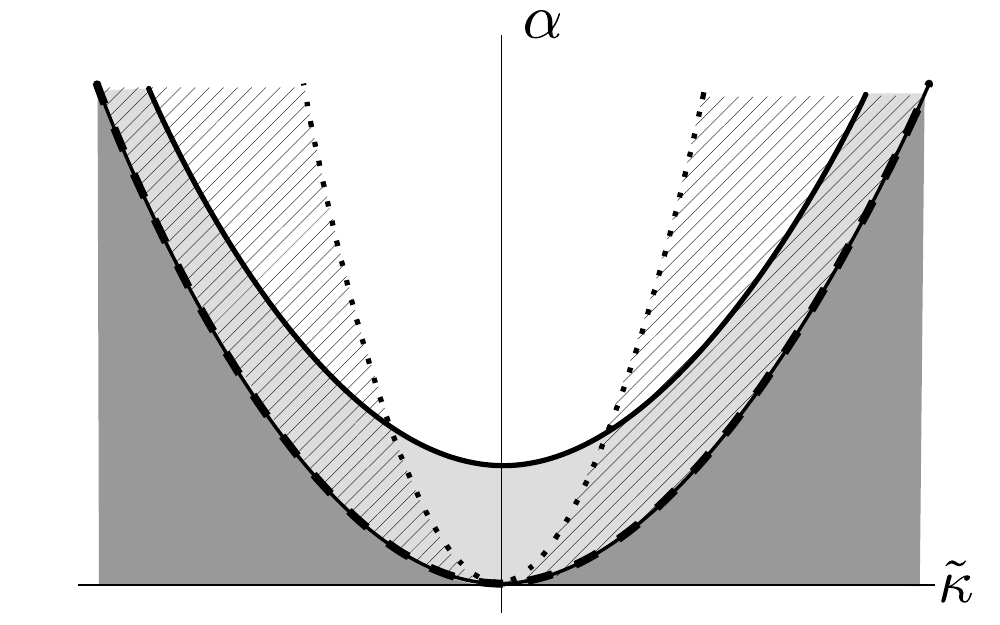}\label{f:hexEckbet0qn0}}
\hfil
\subfloat[$\beta\neq0$, $q=0$]{\includegraphics[width=0.33\linewidth]{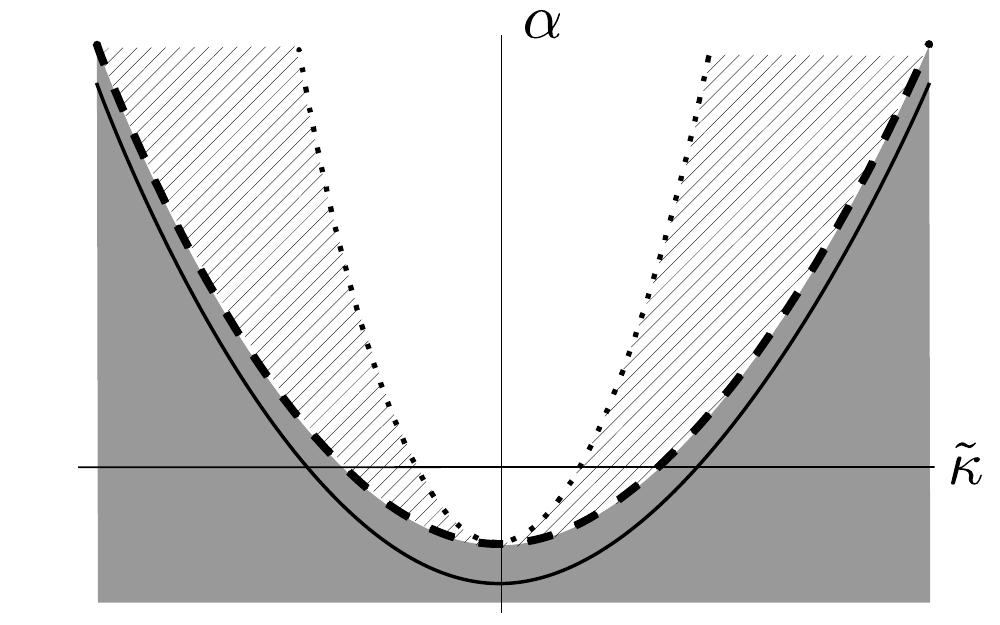}\label{f:hexEckbetn0q0}}
\hfil
\subfloat[$0<|\beta|<\beta_\tp$, $q\neq0$]{\includegraphics[width=0.33\linewidth]{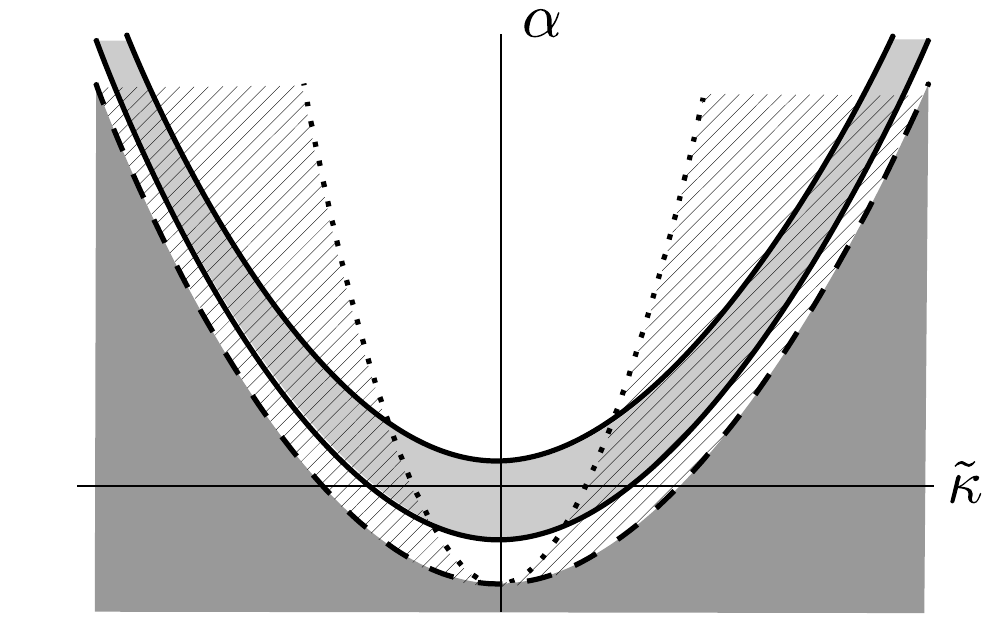}\label{f:hexEckbetn0qn0gr}}
\hfil
\subfloat[$|\beta|=\beta_\tp$, $q\neq0$]{\includegraphics[width=0.33\linewidth]{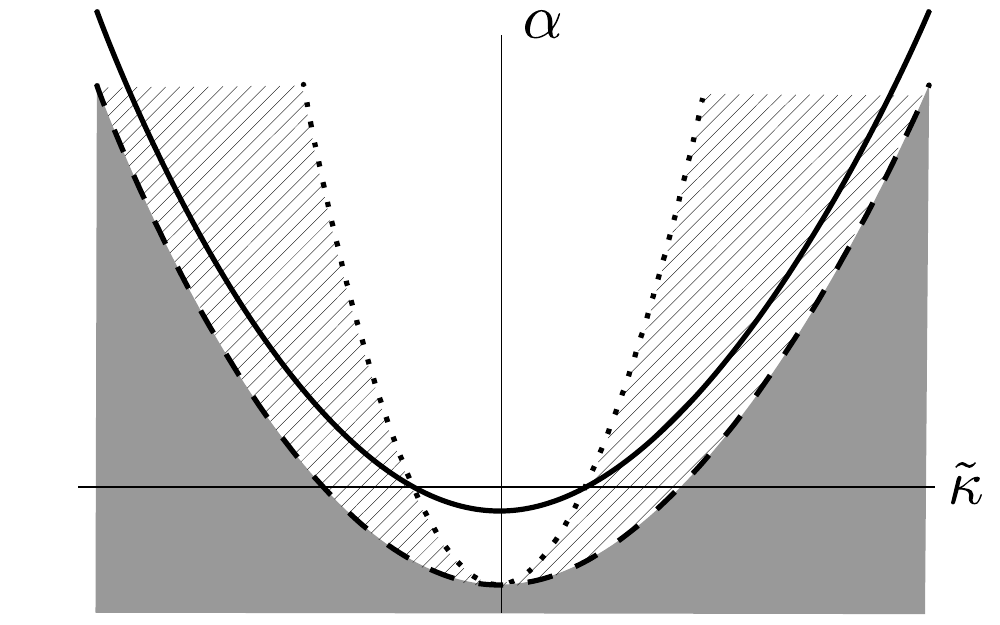}\label{f:hexEckbetn0qn0eq}}
\hfil
\subfloat[$|\beta|>\beta_\tp$, $q\neq0$]{\includegraphics[width=0.33\linewidth]{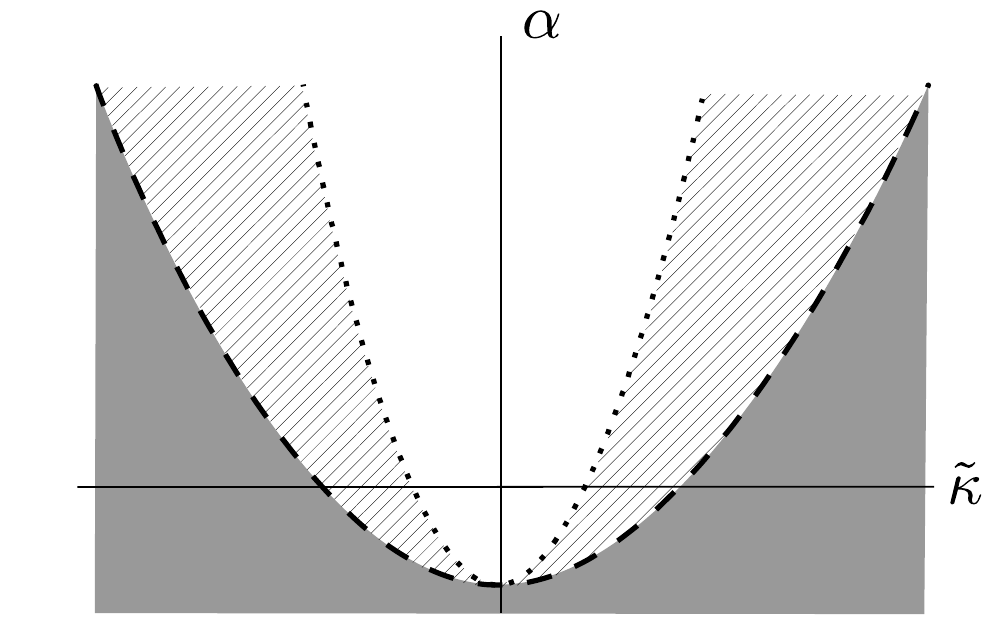}\label{f:hexEckbetn0qn0le}}
	\caption{Sketches of the stability regions in the $(\kap,\alpha)$-plane. Stripes exist in the complement of the dark grey regions; light grey: hex-unstable; hatched regions: Eckhaus-unstable; white: hex-stable. Bifurcation $\alpha=\calB(\kap,\beta)$ (dashed); Eckhaus boundary $\alpha=\calE(\kap,\beta)$ (dotted); hex-stability boundaries (solid) in (a) $\alpha=\calH_\pm(\kap,0,0) = \calB(\kap,0)$, cf. \eqref{e:Eckhausstabbndq0} (b) $\alpha=\calH_+(\kap,0,q)>\calB(\kap,0)$ and $\alpha=\calH_-(\kap,0,q) = \calB(\kap,0)$, cf. \eqref{e:hexbndbet0qn0large} (c) $\alpha=\calH_\pm(\kap,\beta,0)<\calB(\kap,\beta)$, cf. \eqref{e:hexbetn0q0} (d) $\alpha=\calH_\pm(\kap,\beta,q)>\calB(\kap,\beta)$, cf. \eqref{e:hexbnd} (e) $\alpha=\calH_\pm(\kap,\beta_\tp,q)$, (f) no solutions to \eqref{e:hexstabbnd} exist.}
\label{f:hexkapalp}
\end{figure}

\medskip
In the $(\kap,\alpha)$-plane the situation is more involved and can be compared with the Eckhaus instability. In Fig.~\ref{f:hexkapalp} we plot all cases in terms of $\beta$ and $q$, and derive these next.

\paragraph{Case $\beta=0$, $q=0$ (Fig.~\ref{f:hexEckbet0q0})} The hex-stability boundary is given by
\begin{align}\label{e:Eckhausstabbndq0}
\alpha = \calH_\pm(\kap,0,0) = -\rho_\kap\kap^2,
\end{align}
which coincides with the bifurcation curve since $\calH_\pm(\kap,0,0) = \calB(\kap,0)$. Hence the stripes are hex-stable, and the dominant instability mechanism is the Eckhaus boundary.

\paragraph{Case $\beta=0$, $q\neq0$ (Fig.~\ref{f:hexEckbet0qn0})} The hex-stability boundaries are given by
\begin{equation}\label{e:hexbndbet0qn0large}
\begin{aligned}
\alpha &= \calH_+(\kap,0,q) = -\frac{4}{3\hK}q^2 - \rho_\kap\kap^2,\\
\alpha &= \calH_-(\kap,0,q) = -\rho_\kap\kap^2,
\end{aligned}
\end{equation}
where $\calH_-(\kap,0,q) = \calB(\kap,0)$. Hence the stripes are hex-unstable near the bifurcation, which is thus the dominant mechanism near onset. In addition, the curvature of each of the hex-stability boundaries is smaller than that of Eckhaus boundary since $\partial_\kap^2\calH_\pm < \partial_\kap^2\calE$.

\paragraph{Case $\beta\neq0$, $q=0$ (Fig.~\ref{f:hexEckbetn0q0})} The hex-stability boundary is given by
\begin{align}\label{e:hexbetn0q0}
\alpha = \calH_\pm(\kap,\beta,0) = -\frac{7}{4}\rho_\beta\beta^2 - \rho_\kap\kap^2.
\end{align}
Since $\calH_\pm(\kap,\beta,0)<\calB(\kap,\beta)$, the bifurcating stripes are always hex-stable, and the Eckhaus instability is dominant, again in accordance with Remark~\ref{r:1Dstab}.

\paragraph{Case $\beta\neq0$, $q\neq0$ (bottom row of Fig.~\ref{f:hexkapalp})} The hex-stability boundaries are given by \eqref{e:hexbnd}, and roots of the discriminant $\discrh=0$ from \eqref{e:discrh}, lie at
\begin{align}\label{e:betatp}
\beta=\beta_\tp:=\frac{2|q|}{3\sqrt{-\hK\rho_\beta}}.
\end{align}
We summarise the stability results in terms of $\beta$ for fixed $q\neq0$ as follows.

\begin{itemize}
\item[(1)] $|\beta|<\beta_\tp$ (Fig.~\ref{f:hexEckbetn0qn0gr}): hex-stability boundaries satisfy $\calH_\pm(\kap,\beta,q)>\calB(\kap,\beta)$ so that stripes are hex-stable near onset, but there is a hex-unstable `band' which intersects the $\alpha$-axis on the interval $[\calH_-(0,\beta,q),\calH_+(0,\beta,q)]$.
\item[(2)] $|\beta|=\beta_\tp$ (Fig.~\ref{f:hexEckbetn0qn0eq}): The hex-stability boundaries collapse along
\[
\alpha = \calH_\pm(\kap,\beta_\tp,q) =\frac{q^2}{9\hK}-\rho_\kap\kap^2,
\]
which intersects $\alpha$-axis at $\alpha_\tp = -\rho_\beta\beta^2_\tp/4 = q^2/(9\hK)$, cf.\ \eqref{e:qtp}. Notably, this degenerate case does not occur for quasi-hexagonal lattices discussed below.
\item[(3)] $|\beta|>\beta_\tp$ (Fig.~\ref{f:hexEckbetn0qn0le}): $\calH_\pm$ are complex, so there is no hex-stability boundary in the real parameter space and the stripes are hex-stable.
\end{itemize}

In addition, we recall the threshold $q_\tp$, cf.~\eqref{e:qtp}, and highlight the relation $\sgn(|\beta|-\beta_\tp) = -\sgn(|q|-q_\tp)$. Therefore, by increasing $|q|$ for fixed $\beta\neq0$ the hexagonal boundaries change as from Fig.~\ref{f:hexEckbetn0qn0le} to \ref{f:hexEckbetn0qn0gr}.

\begin{figure}[!t]
\centering
\includegraphics[width=0.35\linewidth]{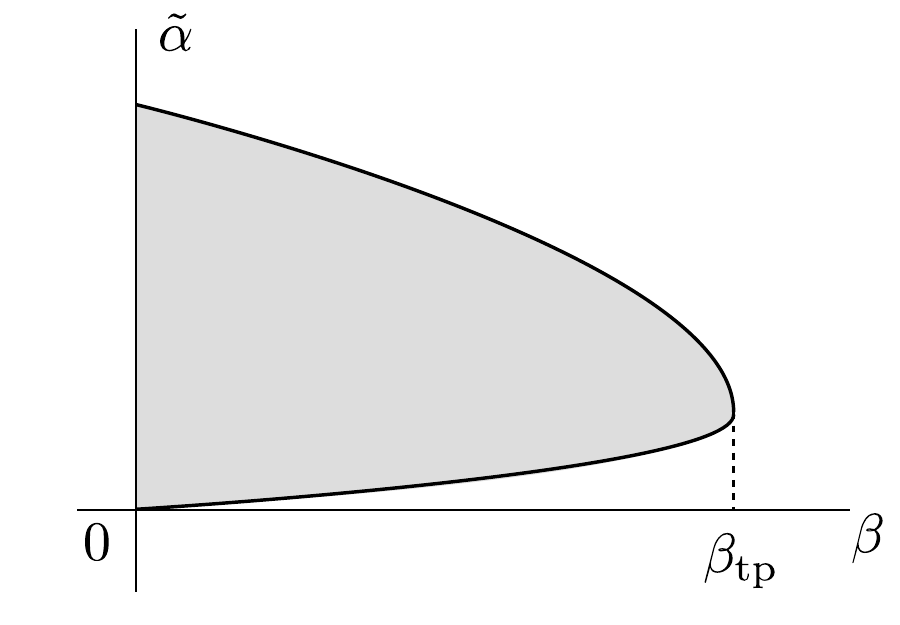}
\caption{Sketch of the hex-unstable region in the $(\beta,\talpha)$-plane for fixed $q\neq0$. Solid curve: hex-boundary $\talpha = \tilde\calH_\pm(\beta)$, cf. \eqref{e:hexbettalpha}; grey: hex-unstable; white: hex-stable.}
\label{f:stabbndalpbeta}
\end{figure}

We consider the width of the unstable band for fixed $q$ in Fig.~\ref{f:hexEckbetn0qn0gr} by setting $\talpha:=\alpha+\rho_\beta\beta^2+\rho_\kap\kap^2$ so that stripe bifurcations occur at $\talpha=0$. Then the hex-stability boundaries in the $(\beta,\talpha)$-plane are
\begin{align}\label{e:hexbettalpha}
\talpha = \widetilde\calH_\pm(\beta) := -\frac{3}{4}\rho_\beta\beta^2 -\frac{1}{3\hK}\left(2q^2 \pm \sqrt{4q^4 + 9\hK q^2\rho_\beta\beta^2}\right),
\end{align}
see Fig.~\ref{f:stabbndalpbeta}. In particular, $\widetilde\calH_-(0) = 0$, $\widetilde\calH_+(0) = -4q^2/(3\hK)$ and $\widetilde\calH_{\pm}(\beta_\tp) = -q^2/(3\hK)$ so the width of hex-unstable band is smaller for larger $|\beta|$, showing the stabilisation of the advection. Note that the width of the unstable band is independent of $\kap$, which will be different for the quasi-hexagonal lattice modes considered next.

\subsection{Stability against quasi-hexagonal perturbations}\label{s:rhomb}
We consider the stability of stripes against \emph{quasi-hexagonal perturbation}, which are nearly hexagonal perturbations that still possess triads $\k_1+\k_2+\k_3=0$ in terms of the spatial scalings.

\medskip
We consider \eqref{e:RDS} with periodic boundary conditions on the rectangular domain 
\[
\Omega_\rhomb:=[0,4\pi/\kappa]\times[0,4\pi/(\sqrt 3 \ell)],\; \kappa:=\kc+\kap,\, \ell:=\kc+\tel,\, \tel\neq\kap,
\]
with the scaling $\tel = \scal\tel'$ so that $\tel = \calO(\scal)$ analogous to $\kap$. Rescaling the spatial variables anisotropically with $\tilde x=x/\kappa$ and $\tilde y=y/\ell$, the rectangular domain becomes $\Omega_3=[0,4\pi]\times[0,4\pi/\sqrt{3}]$ with dual lattice wavevectors $\k_j=(k_j,\ell_j)\in\R^2$, and the perturbation on the six-dimensional kernel is given by $U_\hex(x)$, cf. \S\ref{s:hex}. Analogous to Theorem~\ref{t:cmfhex}, we have the following.
\begin{theorem}\label{t:rhombmatrix}
Consider \eqref{e:RDS} with periodic boundary conditions on rectangular domain $\Omega_\rhomb$. Assume the conditions and notations of Theorem~\ref{t:cmf} for the domain $\Omega_3$ with periodic boundary conditions and the parameter scaling \eqref{e:scaling} for $\mu$. Let the velocity parameter $c=c(\mu)$ be as in \eqref{e:stripev}. The subspace $\{ u_j=0, j=\pm 2,\pm 3\}$ is invariant for the reduced ODE and contains the stripes as equilibria. The linearisation in stripes in the index ordering $(1, -1, 2, -3, 3, -2)$ has a block diagonal matrix of the form $L_\rhomb = \diag(L_1, L_2^\rhomb, L_2^\rhomb)+\calO(\scal^3)$ with $2\times 2$-subblocks

\begin{align*}
L_1 =
A^2\begin{pmatrix}\rho_\nl & \rho_\nl\\ \rho_\nl & \rho_\nl\end{pmatrix}, \quad
L_2^\rhomb =\scal^2
\begin{pmatrix}
\lambda_{\mu,\tel}' + A'^2 \eta & 2A'\frac{q}{\scal} + A' p(\mu_1,\tel')\\
2A'\frac{q}{\scal} + A' \overline{p(\mu_1,\tel')} & \lambda_{\mu,\tel}' + A'^2 \eta
\end{pmatrix}
\end{align*}
where $\eta$ is as in Theorem~\ref{t:cmfhex} and
\begin{align*}
\lambda_{\mu,\tel}' :=&\ \alpha' + \frac{1}{4}\rho_\beta\beta'^2 + \frac{\rho_\kap}{16}(\kap' + 3 \tel')^2 + \calO(\scal),\\
p(\mu_1,\tel') :=&\ \langle \Q[\rmi\beta'w_{A\beta} + (\tfrac{5}{2}\kap' + \tfrac{3}{2}\tel')w_{A\kap},\E_0] -(\rmi \beta'\kc B + (\kap'+3\tel') \kc D )\Q_1,E_0^*\rangle.
\end{align*}
\end{theorem}

\begin{proof}
The rescaled linear operator of \eqref{e:RDS} is given by
\[
\calL_\mu^\rhomb := \kappa^2 D\partial_x^2 + \ell^2 D\partial_y^2 + \A + \calpha\M + \beta\kappa\B\partial_x.
\]
Analogous to the proof of Theorem~\ref{t:cmfhex}, the linearisation in stripes gives the same matrix $L_1$ since the rescaling in $y$-direction does not influence the one-dimensional stability. The matrix $L_2^\rhomb$, however, is different from $L_2^\hex$. The eigenvalue $\lambda_{\mu,\tel}'$ is that of the linearisation in the trivial equilibrium whose expansion can be determined  analogous to Lemma~\ref{l:Turbeta}. The term $p(\mu_1)$ is replaced by $p(\mu_1,\tel')$ by straightforward calculation, which is analogous to the proof in Appendix~\ref{s:hexmat}.
\end{proof}

Concerning the subblock $L_2^\rhomb$, we first note the general form of eigenvalues.
\begin{lemma}\label{l:rhombevals}
Under the assumptions of Theorem~\ref{t:rhombmatrix}, the eigenvalues of the matrix $L_2^\rhomb$ are
\[
\lambda_\pm = \scal^2\left(A'^2(3\hK - q_2 + 8 q_1) - \frac 3 4 \rho_\beta \beta'^2 + \omega' \pm A'\left|\frac{2q}{\scal}+A'p(\mu_1,\tel')\right|\right) + \calO(\scal^3),
\]
where $A'=\sqrt{-(\alpha'+\rho_\beta \beta'^2 + \rho_\kap \kap'^2)/\rho_\nl} + \calO(\scal)$ and $\omega' := (9\tel'+15\kap')(\tel'-\kap')\rho_\kap/16$. The most unstable quasi-hexagonal perturbation with respect to $\tel$ occurs at $\tel=-\kap/3$ for which $\omega = -\rho_\kap\kap^2\geq 0$, and $\tel=0$ gives $\omega=0$ and $\lambda_\pm=\lambda_\hex^\pm$.
\end{lemma}
\begin{proof}
The eigenvalues are derived as in Lemma~\ref{l:hexevals}. 
As a function of $\tel$, the parabola $\omega = \omega(\tel)$ has positive maximum $\max_{\tel\in\R}\omega = -\rho_\kap\kap^2$ at $\tel=-\kap/3$.
\end{proof}

\begin{remark}
We note a relation of the most unstable quasi-hexagonal modes at $\tel=-\kap/3$ to the critical circle of spectrum $S_\kc$ at the onset of the Turing instability. Indeed, it follows from $(\frac{1}{2}(\kc+\kap))^2+(\frac{\sqrt{3}}{2}(\kc+\tel))^2=\kcsq$ that $\tel=-\frac{1}{3}\kap - \frac{2}{9\kc}\kap^2+\calO(\kap^3)$. Therefore, the locations of the most unstable oblique wavevectors are to leading order on the critical circle $S_\kc$. 
\end{remark}

In the remainder of this section, we focus on the quasi-hexagonal perturbation that are more unstable than the hexagonal ones, i.e. in case $\omega>0$, and parametrise $\omega\in(0,-\rho_\kap \kap^2]$ by $\theta\in(0,1]$ via
\[
\omega = -\theta\rho_\kap\kap^2,
\]
so $\theta=1$ is the most unstable quasi-hexagonal perturbation and the limit $\theta=0$ yields the hexagonal one. The previous lemma shows that as for hexagonal perturbations, a smallness assumption on $q$ is required, and as in Theorem~\ref{t:evhex} this changes $A'$ in Lemma~\ref{l:rhombevals} to where $\tA':=\sqrt{-(\alpha'+\rho_\beta\beta'^2+\rho_\kap\kap'^2)/(3\hK)}$. However, unlike the hexagonal stability, for $\Q=o(\scal)$ the stripes are not necessarily stable against quasi-hexagonal perturbations. The previous lemma then directly gives

\begin{theorem}
Under the assumptions of Theorem~\ref{t:rhombmatrix} and $\theta\in(0,1]$ the quasi-hexagonal stability boundary, i.e., zero real part of the eigenvalues of the matrix $L_2^\rhomb$ is to leading order given as follows.
 
If $\Q=g(\scal)\Q'$, $g(\scal)=o(\scal)$ this stability boundary reads
\begin{align}\label{e:rhombhigh}
\alpha =\calM_\rhomb(\kap,\beta;\theta) := -\frac{7}{4}\rho_\beta\beta^2 -(\theta+1) \rho_\kap\kap^2.
\end{align}
Under Hypothesis~\ref{h:Qscale} this stability boundary is given by the two curves
\begin{align}
\{\alpha=\calM_\Rhomb^\pm(\kap,\beta,q;\theta): \alpha\geq\calM_\rhomb(\kap,\beta;\theta),\; \discr\geq 0\},\label{e:rhombbnd1}\\
\{\alpha=\calM_\Rhomb^\pm(\kap,\beta,q;\theta): \alpha\leq\calM_\rhomb(\kap,\beta;\theta),\; \discr\geq 0\},\label{e:rhombbnd2}
\end{align}
corresponding to the two eigenvalues $\lambda_\sigma$ with possibly different $\sigma=\pm$, where
\begin{align}
\calM_\Rhomb^\pm(\kap,\beta,q;\theta)&:= -\frac{7}{4}\rho_\beta\beta^2- (\theta+1)\rho_\kap\kap^2\nonumber
 -\frac{1}{3\hK}\left(2q^2   \pm\sqrt{\discr}\right),\\\label{e:discr}
 \discr &:= 4q^4 + 9\hK\rho_\beta\beta^2 q^2 + 12\hK\theta\rho_\kap\kap^2q^2.
\end{align}
\end{theorem}

\subsubsection{Isotropic case $\beta=0$} 
For $\Q=\calO(\scal)$, the quasi-hex-stability boundary is given to leading order by the following two parts, see Fig.~\ref{f:rhombetn0kapn0le}.
\begin{align}
\alpha &= \calM_\Rhomb^+(\kap,0,q;\theta) = -(\theta+1)\rho_\kap\kap^2 -\frac{2}{3\hK}\left(q^2 + \sqrt{q^4  + 3\hK\theta\rho_\kap\kap^2q^2}\right),\label{e:bndrhomiso}\\
\alpha &= \calM_\Rhomb^-(\kap,0,q;\theta) = -(\theta+1)\rho_\kap\kap^2 -\frac{2}{3\hK}\left(q^2 - \sqrt{q^4  + 3\hK\theta\rho_\kap\kap^2q^2}\right).\label{e:bndrhomiso2}
\end{align}

In the $(q,\alpha)$-plane, the boundary $\alpha=\calM_\Rhomb^+(\kap,0,q;\theta)$ intersects the $\alpha$-axis at $\calM_\Rhomb^+(\kap,0,0;\theta)$, where $\calM_\Rhomb^+(\kap,0,0;\theta)=-(\theta+1)\rho_\kap\kap^2$. Since $\theta>0$, we have $\calM_\Rhomb^+(\kap,0,q;\theta)>\calB(\kap,0)$. Thus the stripes are quasi-hex-unstable near onset. 
Note that for $\Q=o(\scal)$, the stability boundary \eqref{e:rhombhigh} is independent of $q$.

\medskip
In the $(\kap,\alpha)$-plane, the following cases for the quasi-hex-stability boundary occur.

\begin{figure}[!t]
\centering
\subfloat[$\beta=0,q=0$]{
\includegraphics[width=0.33\linewidth]{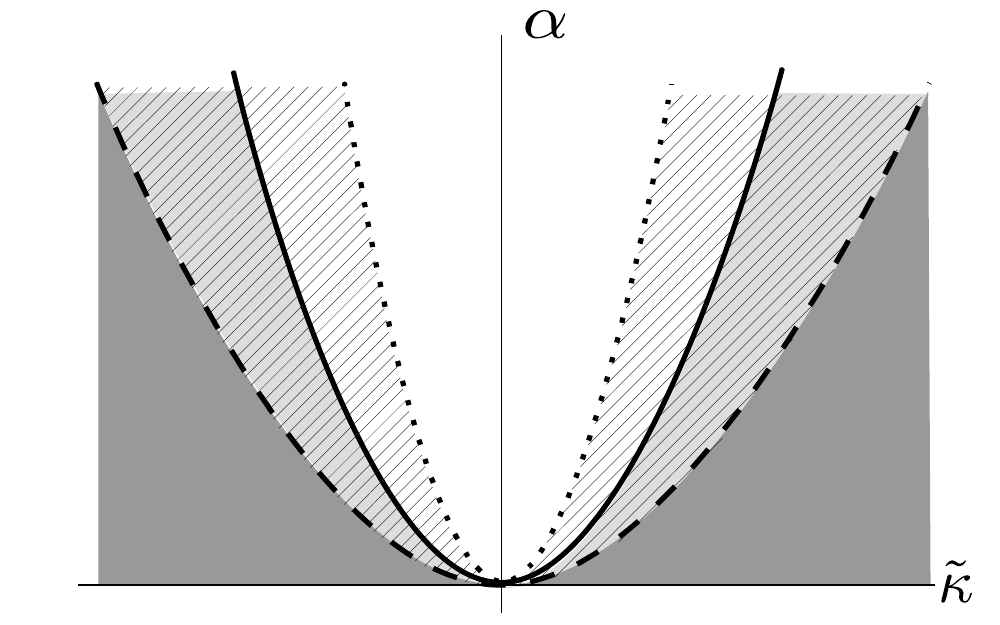}\label{f:rhombet0q0}}
\hfil
\subfloat[$\beta=0,q\neq0$]{\includegraphics[width=0.33\linewidth]{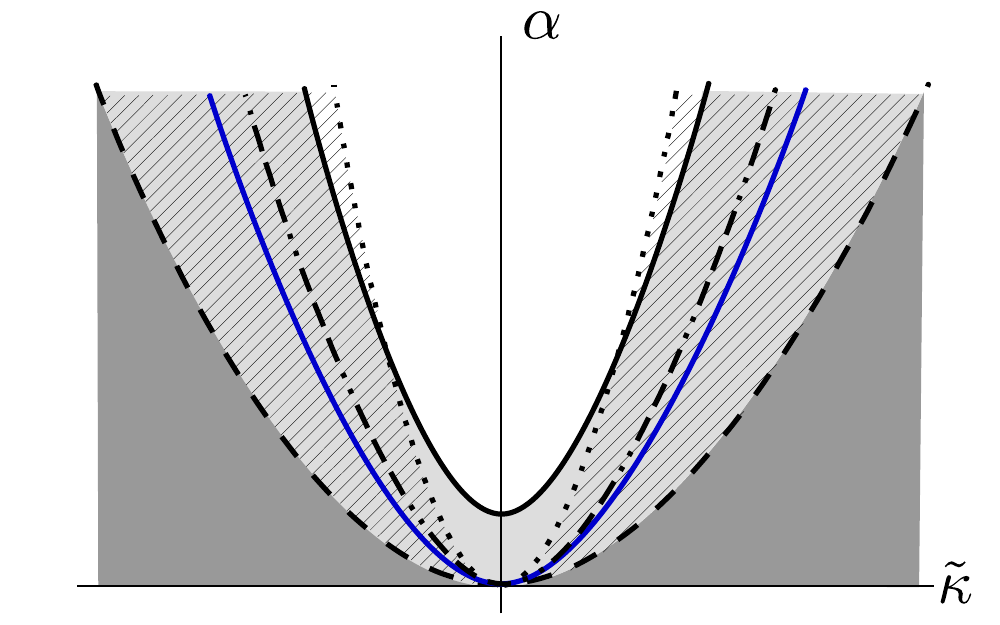}\label{f:rhombet0qn0vthng1}}
\hfil
\subfloat[$\beta\neq0,q=0$]{\includegraphics[width=0.33\linewidth]{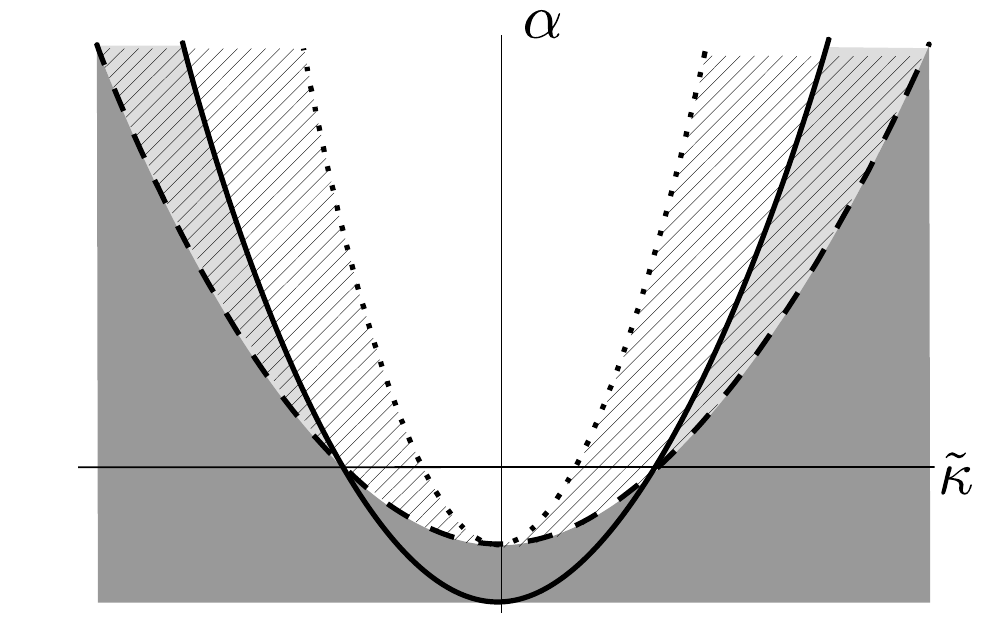}\label{f:rhombetn0q0}}
\caption{Sketches of the stability regions in the $(\kap,\alpha)$-plane for $\theta\in(0,1]$. Stripes exist in the complements of the dark grey regions. Light grey: quasi-hex-unstable; hatched region: Eckhaus-unstable; white: stable; dashed curve: bifurcation curve $\alpha=\calB(\kap,0)$; dotted curve: Eckhaus boundary $\alpha = \calE(\kap,0)$. (a) Quasi-hex-boundary \eqref{e:bndrhomq0} (black solid). (b) Quasi-hex-boundary for $\Q=\calO(\scal)$ \eqref{e:bndrhomiso} (black solid), \eqref{e:bndrhomiso2} (blue solid) and for $\Q=o(\scal)$ \eqref{e:rhombhigh} (dotted dashed). (c) Quasi-hex-boundary \eqref{e:bndrhombetn0q0} (black solid).}
\end{figure}

\paragraph{Case $\beta=0$, $q=0$ (Fig.~\ref{f:rhombet0q0})} The quasi-hex-stability boundaries are independent of $q$. Hence for both $\Q=\calO(\scal)$ and $\Q=o(\scal)$, the quasi-hex-stability boundaries read
\begin{align}\label{e:bndrhomq0}
\alpha =\calM_\Rhomb^\pm(\kap,0,0;\theta) = -(\theta+1)\rho_\kap\kap^2 = \calM_\rhomb(\kap,0;\theta).
\end{align}
Note that since $\theta\in(0,1]$, we have $\calB(\kap,0)\leq\calM_\Rhomb^\pm(\kap,0,0;\theta) = \calM_\rhomb(\kap,0;\theta)\leq\calE(\kap,0)$.

\paragraph{Case $\beta=0$, $q\neq0$ (Fig.~\ref{f:rhombet0qn0vthng1})} For $\Q=\calO(\scal)$, the quasi-hex-stability boundary is on the one hand given by \eqref{e:rhombbnd1} as \eqref{e:bndrhomiso} at $\beta=0$, which intersects the $\alpha$-axis at $\calM_\Rhomb^+(0,0,q;\theta)=-4q^2/(3\hK)$. The ordinate of intersections of Eckhaus boundary and quasi-hex-stability boundary is given by $\alpha_{\rm sec} := -\frac{8q^2}{\hK(2-\theta)^2} = \calO(\scal^2)$. Hence, the quasi-hexagonal instability is the dominant instability mechanism near onset. On the other hand, \eqref{e:rhombbnd2} gives as \eqref{e:bndrhomiso2} at $\beta=0$ a quasi-hex-stability boundary which passes through the origin with curvature larger than that of Eckhaus boundary. 

For $\Q=o(\scal)$, the quasi-hex-stability boundary is given by $\alpha = \calM_\rhomb(\kap,0;\theta)$. We note that $\calM_\Rhomb^-(\kap,0,q;\theta)\leq \calM_\rhomb(\kap,0;\theta) <\calM_\Rhomb^+(\kap,0,q;\theta)$.

\subsubsection{Anisotropic case $\beta\neq0$} 
We first consider the quasi-hex-stability in the $(q,\alpha)$-plane. Since $\Q=o(\scal)$ the quasi-hex-stability boundary is independent of $q$, we omit this case.

For $\Q=\calO(\scal)$, roots of $\discr=0$ from \eqref{e:discr} occur as a function of $q$ precisely when $9\rho_\beta\beta^2 + 12\theta\rho_\kap\kap^2 \geq 0$ so that the threshold in terms of $\beta$ lies at
\begin{align}\label{e:bifrhombep}
\beta = \beta_\ep:=2|\kap|\sqrt{-\theta\rho_\kap/(3\rho_\beta)}.
\end{align}
Notably, $\beta_\ep=0$ for the hexagonal modes, i.e., $\theta=0$, which is consistent with Fig.~\ref{f:hexalpq}.

\medskip
We summarise the quasi-hex-stability boundaries in the $(q,\alpha)$-plane as follows.

\begin{itemize}
\item[(1)] $|\beta|<\beta_\ep$ (Fig.~\ref{f:rhombetn0kapn0le}):  $\alpha = \calM_\Rhomb^+(\kap,\beta,q;\theta)$ has minimum at $\kap=0$ where $\calM_\Rhomb^+(\kap,\beta,0;\theta)> \calB(\kap,\beta)$. Thus the stripes are unstable near the bifurcation.

\medskip
\item[(2)] $|\beta|=\beta_\ep$ (Fig.~\ref{f:rhombetn0kapn0eq}): $\alpha = \calM_\Rhomb^+(\kap,\beta_\ep,q;\theta)$ is a parabola in $q$ which touches the bifurcation line at $q=0$ since $\calM_\Rhomb^+(\kap,\beta_\ep,0;\theta) = \calB(\kap,\beta_\ep)$.

\medskip
\item[(3)] $|\beta|>\beta_\ep$ (Fig.~\ref{f:rhombetn0kapn0gr}): $\alpha=\calM_\Rhomb^\pm(\kap,\beta,q;\theta)$ connected stable region: the stripes are unstable for $\alpha\in(\calM_\Rhomb^-(\kap,\beta,q;\theta),\calM_\Rhomb^+(\kap,\beta,q;\theta))$ and stable elsewhere; in particular the stripes are stable near onset and there are two turning points given by $(\pm q_{\tp,\theta},\alpha_{\tp,\theta})$ where 
\begin{align}\label{e:qtpthe}
q_{\tp,\theta}:=\frac{1}{2}\sqrt{-12\hK\theta\rho_\kap\kap^2-9\hK\rho_\beta\beta^2},\qquad \alpha_{\tp,\theta}:= -\frac{1}{4}\rho_\beta\beta^2 +(\theta-1)\rho_\kap\kap^2.
\end{align}
In particular, the stripes are stable for $|q|<q_{\tp,\theta}$ and all $\alpha$, cf. Fig.~\ref{f:betaq}. The stable regions connect later for larger $|\kap|$ and the connection is wider for larger $|\beta|$.
\end{itemize}

\medskip
Next, we discuss the quasi-hex-stability boundary in the $(\kap,\alpha)$-plane.

\paragraph{Case $\beta\neq0$, $q=0$ (Fig.~\ref{f:rhombetn0q0})} The quasi-hex-stability boundary is independent of $q$. Hence for both $\Q=\calO(\scal)$ and $\Q=o(\scal)$, the quasi-hex-stability boundary is given by
\begin{align}\label{e:bndrhombetn0q0}
\alpha = \calM_\Rhomb^\pm(\kap,\beta,0;\theta) = -\frac{7}{4}\rho_\beta\beta^2 - (\theta+1)\rho_\kap\kap^2 = \calM_\rhomb(\kap,\beta;\theta),
\end{align}
which is a parabola in $\kap$ and is shifted downwards by increasing $|\beta|$. Its curvature is smaller than that of Eckhaus boundary since $\partial_\kap^2\calM_\Rhomb^\pm <\partial_\kap^2\calE$, and thus the Eckhaus instability is dominant.

\begin{figure}[!t]
\centering
\subfloat[$0<|\beta|<\beta_\tp$]{\includegraphics[width=0.35\linewidth]{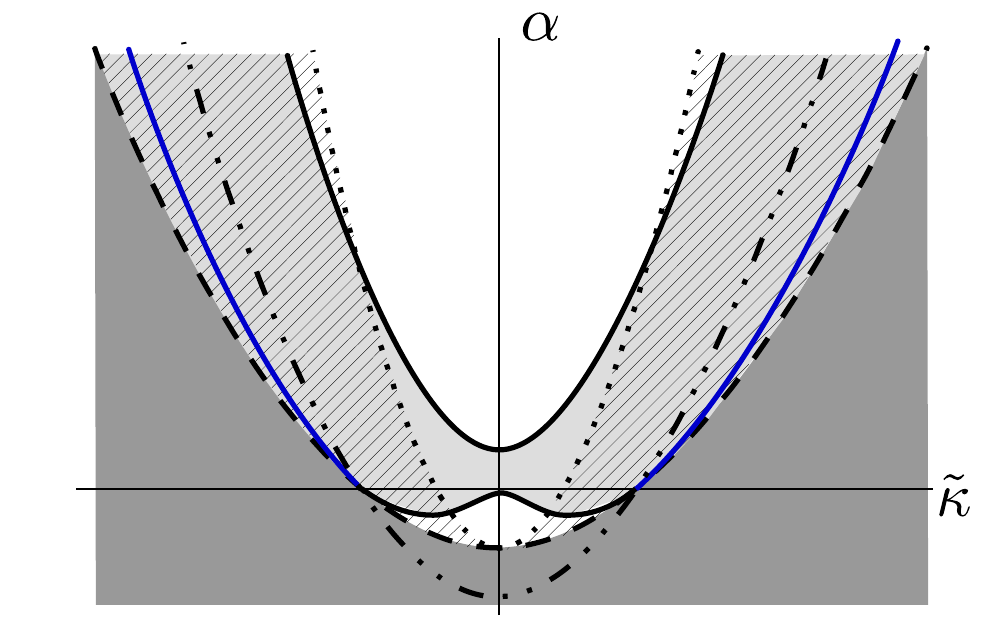}\label{f:Eckrhombetn0qn0gr}}
\hfil
\subfloat[$|\beta|=\beta_\tp$]{\includegraphics[width=0.35\linewidth]{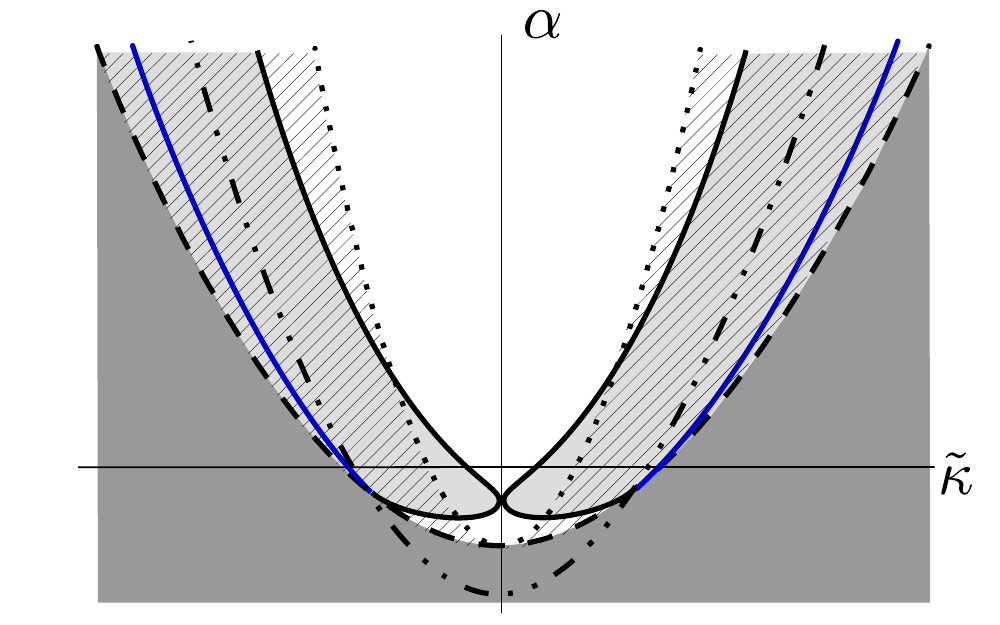}\label{f:Eckrhombetn0qn0eq}}
\hfil
\subfloat[$\beta_\tp<|\beta|<\beta_\ex$]{\includegraphics[width=0.35\linewidth]{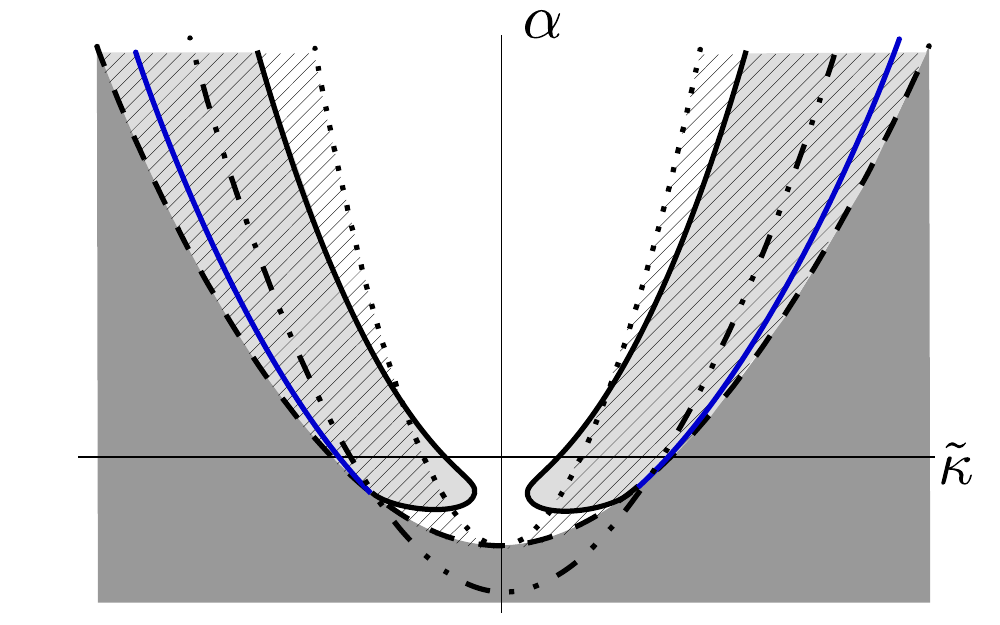}\label{f:Eckrhombetn0qn0rh}}
\hfil
\subfloat[$|\beta|\geq\beta_\ex$]{\includegraphics[width=0.35\linewidth]{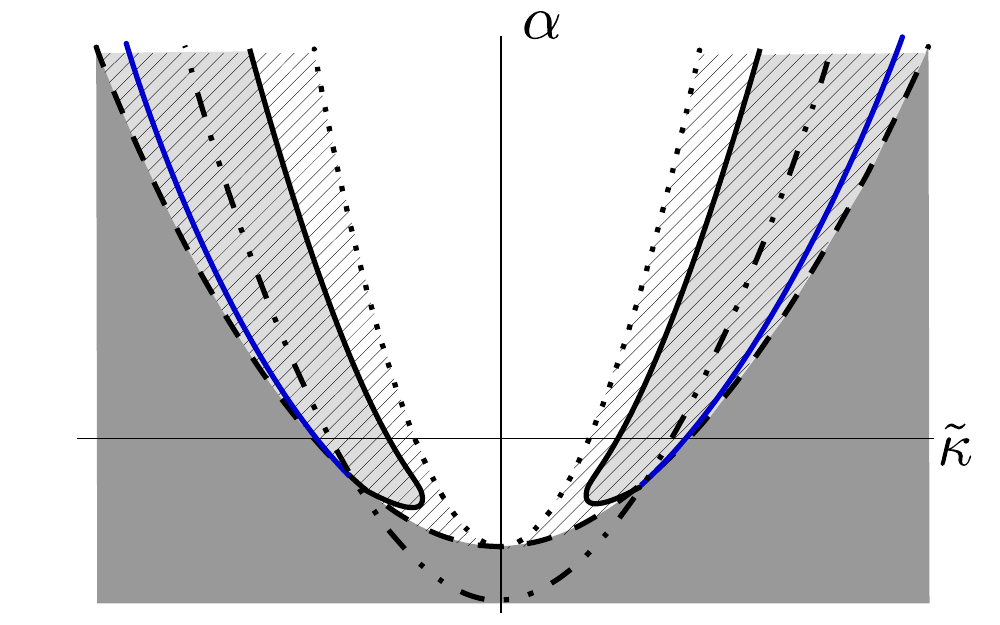}\label{f:Eckrhombetn0qn0le}}
\caption{Sketches of the stability boundaries and Eckhaus boundary $\calE$ for $\beta\neq0$, $q \neq 0$ and $\theta\in(0,1]$. Stripes exist in the complement of the dark grey regions; light grey: quasi-hex-unstable; hatched region: Eckhaus-unstable; white: stable. Dashed curve: bifurcation curve $\alpha=\calB(\kap,\beta)$; dotted curve: $\alpha=\calE(\kap,\beta)$; quasi-hex-boundary for $\Q=\calO(\scal)$ \eqref{e:rhombbnd1} (black solid), \eqref{e:rhombbnd2} (blue solid), quasi-hex-boundary for $\Q=o(\scal)$ \eqref{e:rhombhigh} (dashed-dotted). Zigzag instability occurs for $\kap<0$.}
\label{f:Eckrhombetn0qn0}
\end{figure}

\paragraph{Case $\beta\neq0$, $q\neq0$ (Fig.~\ref{f:Eckrhombetn0qn0})} For $\Q=o(\scal)$, the quasi-hex-stability boundary is given by $\alpha=\calM_\rhomb(\kap,\beta;\theta)$, cf. \eqref{e:rhombhigh}, which is a parabola in $\kap$.

For $\Q=\calO(\scal)$, we recall that the quasi-hex-stability boundaries are given by \eqref{e:rhombbnd1} and \eqref{e:rhombbnd2}, respectively. The boundaries \eqref{e:rhombbnd1} have been shown in Fig.~\ref{f:rhombkapalp2}--\ref{f:rhombkapalp5}. For the completeness of the stability diagrams, however, we replot them in Fig.~\ref{f:Eckrhombetn0qn0}. Solving $\calE(\kap,\beta) =\calM_\Rhomb^\pm(\kap,\beta,q;\theta)$ we find the critical value $\beta_\ex$ such that $\calE$ and $\calM_\Rhomb^\pm$ have only two intersection points for $\beta=\beta_\ex$, where
\begin{align}\label{e:eckrhombex}
\beta_\ex = \frac{2}{3}|q|\sqrt{\frac{2}{\hK(\theta-2)\rho_\beta}}
\end{align}
and $\beta_\ex>\beta_\tp$ where $\beta_\tp$ is given by \eqref{e:betatp}.
This gives the following subcases:

\begin{itemize}
\item[(1)] $|\beta|<\beta_\tp$ (Fig.~\ref{f:Eckrhombetn0qn0gr}): The quasi-hex-stability boundary is given by \eqref{e:rhombbnd1} and composed of two curves. The lower curve touches the bifurcation curve at the endpoints $(\pm\kap_\ep,\alpha_\ep)$ where
\[
\kap_{\ep}= |\beta|\sqrt{-\frac{3\rho_\beta}{4\theta\rho_\kap}}>0,\qquad \alpha_{\ep}:= \left(\frac{3}{4\theta}-1\right)\rho_\beta\beta^2.
\]
In particular, the stripes are quasi-hex-stable near the onset for $|\kap|<\kap_\ep$ only, and these endpoints diverge $\theta\to0$, thus limiting to the hexagonal case, cf. Fig.~\ref{f:hexEckbetn0qn0gr}. In addition, the stability boundaries intersect the $\alpha$-axis at $\calM_\Rhomb^\pm(0,\beta,q;\theta) = \calH_\pm(0,\beta,q)$.
Moreover, the ordinate of intersections of quasi-hex-stability and Eckhaus boundary is given by
\begin{align*}
\alpha_{{\rm sec},\beta}^\pm  =&\ -\frac{1}{4\hK(2-\theta)^2}\bigg(16 q^2+\hK\rho_\beta\beta^2(4\theta^2-25\theta+34) \\
&\quad \pm 4\sqrt{16 q^4+18\hK(2-\theta) q^2\rho_\beta\beta^2}\bigg) = \calO(\scal^2).
\end{align*}
Compared to the isotropic case (cf. Fig.~\ref{f:rhombet0qn0vthng1}), nonzero $\beta$ creates a stable region near the bifurcation and moves the upper boundary downwards, thus the advection stabilises the stripes.
\item[(2)] $|\beta|=\beta_\tp$ (Fig.~\ref{f:Eckrhombetn0qn0eq}): The quasi-hex-stability boundaries intersect the $\alpha$-axis at a single point $\calM_\Rhomb^\pm(0,\beta_\tp,q;\theta)=\calH_\pm(0,\beta_\tp,q) = -\frac{1}{4}\rho_\beta\beta_\tp^2 = q^2/(9\hK)$.
\item[(3)] $\beta_\tp<|\beta|<\beta_\ex$ (Fig.~\ref{f:Eckrhombetn0qn0rh}): The quasi-hex-stability boundary 
consists of two curves whose turning points are given by $(\pm\kap_\mp,\alpha_\mp)$, where 
\begin{align}\label{e:kapmp}
\kap_\mp:= \sqrt{-\frac{3\rho_\beta\beta^2}{4\theta\rho_\kap} - \frac{q^2}{3\hK\theta\rho_\kap}} \ >0,\; \alpha_\mp:= \frac{4q^2(1-\theta)+3\hK(3-4\theta)\rho_\beta\beta^2}{12\hK\theta}.
\end{align}
The stripes are quasi-hex-stable for $|\kap|<\kap_\mp$ and all $\alpha$, cf. Fig.~\ref{f:betakap}. The stable regions connect later for larger $|q|$ and the connection is wider for larger $|\beta|$.
The turning points diverge as $\theta\to0$ and so do the endpoints $(\pm\kap_\ep,\alpha_\ep)$, thus limiting to the hexagonal case, cf. Fig.~\ref{f:hexEckbetn0qn0le}. In contrast to the hexagonal case, here we have two regions where the stripes are quasi-hex-unstable but Eckhaus stable. 

\medskip
\item[(4)] $|\beta|\geq\beta_\ex$ (Fig.~\ref{f:Eckrhombetn0qn0le}): The quasi-hex-stability boundaries touch the Eckhaus boundary for $|\beta|=\beta_\ex$ and lie inside the Eckhaus unstable region.
\end{itemize}

Notably, in each case the Eckhaus instability is dominant near the bifurcation as predicted in Remark~\ref{r:1Dstab}. We recall the threshold $q_\tp$, cf.~\eqref{e:qtp} and have $\sgn(|\beta|-\beta_\tp) = -\sgn(|q|-q_\tp)$, also $\sgn(|\beta|-\beta_\ex) = -\sgn(|q|-q_\ex)$, where 
\[
q_\ex:= \frac{3}{2}|\beta|\sqrt{\hK(2-\theta)\rho_\beta/2},
\]
and $q_\ex<q_\tp$. Therefore, by increasing $|q|$ for fixed $\beta\neq0$ the quasi-hex-stability boundaries change as from Fig.~\ref{f:Eckrhombetn0qn0le} to \ref{f:Eckrhombetn0qn0gr}.

\section{Examples}\label{s:example}

\subsection{Designed example}\label{s:exact}
For illustration of the stabilities we consider the same concrete system as in \cite{Yang2019a}, except the flexible coefficient $\epsilon$ of the quadratic nonlinearity,
\begin{equation}\label{e:example}
	\begin{aligned}
	u_t &= \Delta u + 3u - v + \calpha u + 4\calpha v + \beta u_x + \epsilon u^2+\frac{1}{4}\epsilon v^2- uv^2\\
	v_t &= \frac{7}{2}\Delta v + 14u - \frac{7}{2}v  -\frac{1}{5}\calpha u +\calpha v + \epsilon u^2+\frac{1}{4}\epsilon v^2 + uv^2
	\end{aligned}
\end{equation}
where $U:= (u,v)^T$, $D = \diag(1,7/2)$,
\begin{align*}
	\A = \begin{pmatrix} 3 & -1\\ 14 & -\frac{7}{2} \end{pmatrix}\!,\; \M = \begin{pmatrix} 1 & 4\\ -\frac{1}{5} & 1 \end{pmatrix}\!,\; \Q[U,U] = \epsilon\begin{pmatrix} u^2+\frac{1}{4}v^2 \\ u^2+\frac{1}{4}v^2 \end{pmatrix}\!,\; \K[U,U,U] = \begin{pmatrix} -uv^2 \\ uv^2 \end{pmatrix}\!.
\end{align*}
The generic form of $\Q$ is given by $\Q[U_1,U_2] = (\Q_|[U_1,U_2],\Q_{||}[U_1,U_2])^T$ with
\[
\Q_|[U_1,U_2] = \Q_{||}[U_1,U_2] = \epsilon U_1^T \begin{pmatrix}
1 & 0 \\ 0 & \frac{1}{4}
\end{pmatrix} U_2,
\]
where $U_j:= (u_j,v_j)^T$, $j = 1,2,3$. 

\begin{figure}[!t]
\centering
\subfloat[$\beta=0$]{\includegraphics[width=0.3\linewidth]{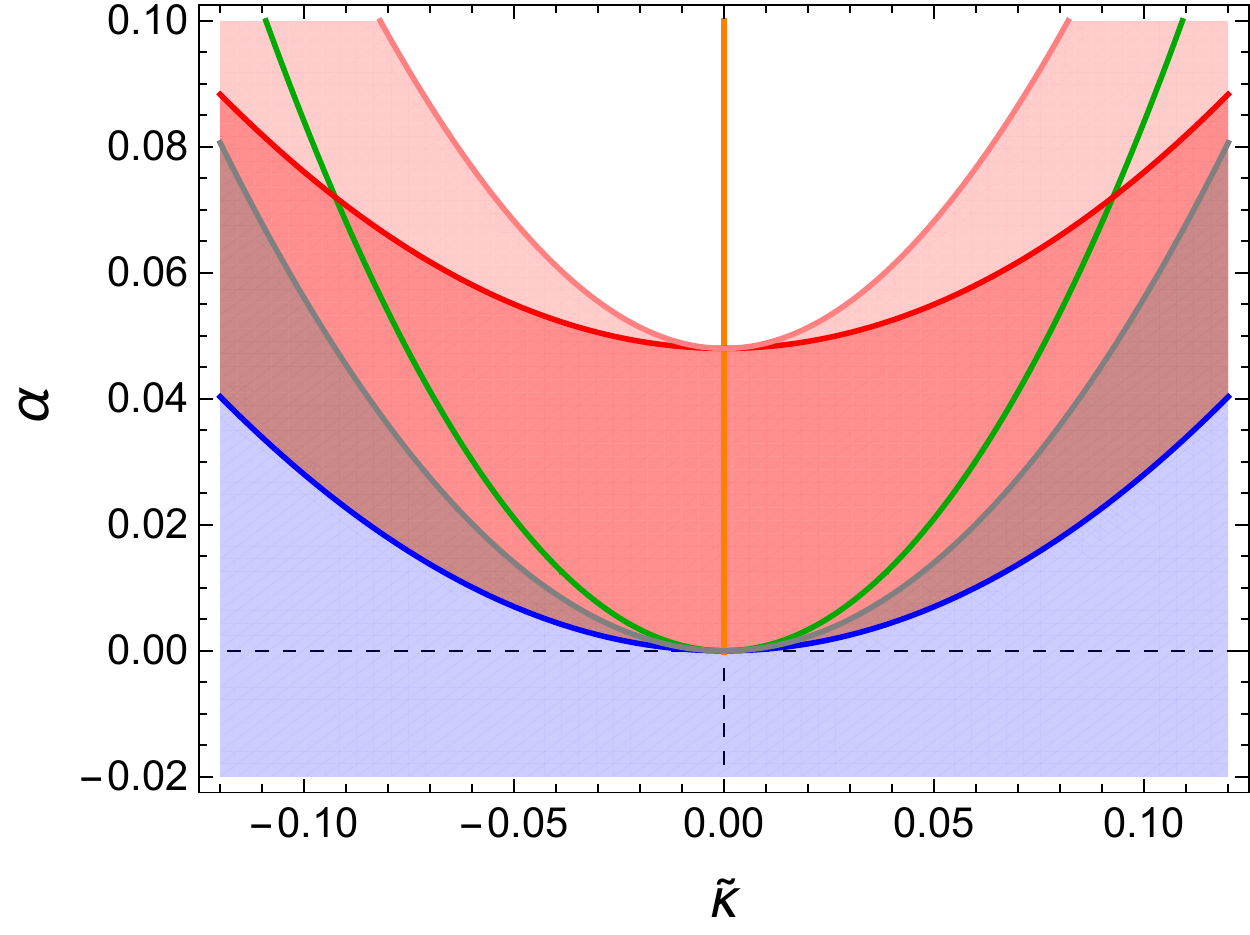}\label{f:egbet0ep0p8}}
\hfil
\subfloat[$\beta=0.36$]{\includegraphics[width=0.3\linewidth]{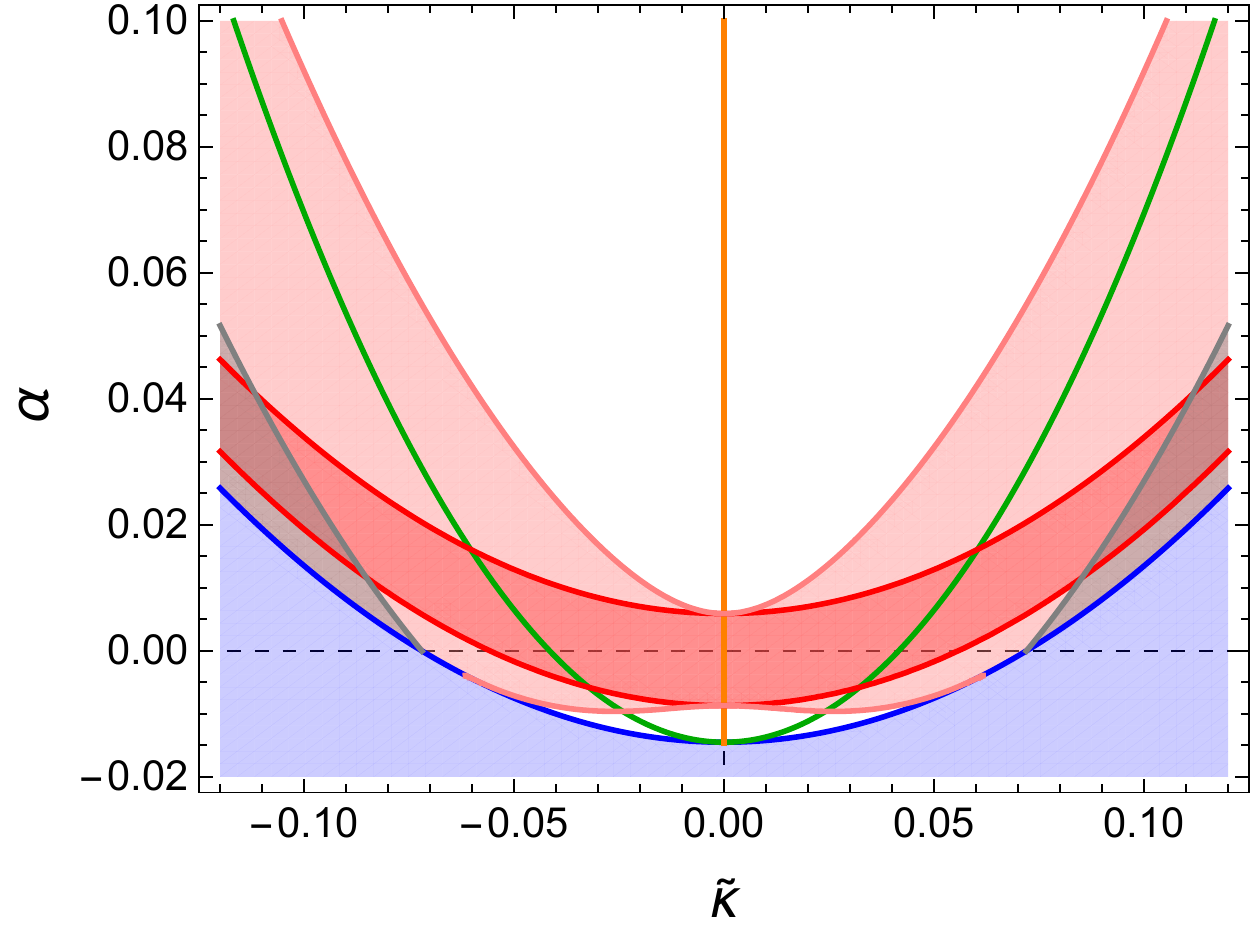}\label{f:egbet0p36ep0p4}}
\hfil
\subfloat[$\beta=0.378$]{\includegraphics[width=0.3\linewidth]{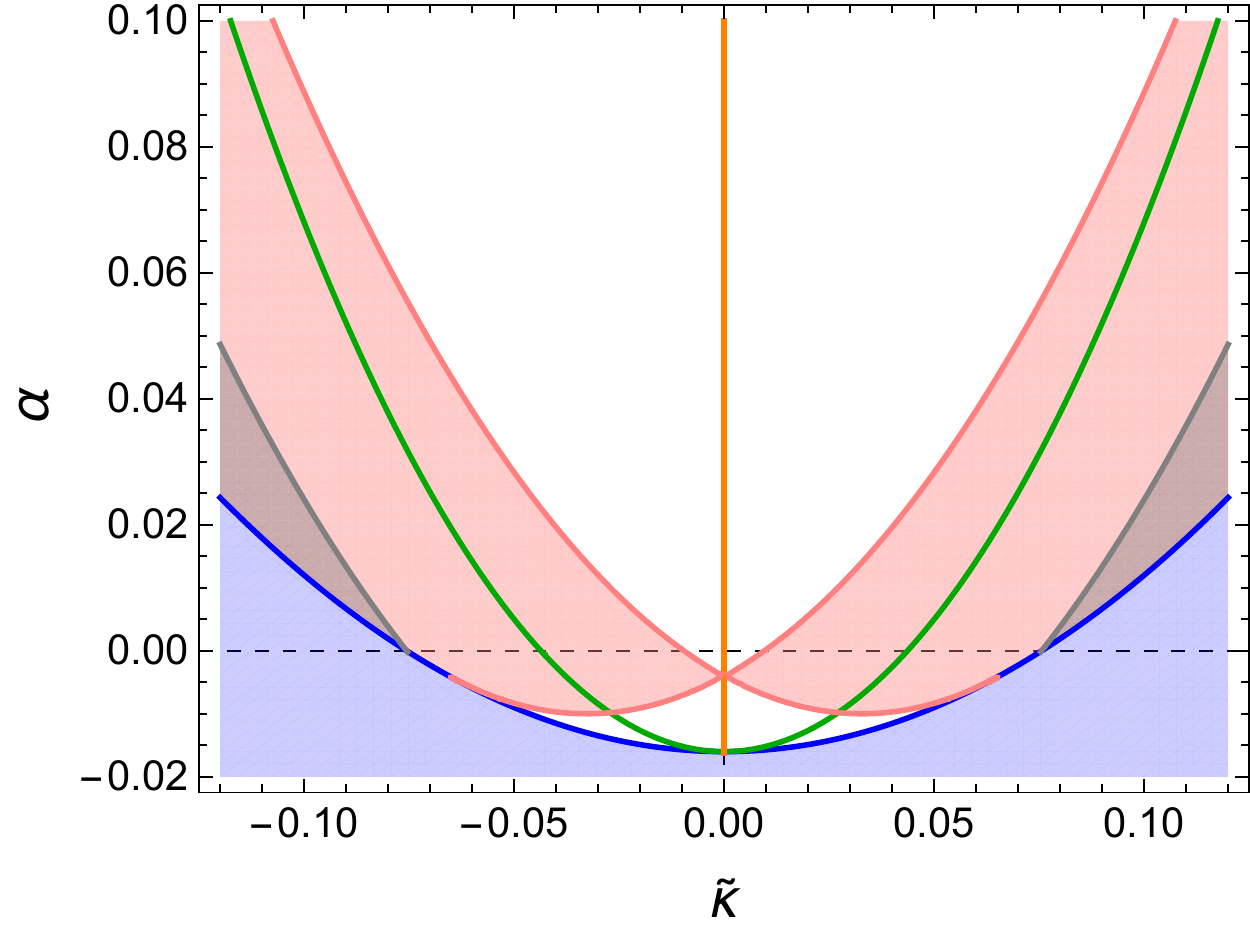}\label{f:egbet0p378ep0p4}}
\hfil
\subfloat[$\beta=0.39$]{\includegraphics[width=0.3\linewidth]{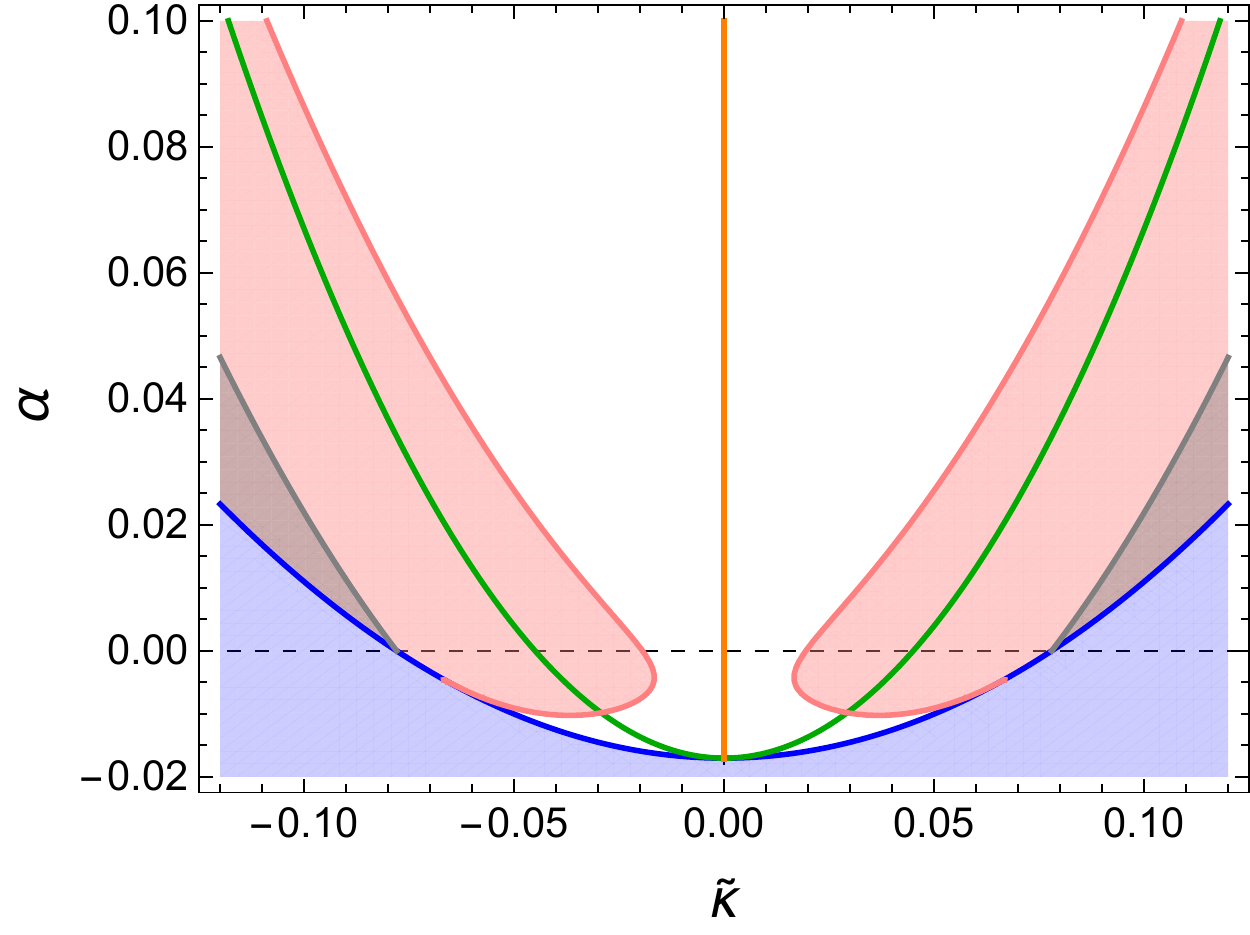}\label{f:egbet0p39ep0p4}}
\hfil
\subfloat[$\beta=0.56$]{\includegraphics[width=0.3\linewidth]{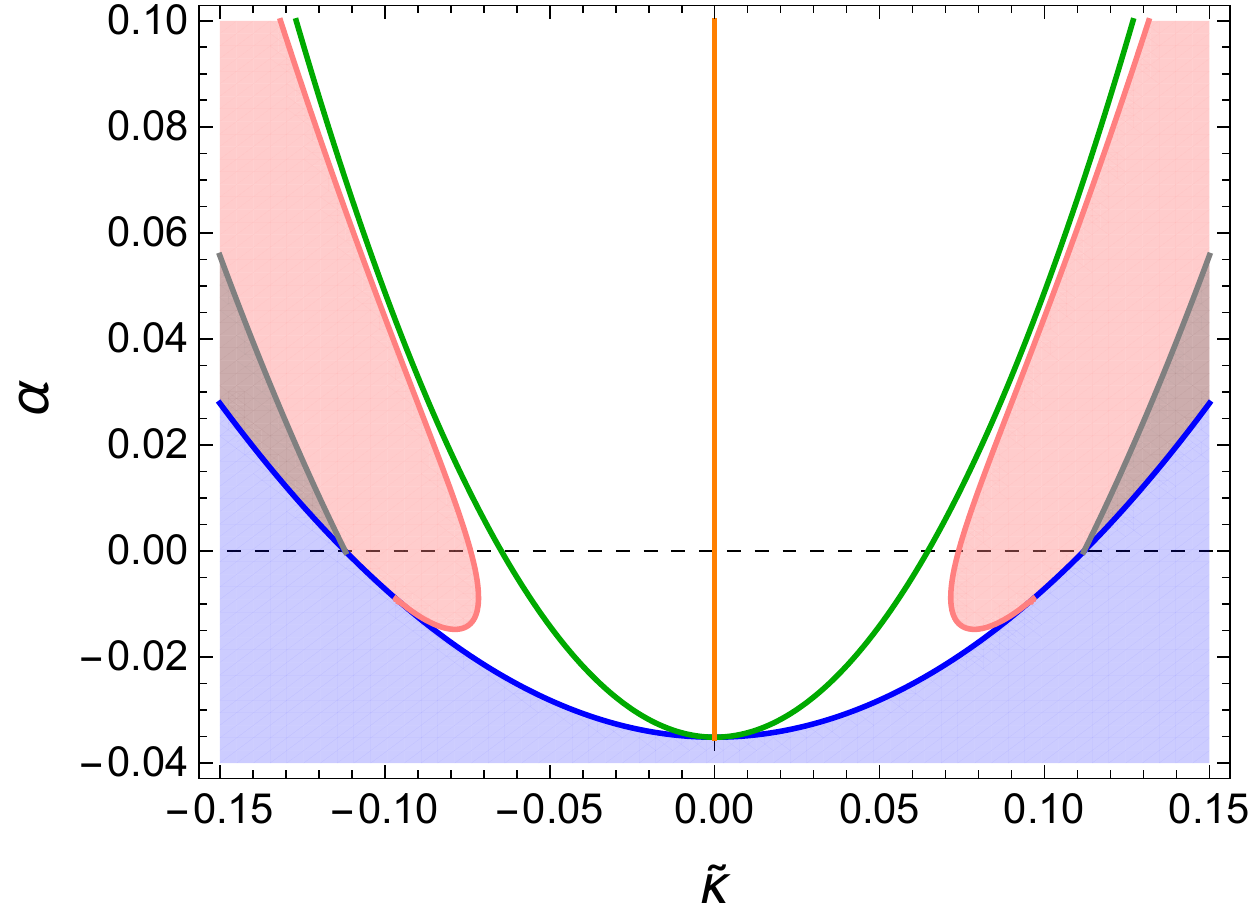}\label{f:egbet0p56ep0p4}}
\caption{Numerical computations based on the analytic leading order formulae of the leading order of instability regions and boundaries of the stripes for \eqref{e:example} in the $(\kap,\alpha)$-plane. Stripes exist in the complement of the blue regions. Bifurcation \eqref{e:egbif} (blue); Eckhaus boundaries \eqref{e:egeh} (green, unstable below); zigzag boundaries (unstable to the left) $\kap=0$ (orange); grey regions: most quasi-square-unstable ($\tel=0$); red regions: hex-unstable; pink regions: most quasi-hex-unstable ($\tel=-\kap/3$); otherwise stable. Here $q=-0.215$. (b) $\beta\in(0,\beta_\tp)$. (c) $\beta=\beta_\tp$ (d) $\beta\in(\beta_\tp,\beta_\ex)$. (e) $\beta>\beta_\ex$.}
\label{f:egalpkap}
\end{figure}

In this system, the Turing conditions are fulfilled with critical wavevectors $(k,\ell)\in S_\kc$ with $\kc = 1$. We have
\begin{align*}
	\hat\calL_0 := -\kcsq D+\A = 
	\begin{pmatrix}
	2 & -1\\
	14 & -7
	\end{pmatrix}.
\end{align*}
From Remark~\ref{r:ev} the rescaled kernel and adjoint kernel eigenvectors of $\hat\calL_0$ and $\hat\calL_0^*$ are given by
\begin{align*}
	\E_0 = -\frac{1}{\sqrt{5}}(1,2)^T,\quad 
	\E_0^* = \frac{1}{\sqrt{5}}(-7,1)^T,
\end{align*}
respectively. Based on the coefficients in \eqref{e:defs}, \eqref{e:stripecoeffalt} we compute, cf.\ Fig.~\ref{f:egalpkap},
\begin{align}
\text{bifurcation curve:}\quad & \alpha = -0.112\beta^2+2.8\kap^2,\label{e:egbif}\\
\text{Eckhaus boundary:}\quad & \alpha=-0.112\beta^2+8.4\kap^2,\label{e:egeh}
\end{align}
and due to the scaling \eqref{e:scaling} the zigzag boundary is at $\kap=0$ to leading order. 
The striped solutions exist for $\alpha>-0.112\beta^2+2.8\kap^2$. Notably, $\alpha = 12.24\calpha$.

\medskip
Fig.~\ref{f:egalpkap} illustrates the stabilities of stripes against the lattice modes. We consider the most unstable quasi-square mode (i.e., $\tel'=0$, cf.~\eqref{e:evrect}), hexagonal mode and the most unstable quasi-hexagonal mode (i.e., $\tel'=-\kap'/3$, cf.~Lemma~\ref{l:rhombevals}). The quadratic coefficient $q$ is linearly dependent on the coefficient $\epsilon$, i.e., $q=-0.537\epsilon$. We choose $\epsilon=0.4$, which leads to $q = -0.215$. The critical value $\beta_\tp \approx 0.378$ (cf.~\eqref{e:betatp}) is such that the quasi-hex-stable region is connected and the hex-unstable band vanishes for $\beta>\beta_\tp$. The critical value $\beta_\ex \approx 0.535$ (cf.~\eqref{e:eckrhombex} with $\theta=1$  most unstable) is such that the quasi-hex-unstable regions are completely covered by the Eckhaus-unstable regions for $\beta>\beta_\ex$. The quasi-hexagonal mode is more unstable than the quasi-square and hexagonal modes.

\begin{figure}[!t]
\centering
\subfloat[$\beta=0$]{\includegraphics[width=0.3\linewidth]{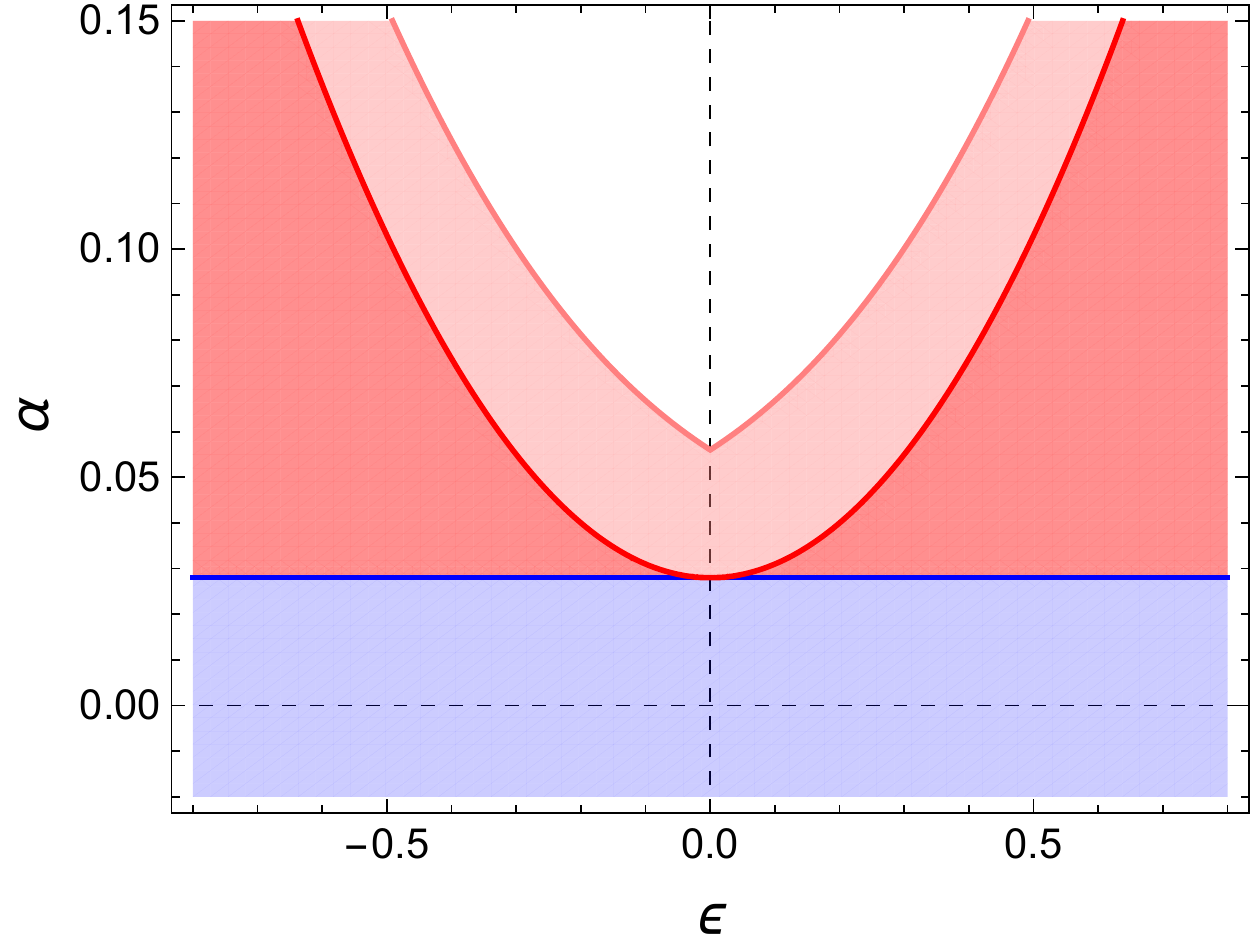}\label{f:egbet0kap0p1}}
\hfil
\subfloat[$\beta=0.3$]{\includegraphics[width=0.3\linewidth]{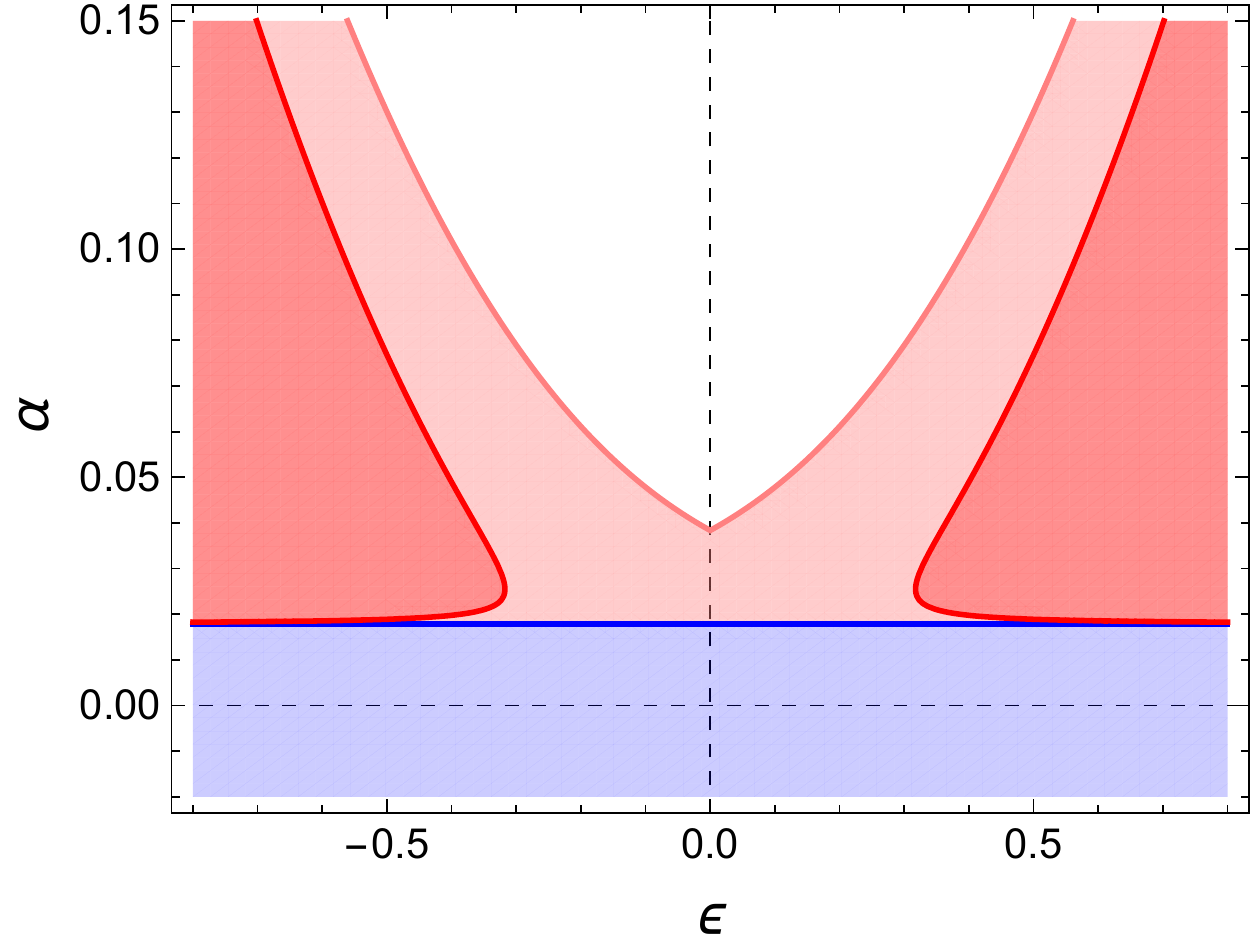}\label{f:egbet0p3kap0p1}}
\hfil
\subfloat[$\beta=0.6$]{\includegraphics[width=0.3\linewidth]{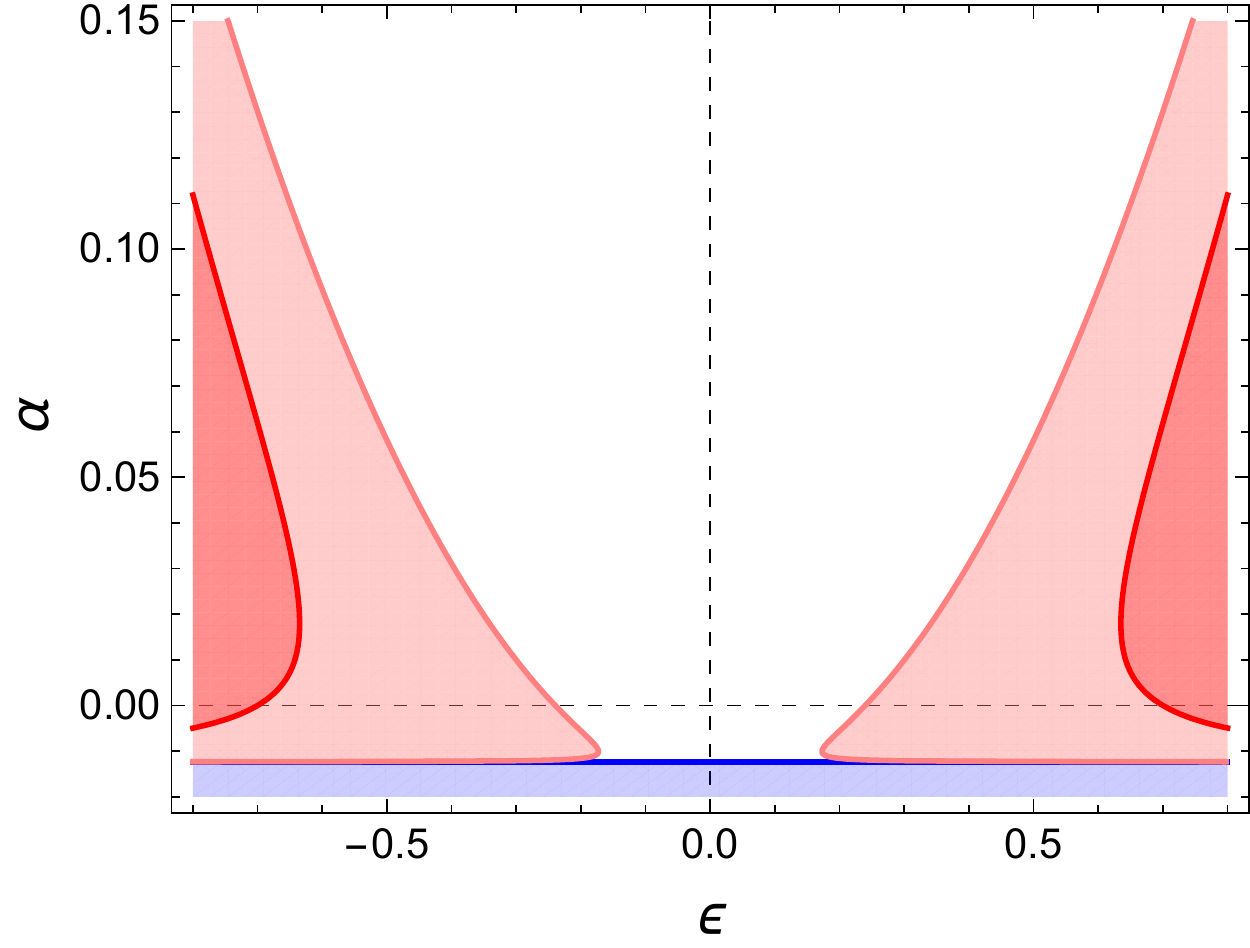}\label{f:egbet0p6kap0p1}}
\caption{Numerical computations  based on the analytic leading order formulae of the leading order of hex- and quasi-hex-instability regions and boundaries of the stripes for \eqref{e:example} in the $(\epsilon,\alpha)$-plane. Colours as in Fig.~\ref{f:egalpkap}. The off-critical parameter $\kap=0.1$. (b) $\beta\in(0,\beta_\ep)$. (c) $\beta>\beta_\ep$.}
\label{f:egalpepsilon}
\end{figure}

\medskip
In Fig.~\ref{f:egalpepsilon} the stability of stripes against quasi-/hexagonal perturbations are depicted. For convenience, and with some abuse of notation we use the coefficient $\epsilon$ as the horizontal axis rather than the quadratic coefficient $q$. The threshold $\beta_\ep\approx0.577$ (cf.~\eqref{e:bifrhombep} with $\theta=1$ most unstable) is such that the quasi-hex-stable region is connected for $\beta>\beta_\ep$. In particular, the hex-stable region is connected for $\beta>0$. The quasi-hexagonal mode is more unstable than the hexagonal one.

\subsection{Numerical example: extended Klausmeier model}\label{s:Klausmeier}

In contrast to the designed example in the previous section, here we consider a model from applications that is not in the normal form \eqref{e:RDS} and where parameters are not in the scaling form \eqref{e:scalingintro} that allows a clean separation of leading order terms. Indeed, the numerical results do not have the symmetry that is induced to leading order by this scaling assumption.
With this in mind, the results presented here can nevertheless be directly explained and related to our analytical results. The study of the Klausmeier model in \cite{Siero2015}, which we refine hereby, was the main motivation for this paper, which now provides a more complete understanding of stripe stability near onset for weak anisotropy.

\medskip
The extended Klausmeier in two space dimensions \cite{Klausmeier1999,Siero2015}, in scaled form, is given by:
\begin{equation}\label{e:Klausmeier}
       \begin{aligned}
       u_t &= d\Delta u + \beta u_x+ a - u - uv^2,\\
       v_t &= \Delta v - mv + uv^2.
       \end{aligned}
\end{equation} where we fix $m=0.45$ and $d=500$.

We complement the computations of the zigzag and sideband stability boundaries near onset from \cite{Yang2019a} for this model by numerical computations of the quasi-square and quasi-hexagon instability boundaries. More generally, we compute the rhomb-breakup boundaries for small advection $\beta$, thereby refining results from \cite{Siero2015}. By rectangle-breakup we denote instability on any rectangular lattice: $\k^\rect_1=(\kappa,0)$; $\k^\rect_2=(0,\ell)$ for any $\ell>0$. So this includes the zigzag instability ($\ell\approx 0$) and the quasi-square lattice ($\ell=\kc$). By rhomb-breakup we denote instability on any rhombic lattice: $\k^\rhb_1=(\kappa,0)$; $\k^\rhb_2=(-\kappa/2,\ell)$; $\k^\rhb_3=-(\kappa/2,\ell)$ for any $\ell>0$. This includes the hexagonal ($\ell=\sqrt{3}\kappa/2$) and quasi-hexagonal ($\ell=\sqrt{\kcsq-\kappa^2/4}$) lattices.

In \cite{Siero2015}, continuation of rhombic and rectangular (in)stability curves was performed by imposing tangency conditions on the spectrum while computing stripe patterns with the software package AUTO \cite{Doedel2007,Rademacher2007}. In \cite{Yang2019a}, a more global brute force approach was used to compute stripe solutions and their spectra in terms of Floquet-Bloch modes, using the software package pde2path \cite{Uecker2014}.

\begin{figure}[t!]
\centering
	{\includegraphics[clip=true, trim=70 0 65 5, width=0.43\linewidth]{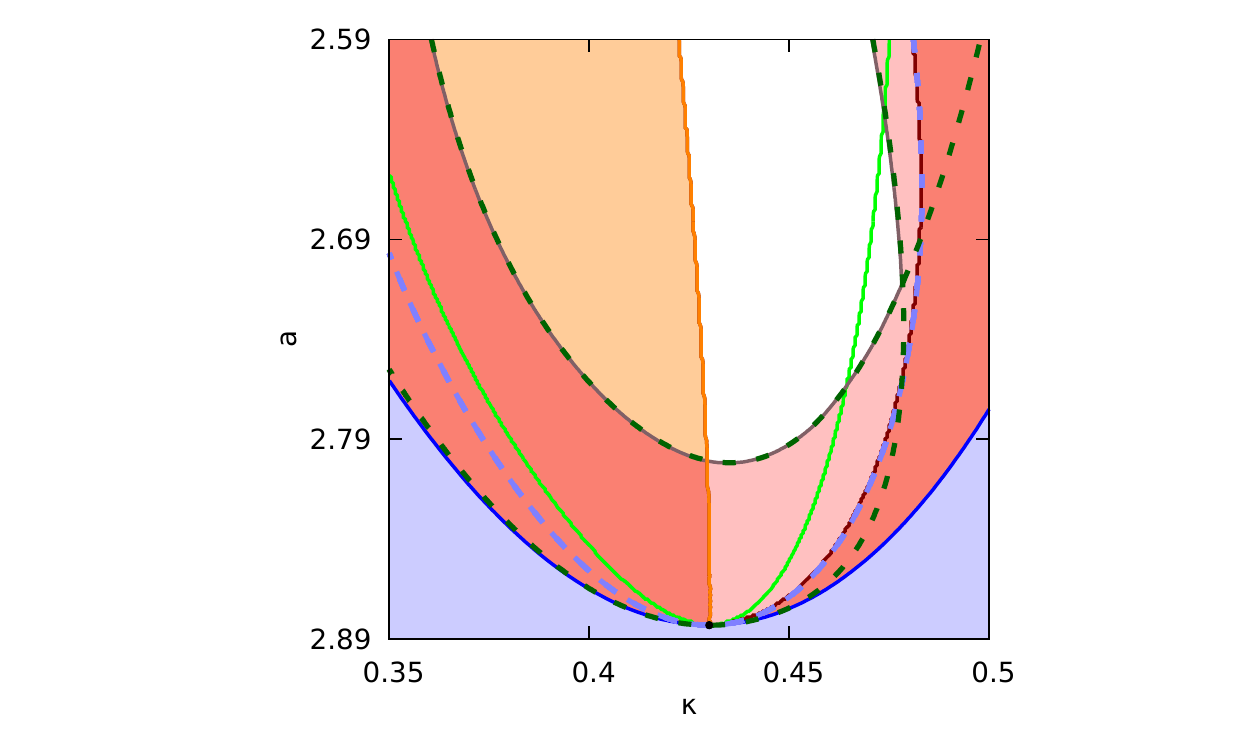}}
	{\includegraphics[clip=true, trim=185 450 170 130, width=0.49\linewidth]{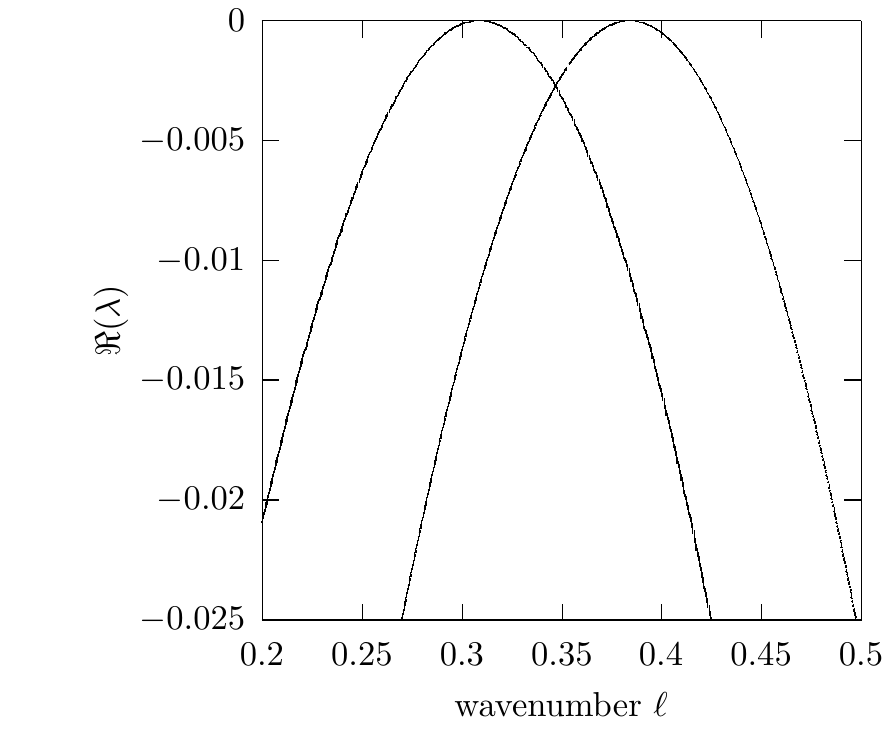}}
	\caption{Left panel: regions of (in)stability of stripes for \eqref{e:Klausmeier} in the $(\kappa,a)$-plane for $\beta=0$. Stripes exist in the complement of the blue region; the light-green curve is the Eckhaus stability boundary. In the orange \& salmon regions, stripes are rectangle unstable, either by zigzag instability (orange curve) or on a quasi-square lattice (red curve). The blue dashed curve - computed by tracing a critical part of the spectrum - perfectly matches the red rectangle instability curve. In the pink \& salmon regions, stripes are rhomb-unstable, so stripes bifurcate rhomb-unstably. The green dashed curves - again computed by tracing a critical part of the spectrum - both partially match the grey rhombic instability curve. For the intersection of the green dashed curves ($\kappa\approx 0.4784$, $a\approx 2.712$), the two corresponding pieces of tangential spectrum of the stripe pattern are shown in the right panel.}
	\label{f:Klausmeierflat}
\end{figure}

\medskip
Fig. \ref{f:Klausmeierflat} shows the various stripe stability regions in the extended Klausmeier model for $\beta=0$. The coloured regions, bounded by continuous curves, are the result of the brute force method on an equidistant grid (spacing between neighbouring points $a=0.001$ and $\kappa=0.0005$). The dashed curves are the results of imposing a tangency condition on the spectrum using AUTO \cite{Doedel2007,Rademacher2007}. The computation of tangent spectrum is computationally more efficient, but only establishes (in)stability of a local piece of spectrum. In the left panel there are two green dashed curves corresponding to rhomb-instability. Where they cross, there are indeed two critical curves of spectrum (right panel). Focusing on only one of these yields an incomplete stability picture \cite{Siero2015}.

\begin{figure}[t!]
\centering
	\subfloat[$\beta=0$]{\includegraphics[clip=true, trim=70 0 65 5, width=0.33\linewidth]{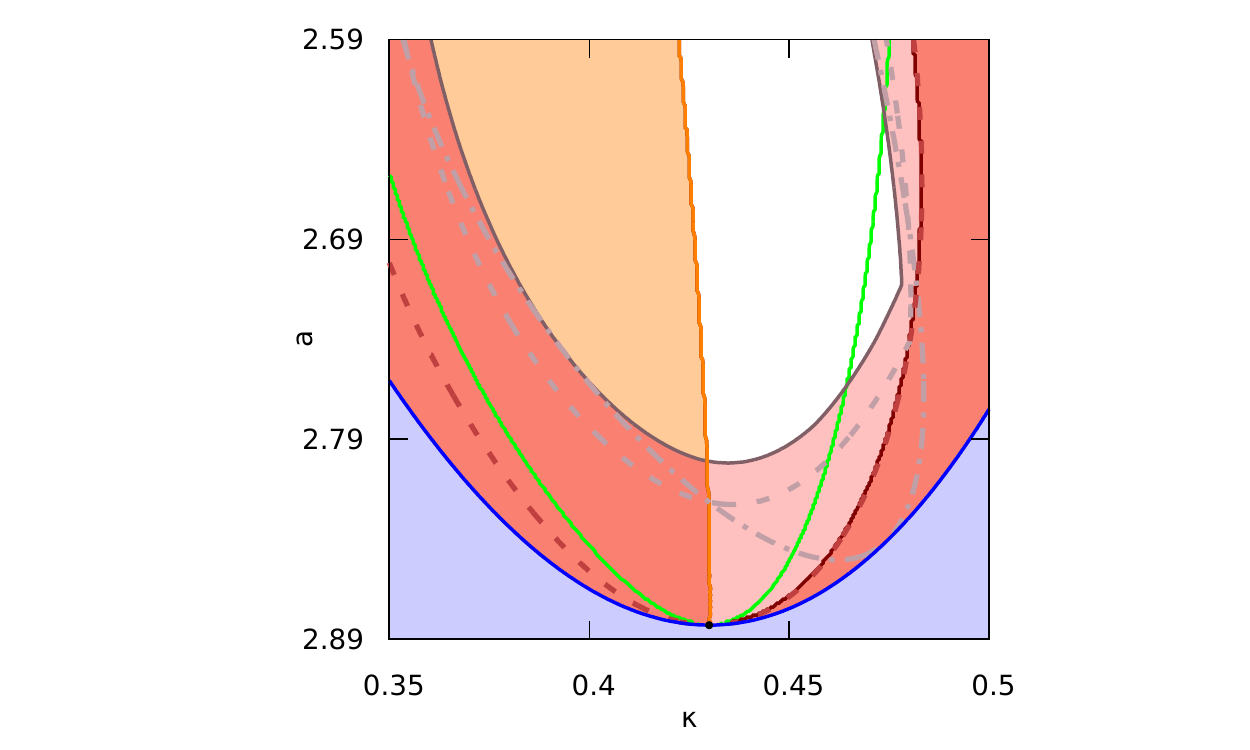}\label{f:Klausmeierbet0}}
	\subfloat[$\beta=30$]{\includegraphics[clip=true, trim=70 0 65 5, width=0.33\linewidth]{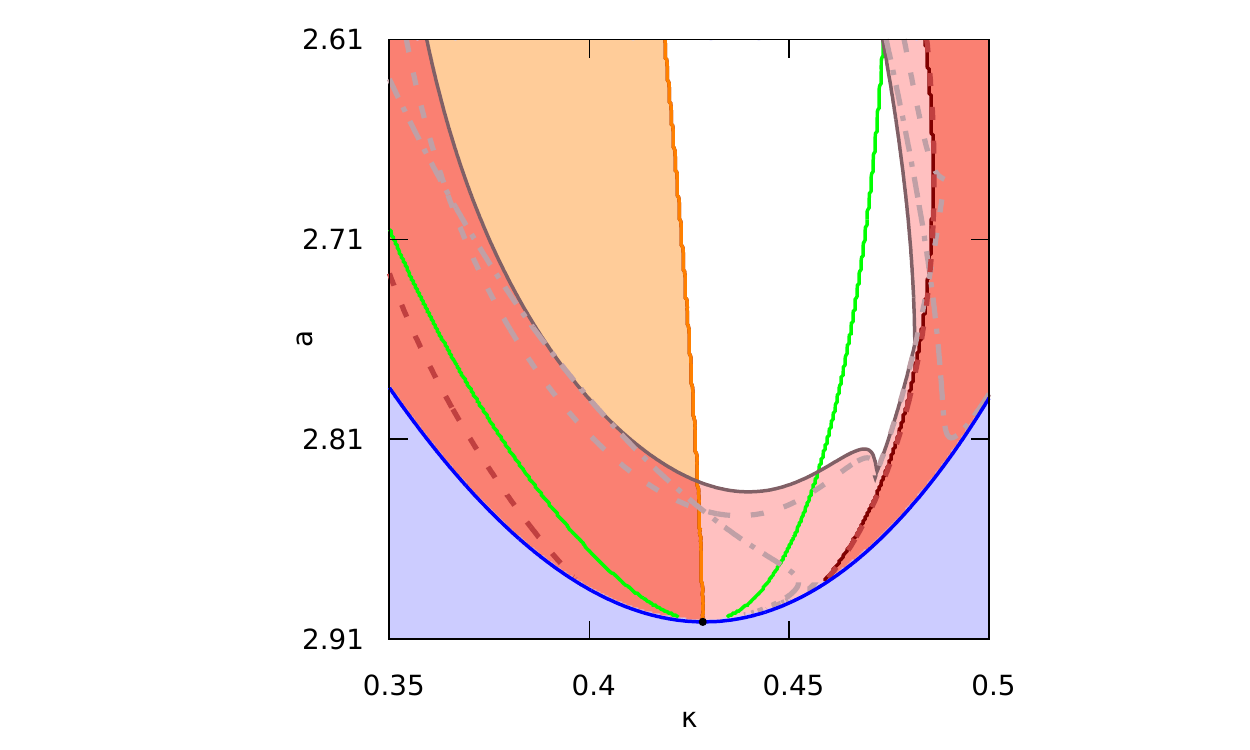}}
	\subfloat[$\beta=30$]{\includegraphics[clip=true, trim=70 0 65 5, width=0.33\linewidth]{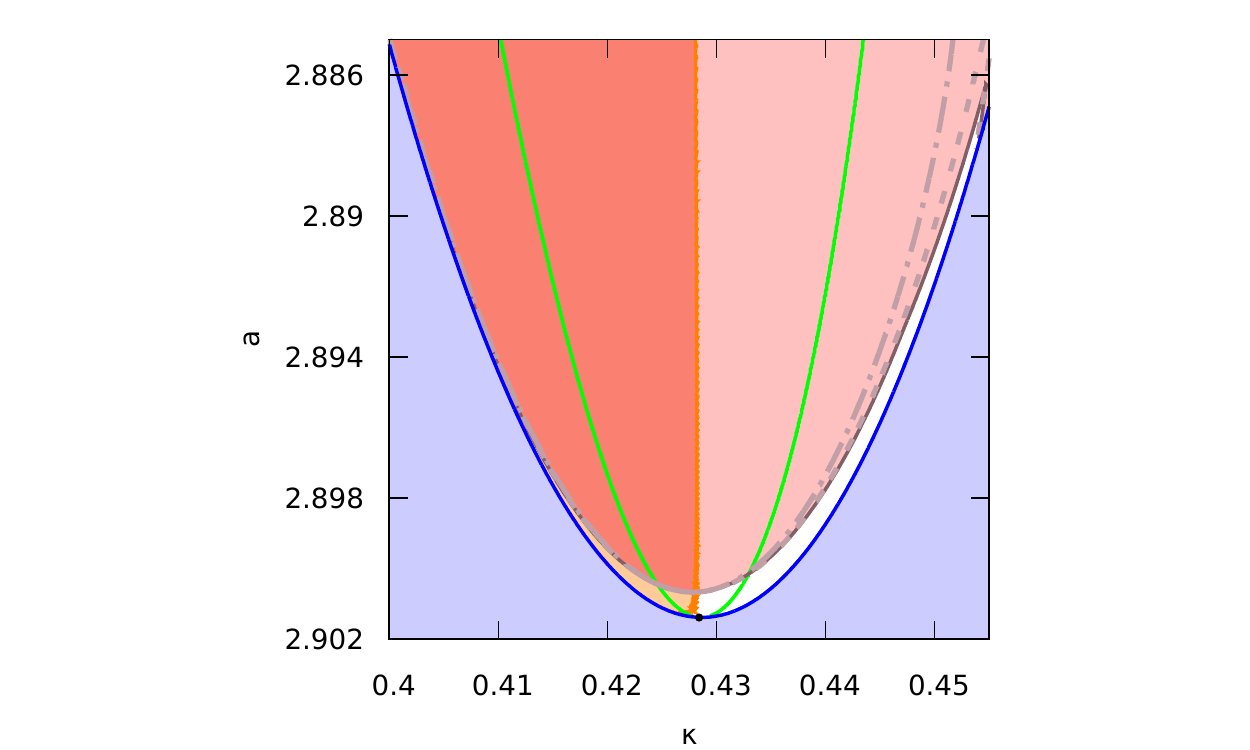}} \\
	\subfloat[$\beta=40$]{\includegraphics[clip=true, trim=70 0 65 5, width=0.33\linewidth]{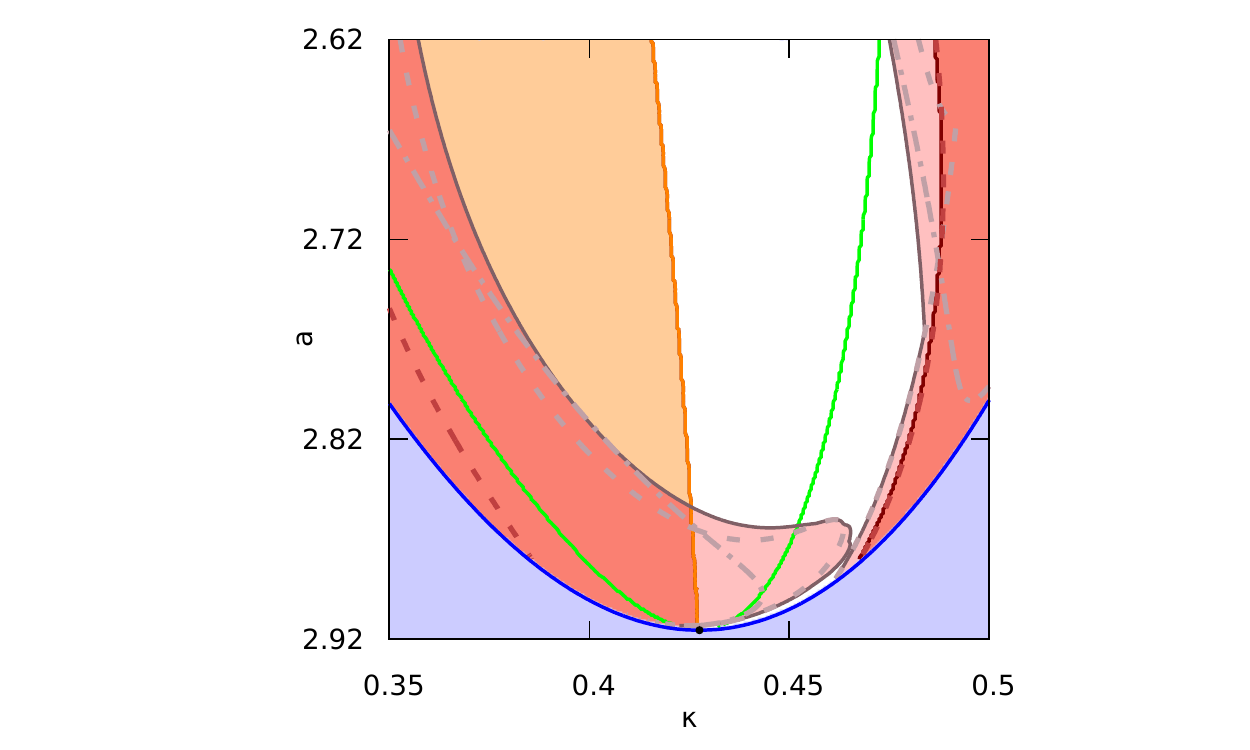}}
	\subfloat[$\beta=50$]{\includegraphics[clip=true, trim=70 0 65 5, width=0.33\linewidth]{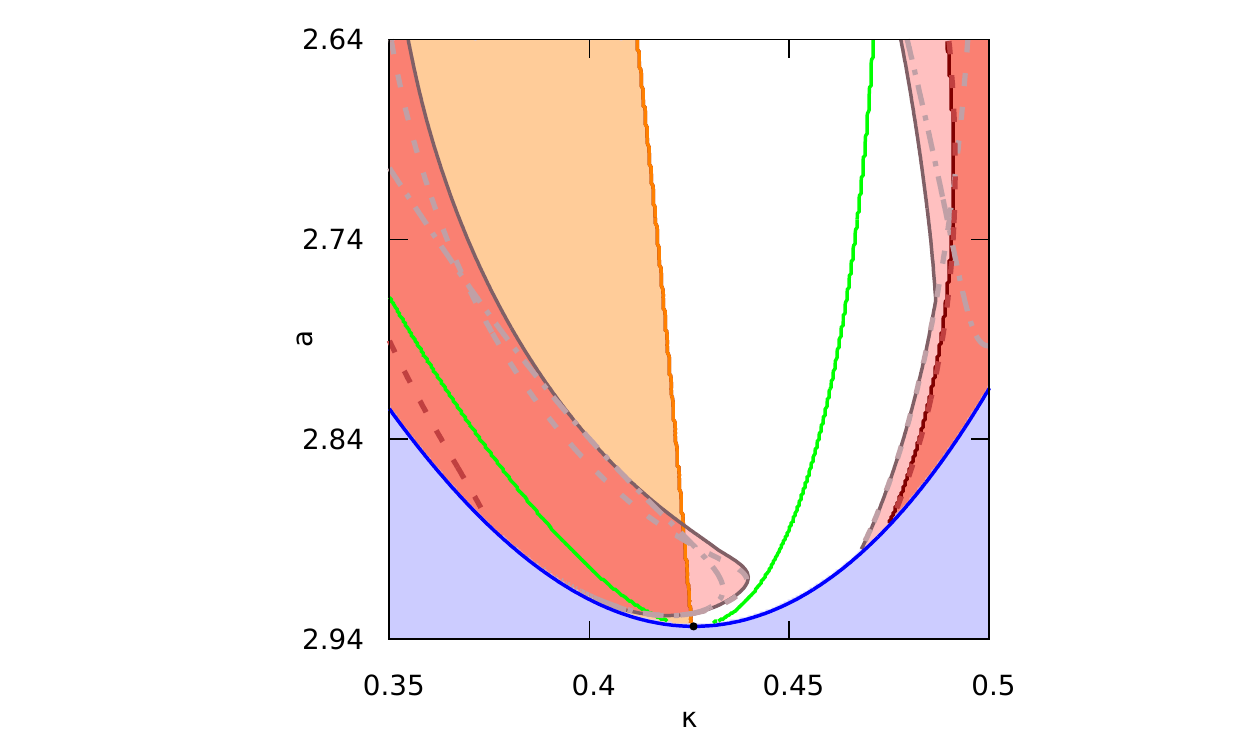}}
	\subfloat[$\beta=100$]{\includegraphics[clip=true, trim=70 0 65 5, width=0.33\linewidth]{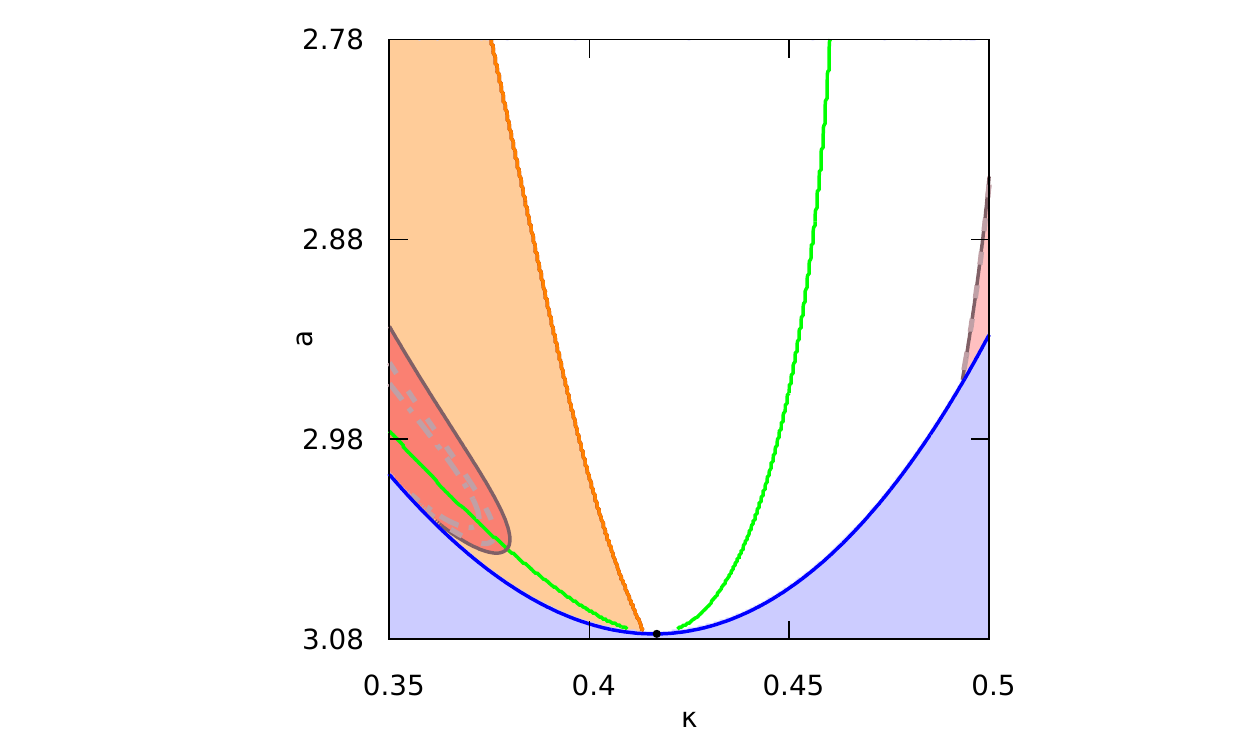}}

	\caption{Regions of (in)stability of stripes for \eqref{e:Klausmeier} in the $(\kappa,a)$-plane for $\beta=0$. Stripes exist in the complement of the blue region; the light-green curve is the Eckhaus stability boundary. In the orange \& salmon regions, stripes are rectangle unstable, either by zigzag instability (orange curve) or on a quasi-square lattice (red curve). The light-red dashed quasi-square lattice instability curve perfectly matches the red rectangle instability curve on the right. In the pink \& salmon regions, stripes are rhomb-unstable. The grey dashed quasi-hexagonal lattice curves are close to the rhombic instability curves near pattern formation onset, but they diverge further away. The dash-dotted grey curve is the instability curve for the hexagonal lattice. (a) $\beta=0$. (b-c) $\beta=30$. (d) $\beta=40$. (e) $\beta=50$. (f) $\beta=100$. For the zoom in (c) of (b), a finer grid was used (with a spacing between neighbouring points of $a=0.00002$ and $\kappa=0.0001$).
		}\label{f:Klausmeier}
\end{figure}

\medskip
In Fig. \ref{f:Klausmeier} the impact of advection on rhombic and rectangle stability of stripes is depicted. Panel (a) shows the same situation $\beta=0$ as in Fig. \ref{f:Klausmeierflat}. For relatively small advection $\beta=30$, the stability regions have deformed a bit (panel b) and, attached to the Turing-Hopf, a separate small region of stable stripes has appeared (panel c). At $\beta=40$, the pink region of rhomb-unstable stripes has just split into two separate regions and a stable connection for the lattice modes has `opened'. For $\beta=50$, the white regions are also connected above the green sideband curve, implying that the two regions of stable stripes for small amplitude and larger amplitude are connected also in terms of large-wavelength modes. At the largest value, $\beta=100$, only the large-wavelength instabilities bound stripe stability near onset and stripes bifurcating with critical wavenumber remain stable after onset.

\appendix

\section{Proof of Theorem~\ref{t:cmfrect}}\label{s:rectmat}
We recall the simplified linearisation \eqref{e:genlin1D}, namely
\begin{align}\label{e:linsimp}
\partial_u f(u_c;\mu) =\ &\scal^2 P\Big(\Lc(\mu_2)+\Lc(\mu_1)\Psi_{11}[\mu_1,\cdot] + 2A'^2\Q[\Psi_{20}[u_1,u_1],\cdot]\\ 
& + 4A'^2\Q[u_1,\Psi_{20}[u_1,\cdot]] + 3A'^2\K[u_1,u_1,\cdot] \Big) + \calO(\scal^3).
\end{align}

The matrix $L_1$ is known a priori from Theorem~\ref{t:cmfstripe}, but for completeness, we derive it here directly.
Setting $u_1=0$ gives the linearisation in the zero state so that the first two terms in the bracket generate the eigenvalue from \eqref{e:evlinearop}, $\alpha'+\rho_\beta\beta'^2+\rho_\kap\kap'^2$, which is of the form $-\rho_\nl A'^2 + \calO(\scal)$, cf.\ \eqref{e:stripeeqn}.

More specifically these contribute the diagonal 2-by-2 matrix $A_L:=-\rho_\nl\Id$ to the linearisation at order $\scal^2$. As to the nonlinear terms, the simplest is $K[u_1,u_1,\cdot]$ and with $u_1=(\rme^{\rmi x}+\rme^{-\rmi x})E_0$ we find  the 2-by-2 matrix with entries generated by choosing $\sigma_1,\sigma_2\in\{\pm1\}$ as
\[
\langle K[u_1,u_1,\rme^{\sigma_1 \rmi x}E_0], \rme^{\sigma_2 \rmi x}E_0^*\rangle
=\frac{\hK}{|\Omega_1|} \int_{\Omega_1}(\rme^{2\rmi x}+2+\rme^{-2\rmi x})\rme^{\rmi (\sigma_1-\sigma_2) x}\rmd x.
\]
This results in the matrix $A_K:=\hK\begin{pmatrix} 2 & 1\\1&2 \end{pmatrix}$.
The contributions from the quadratic term depend on $\Psi_{20}$, which can be computed from the general centre manifold characteristic equation \cite{Haragus2010}
\[
\partial_u\Psi(u_c;\mu)f(u_c;\mu) = P_h (\calL_\mu (u_c+ \Psi(u_c;\mu)) + F(u_c+\Psi(u_c;\mu))),
\]
which holds for all $u_c\in N$, $|u_c|\ll1$. At the bifurcation point, i.e., $\dot u_c = f(u_c;\mu)=0$, the above equation reduces to the fixed point equation~\cite[Eq. A.9]{Yang2019a}
\[
P_h\calL_\mu\Psi(u_c;\mu) = -P_h F(u_c+\Psi(u_c;\mu)) - P_h (\calL_\mu-\calL_0) u_c.
\]
At order $u_c^2$ we find $P_h(\calL_0\Psi_{20} + 2Q)=0$ on $N$ in analogy to the expansion for~\cite[Eq. A.8]{Yang2019a}, so that $\Psi_{20} = -2\calL_0^{-1}Q$. This means
\begin{align*}
\Psi_{20}[u_1,u_1] &= \Psi_{20}[E_0,E_0](\rme^{2\rmi x}+2+\rme^{-2\rmi x}) = \tQ+\frac 1 2 \hQ(\rme^{2\rmi x}+\rme^{-2\rmi x}),\\
\Psi_{20}[u_1,ae_0+b\overline{e_0}] &= \Psi_{20}[E_0,E_0](a\rme^{2\rmi x}+a+b+b\rme^{-2\rmi x})\\
& = \frac a 2\left(\tQ+ \hQ\rme^{2\rmi x}\right)+\frac b 2\left(\tQ+\hQ\rme^{-2\rmi x}\right).
\end{align*}
The first equation is in fact an immediate consequence of the formula for stripes and $f(u_c;\mu)=0$.
As to the matrix entries this generates, we compute for the first row
\begin{align*}
\langle Q[\Psi_{20}[u_1,u_1],e_0],e_0^* \rangle 
& = \langle Q[\tQ+ \frac 1 2\hQ\rme^{2\rmi x},e_0],e_0^* \rangle
= \tq\\
\langle Q[\Psi_{20}[u_1,u_1],\overline{e_0}],e_0^* \rangle 
& = \langle Q[ \tQ+ \frac 1 2\hQ\rme^{2\rmi x},\overline{e_0}],e_0^* \rangle
=\frac 1 2 \hq,
\end{align*}
whose entries are reversed in the second row so we get $A_Q:=\frac 1 2 \begin{pmatrix} 2\tq  & \hq \\ \hq & 2\tq \end{pmatrix}$. Analogously,
\begin{align*}
\langle Q[u_1,\Psi_{20}[u_1,e_0]],e_0^* \rangle 
& = \langle Q[e_0+\overline{e_0},\frac 1 2\left(\tQ+ \hQ\rme^{2\rmi x}\right)],e_0^* \rangle
=\frac 1 2 (\tq+\hq)\\
\langle Q[u_1,\Psi_{20}[u_1,\overline{e_0}]],e_0^* \rangle 
& = \langle Q[e_0+\overline{e_0},\frac 1 2\left(\tQ+ \hQ\rme^{-2\rmi x}\right)],e_0^* \rangle
=\frac 1 2 \tq,
\end{align*}
whose entries are reversed in the second row so we get $B_Q:=\frac 1 2 \begin{pmatrix} \tq+\hq  & \tq \\ \tq & \tq+\hq \end{pmatrix}$.

In sum, the matrix for the linearisation on the centre manifold is, as claimed,
\[
\partial_u f(u_c;\mu) =\scal^2A'^2( A_L+3A_K+2A_Q+4B_Q)=A^2\rho_\nl\begin{pmatrix} 1  & 1 \\ 1 & 1\end{pmatrix}.
\]

\medskip
The claimed block diagonal structure for the linearisation in stripes is a result of non-resonance between the arising wave vectors; the only relevant resonances away from the subblock $L_1$ are $\krec_2+\krec_{-2}=0$. Casting $L_\sq$ as a matrix, its entries are 
\[
(L_\sq)_{j,\ell} = \langle \partial_u f(u_c;\mu)e_\ell, e_j^*\rangle,
\;j,\ell=\pm1,\pm2.
\]
Being the linearisation in stripes, multiples of $\krec_{\pm 1}$ enter from $\partial_u f(u_c;\mu)e_{\pm1}$, but (in the chosen ordering) off-diagonal entries give one additional wavevector $\krec_j$ for $j\neq \pm1$, and hence no resonance is possible.  Therefore, the linearisation has block-diagonal form.

\medskip
Concerning the subblock $L_2^\sq$, analogous to $L_1$, due to the lack of resonances, \eqref{e:genlin} simplifies to \eqref{e:linsimp}. Setting $u_1=0$ gives the linearisation in the trivial equilibrium to order $\scal^2$ and 
the eigenvalues arise directly from the Fourier transform in $y$-direction, or by using \eqref{e:evlinearop} with $\beta=0$, $\kap=\tel$, i.e., 
\[
\lambda_\tel = \alpha + \rho_\kap\tel^2 +\calO(\scal^3)\in\R,
\]
so that with \eqref{e:scaling} we have $\lambda_\tel = \scal^2\lambda_\tel'$. Concerning the simplest nonlinear term $K[u_1,u_1,\cdot]$. The stripes $u_1=(\rme^{\rmi \krec_1\cdot \x}+\rme^{-\rmi \krec_1\cdot \x})E_0$ yield
\[
\langle K[u_1,u_1,e_\ell],e_j^*\rangle
=\frac{\hK}{|\Omega_2|} \int_{\Omega_2} (\rme^{\rmi 2\krec_1\cdot\x}+2+\rme^{-\rmi 2\krec_1\cdot\x})\rme^{\rmi (\krec_\ell- \krec_j)\cdot \x}\rmd \x
\] 
which, for $j,\ell\neq \pm 1$, gives a contribution on the diagonal $j=\ell$ only, namely 
$6\hK \scal^2A'^2$.

It remains to consider the contributions from $Q$ and the centre manifold $\Psi_{20}$, i.e. the terms from  Corollary~\ref{c:cmf} at order $\scal^2$: 
\begin{gather*}
	2PQ[\Psi_{20}[u_1,u_1],\cdot], \ 4PQ[u_1,\Psi_{20}[u_1,\cdot]].
\end{gather*}
Since $\Psi_{20}[u_1,u_1]=\frac{1}{2}\hQ (\rme^{2\rmi\krec_1 x}+\rme^{-2\rmi\krec_1 x})+\tQ$, for $\ell,j\neq \pm1$ we find
\[
2\langle Q[\Psi_{20}[u_1,u_1],e_\ell],e_j^*\rangle = 2\tq 
\]
for $j=\ell$ and zero otherwise due to non-resonance with $2 \krec_1$. 

As to $4PQ[u_1,\Psi_{20}[u_1,\cdot]]$ we first compute, since $\ell,j\in\{\pm2\}$ and the only contribution comes from the identity in $P_h = \Id-P$ that 
\begin{align*}
\Psi_{20}[u_1, e_\ell] &= -\calL_0^{-1}P_h Q[u_1,e_\ell] = -\calL_0^{-1}(\Q[u_1,e_\ell] - P\Q[u_1,e_\ell]) \\
&= -\calL_0^{-1} Q[\E_0,\E_0]\rme^{\rmi (\krec_1+\krec_\ell)\cdot\x}\\
&=-(-2\kcsq D+L)^{-1}Q[\E_0,\E_0]\rme^{\rmi (\krec_1+\krec_\ell)\cdot\x}.
\end{align*}
Substituting into $4PQ[u_1,\Psi_{20}[u_1,\cdot]]$ gives for $\ell=j$
\begin{align*}
4\langle Q[u_1,\Psi_{20}[u_1,e_\ell]],e_j^*\rangle &= 8\langle Q[E_0,\Q_{11}],E_0^*\rangle = 8q_{11},
\end{align*}
and zero otherwise.

\section{Proof of Theorem~\ref{t:cmfhex}}\label{s:hexmat}

The 1D subsystem is clearly an invariant subsystem (as are several others) and the form of the block $L_1$ follows from Theorem~\ref{t:cmfstripe}. 

Analogous to Appendix~\ref{s:rectmat}, the claimed block diagonal structure for the linearisation in stripes is a result of non-resonance between the arising wave vectors; the only relevant resonances away from the subblock $L_1$ are triads $\k_1+\k_2=\k_{-3}$. Casting $L_\hex$ as a matrix, its entries are 
\[
(L_\hex)_{j,\ell} = \langle \partial_u f(u_c;\mu)e_\ell, e_j^*\rangle,
\;j,\ell=\pm1,\pm2,\pm3.
\]
Being the linearisation in stripes, multiples of $k_{\pm 1}$ enter from $\partial_u f(u_c;\mu)e_{\pm1}$, but (in the chosen ordering) off-diagonal entries give one additional wavevector $k_j$ for $j\neq \pm1$, and hence no triad is possible.  Therefore, the linearisation has block-diagonal form.

The two equal subblocks $L_2^\hex$ are obtained by symmetry, and Corollary~\ref{c:cmf} gives the relevant terms at order $\scal$ and $\scal^2$. The only term at order $\scal$ is $2\scal A'PQ[u_1,\cdot]$, which contributes through triads on the off-diagonal only as $2\scal A'q$.

Setting $u_1=0$ gives the linearisation in the trivial equilibrium to order $\scal^2$ and 
the eigenvalues are known a priori from Lemma~\ref{l:Turbeta}, see also \eqref{e:stripeeqn2}, and $\alpha, \kap$ can be readily included analogous to \eqref{e:evlinearop}; note that the coefficient of $\kap$ stems from isotropic domain scaling.
By choice of $c$ and with $k_j$ the first component of the wavevectors, these eigenvalues have the form 
\[
\lambda_{\mu,j} = \alpha + k_j^2\rho_\beta\beta^2 + \rho_\kap\kap^2 + \calO(|A|^3)\in\R,
\]
which are all equal for $j\neq \pm 1$ ($k_2=k_3=1/2$) and enter as entries of $L_\hex$ along the diagonal. Due to the scalings \eqref{e:scaling}, we can write $\lambda_{\mu,j} = \scal^2\lambda_{\mu,j}'$.

We proceed analogous to Appendix~\ref{s:rectmat} with the simplest nonlinear term $K[u_1,u_1,\cdot]$. Stripes $u_1=(\rme^{\rmi \k_1\cdot \x}+\rme^{-\rmi \k_1\cdot \x})E_0$ yield
\[
\langle K[u_1,u_1,e_\ell],e_j^*\rangle
=\frac{\hK}{|\Omega_3|} \int_{\Omega_3} (\rme^{\rmi 2\k_1\cdot\x}+2+\rme^{-\rmi 2\k_1\cdot\x})\rme^{\rmi (\k_\ell- \k_j)\cdot \x}\rmd \x
\] 
which, for $j,\ell\neq \pm 1$, gives a contribution on the diagonal $j=\ell$ only, namely 
$6\hK \scal^2A'^2$.

\medskip
It remains to consider the contributions from $Q$ and the centre manifold via $\Psi_{20}, \Psi_{11}$, i.e. the five terms from  Corollary~\ref{c:cmf} at order $\scal^2$: 
\begin{gather*}
	2PQ[\Psi_{20}[u_1,u_1],\cdot], \ 4PQ[u_1,\Psi_{20}[u_1,\cdot]], \ 2P\Lc(\mu_1)\Psi_{20}[u_1,\cdot], \\
	2PQ[\Psi_{11}[\mu_1,u_1],\cdot], \ 2P\Q[u_1,\Psi_{11}[\mu_1,\cdot]].
\end{gather*}
Notably, the first two enter with a factor $A'^2$, while the others only have a factor $A'$.

Since $\Psi_{20}[u_1,u_1]=\frac{1}{2}\hQ (\rme^{2\rmi\k_1 x}+\rme^{-2\rmi\k_1 x})+\tQ$, for $\ell,j\neq \pm1$ we find
\[
2\langle Q[\Psi_{20}[u_1,u_1],e_\ell],e_j^*\rangle = 2\tq 
\]
for $j=\ell$ and zero otherwise due to non-resonance with $2 \k_1$.

As to $4PQ[u_1,\Psi_{20}[u_1,\cdot]]$ we first compute, since $\ell,j\in\{\pm2,\pm3\}$ and the only contribution comes from a triad $\k_1+\k_\ell=\k_{-j}$ that 
\begin{align*}
\Psi_{20}[u_1, e_\ell] &= -\calL_0^{-1}P_h Q[u_1,e_\ell] \\
&= -\calL_0^{-1} (Q[\E_0,\E_0]\rme^{\rmi (\k_1+\k_\ell)\cdot\x} - \langle Q[\E_0,\E_0],E_0^*\rangle E_0\rme^{\rmi \k_j\cdot \x})\\
&=-(-\kcsq D+L)^{-1}(Q[\E_0,\E_0]- \langle Q[\E_0,\E_0],E_0^*\rangle E_0)\rme^{\rmi \k_j\cdot \x} = \Q_1\rme^{\rmi \k_j\cdot \x},
\end{align*}
with $\Q_1$ as in the theorem statement. Substitution into $4PQ[u_1,\Psi_{20}[u_1,\cdot]]$ gives 
\begin{align*}
4\langle Q[u_1,\Psi_{20}[u_1,e_\ell]],e_j^*\rangle &= 8\langle Q[E_0,\Q_1],E_0^*\rangle = 8q_1,
\end{align*}
for $\ell=j$, and zero otherwise.\par

As to the third term, triads $\k_1+\k_\ell=\k_{-j}$ give the only non-trivial term
\begin{align*}
2\langle L(\mu_1)\Psi_{20}[u_1,e_\ell],e_j^*\rangle = 2\langle(-2\kap'\kc D+\rmi \beta'\kc k_\ell B )\Q_1,E_0^*\rangle,
\end{align*}
and its complex conjugate on the anti-diagonal of $L_2$.

For the fourth term $2PQ[\Psi_{11}[\mu_1,u_1],\cdot]$, the characteristic equation of the centre manifold to order $u\mu$ gives $\calL_0\Psi_{11} =-P_h \partial_\mu L(0)$, which means 
\begin{align*}
\Psi_{11}[\mu_1,u_1] &= -\calL_0^{-1} (\rmi \beta'\kc B u_1 -2\kap'\kc Du_1) \\
& = \rmi\beta'w_{A\beta}(\rme^{\rmi x} - \rme^{-\rmi x})+\kap'w_{A\kap}(\rme^{\rmi x} + \rme^{-\rmi x})
\end{align*}
note $PBu_1=0$ by choice of $c$ and $PDu_1=0$ as remarked earlier. Therefore, only triads $\k_1+\k_\ell=\k_{-j}$ give the nontrivial term
\begin{align*}
2\langle Q[\Psi_{11}[\mu_1,u_1],e_\ell],e_j^*\rangle = 
2\langle \Q[\rmi\beta'w_{A\beta}+\kap'w_{A\kap},\E_0],\E_0^* \rangle.
\end{align*}

The final quadratic term is $2P\Q[u_1,\Psi_{11}[\mu_1,\cdot]]$. Here, the triads $\k_1+\k_\ell = \k_{-j}$ give the nontrivial term
\[
2\langle \Q[u_1,\Psi_{11}[\mu_1,e_\ell]], e_j^* \rangle = 2\langle \Q[\E_0,\rmi\beta'k_\ell w_{A\beta}+ \kap' w_{A\kap}],\E_0^* \rangle.
\]

Together with the previous two terms, the anti-diagonal terms generate $p(\mu_1)$ and its complex conjugated, i.e. the matrix 
$\begin{pmatrix}0&p(\mu_1)\\\overline{p(\mu_1)}&0\end{pmatrix}$.

\end{document}